\documentclass[11pt]{amsart} \usepackage{amscd}
\usepackage{amsfonts,amssymb,latexsym} 
\input xy
\xyoption{all}
\xyoption{arc}
\xyoption{ps}
\xyoption{dvips}
\CompileMatrices

\newcommand{\tensor}{tensor}
\newcommand{\tensors}{tensors}

\newcommand{\SB}{{\mathcal{B}}}

\newcommand{\SE}{{\mathcal{E}}}
\newcommand{\SF}{{\mathcal{F}}}
\newcommand{\SG}{{\mathcal{G}}}
\newcommand{\SH}{{\mathcal{H}}}
\newcommand{\SI}{{\mathcal{I}}}

\newcommand{\SK}{{\mathcal{K}}}

\newcommand{\SM}{{\mathcal{M}}}

\newcommand{\SO}{{\mathcal{O}}}
\newcommand{\SP}{{\mathcal{P}}}

\newcommand{\SU}{{\mathcal{U}}}

\newcommand{\SW}{{\mathcal{W}}}

\newcommand{\SZ}{{\mathcal{Z}}}

\newcommand{\PP}{\mathbb{P}}
\newcommand{\ZZ}{\mathbb{Z}}
\newcommand{\CC}{\mathbb{C}}
\newcommand{\RR}{\mathbb{R}}

\newcommand{\DD}{\mathbb{D}}

\newcommand{\isom}{\cong}

\newcommand{\Ext}{\operatorname{Ext}}

\newcommand{\Spec}{\operatorname{Spec}}

\newcommand{\codim}{\operatorname{codim}}

\newcommand{\Hom}{\operatorname{Hom}}
\newcommand{\SHom}{\mathcal{H}{\rm om}}
\newcommand{\Mor}{\operatorname{Mor}}
\newcommand{\Quot}{\operatorname{Quot}}
\newcommand{\Pic}{\operatorname{Pic}}
\newcommand{\Sch}{\operatorname{Sch}}
\newcommand{\Sets}{\operatorname{Sets}}

\newcommand{\id}{\operatorname{id}}

\newcommand{\surj}{\twoheadrightarrow}
\newcommand{\inj}{\hookrightarrow}
\newcommand{\too}{\longrightarrow}
\newcommand{\grad}{\rm grad \,}
\newcommand{\gr}{{\rm gr\,}}
\newcommand{\rk}{\operatorname{rk}}
\newcommand{\End}{\operatorname{End}}
\newcommand{\wt}{\widetilde}
\newcommand{\tr}{\operatorname{tr}}
\newcommand{\gl}{\operatorname{GL}}
\newcommand{\pgl}{\operatorname{PGL}}
\newcommand{\glv}{{\operatorname{GL}(V)}}
\newcommand{\pglv}{\operatorname{PGL}(V)}
\newcommand{\slv}{\operatorname{SL}(V)}
\newcommand{\slw}{\operatorname{SL}(W)}
\newcommand{\sltwo}{\mathfrak{sl}_2}
\newcommand{\glw}{\operatorname{GL}(W)}
\newcommand{\GLR}{\operatorname{GL}(R)}

\newcommand{\glrc}{\operatorname{GL}(R,\CC)}

\newcommand{\tens}{{\rm t}}
\newcommand{\ev}{{\rm ev}}
\newcommand{\gitq}{/\!\!/}
\newcommand{\an}{{\rm an}}
\newcommand{\et}{{\rm et}}

\newcommand{\fg}{\mathfrak{g}}
\newcommand{\fq}{\mathfrak{q}}
\newcommand{\fh}{\mathfrak{h}}
\newcommand{\ft}{\mathfrak{t}}
\newcommand{\fl}{\mathfrak{l}}
\newcommand{\fu}{\mathfrak{u}}
\newcommand{\fgp}{{\mathfrak{\fg}'}}
\newcommand{\fz}{\mathfrak{z}}
\newcommand{\fc}{\mathfrak{c}}
\newcommand{\fa}{\mathfrak{a}}
\newcommand{\autlie}{{\rm Aut}}
\newcommand{\Aut}{{\rm Aut}}
\newcommand{\autfgp}{{\rm Aut}(\fgp)}
\newcommand{\Lie}{{\rm Lie}}
\newcommand{\sk}{{\rm skew}}
\newcommand{\asym}{{\rm asym}}
\newcommand{\ogp}{{\rm O}(\fgp)}
\newcommand{\glgp}{{\rm GL}(\fgp)}
\newcommand{\sogp}{{\rm SO}(\fgp)}
\newcommand{\sone}{{\mathbb{S}^1}}
\newcommand{\groupp}{{G_2}}
\newcommand{\group}{{G_1}}

\newcommand{\topol}{{\'etale\ }}

\newtheorem{proposition}{Proposition}[section]
\newtheorem{theorem}[proposition]{Theorem}
\newtheorem{definition}[proposition]{Definition}

\newtheorem{lemma}[proposition]{Lemma}

\newtheorem{corollary}[proposition]{Corollary}
\newtheorem{remark}[proposition]{Remark}
\newtheorem{construction}[proposition]{Construction}

\numberwithin{equation}{section}

\title[Principal sheaves]
{Moduli space of principal sheaves over projective varieties}
\author[T. G\'omez, I. Sols]{Tom\'as L. G\'omez and Ignacio Sols}
\date{21 January 2003}

\address{
Departamento de Algebra, Facultad de Ciencias Matem\'aticas,
Universidad Complutense de Madrid, 28040 Madrid (Spain)}

\email{tgomez@alg.mat.ucm.es,sols@mat.ucm.es}

\begin{document}

\begin{abstract}
Let $G$ be a connected reductive group. The late Ramanathan
gave a notion of (semi)stable principal $G$-bundle on
a Riemann surface and constructed a
projective moduli space
of such objects. We generalize Ramanathan's notion
and construction to higher dimension, 
allowing also objects which we call semistable 
principal $G$-sheaves,
in order to obtain a projective moduli space:
a principal $G$-sheaf on a projective
variety $X$ is a triple $(P,E,\psi)$, where 
$E$ is a torsion free sheaf on $X$,
$P$ is a 
principal $G$-bundle on the open set $U$ where $E$ is locally free
and 
$\psi$ is an isomorphism
between $E|_U$ and the vector bundle
associated to $P$ by the adjoint representation.

We say it is (semi)stable if all 
filtrations $E_\bullet$ of $E$ as sheaf of (Killing)
orthogonal algebras, i.e. filtrations with
$E^\perp_i=E_{-i-1}$ and $[E_i,E_j]\subset 
E_{i+j}^{\quad\vee\vee}$,
have 
$$
\sum (P_{E_i}\rk E - P_E \rk E_i)\,(\preceq)\,0 ,
$$
where $P_{E_i}$ is the Hilbert polynomial of $E_i$.
After fixing the Chern classes of $E$ and of the 
line bundles associated to the principal bundle $P$
and characters of $G$,
we obtain a projective moduli space of semistable
principal $G$-sheaves.
We prove that, in case $\dim X=1$, our notion of
(semi)stability is equivalent
to Ramanathan's notion.
\end{abstract}

\maketitle

\hspace{7.5cm}\textit{To A. Ramanathan, in memoriam}
\bigskip

\bigskip

\section*{Introduction}

Let $X$ be a smooth projective variety of dimension $n$ over $\CC$, with a
very ample line bundle $\SO_X(1)$, and let $G$ be a connected
algebraic reductive group. 
A principal $\glrc$-bundle over $X$ is equivalent to a vector bundle
of rank $R$. If $X$ is a curve, the moduli space was constructed by 
Narasimhan and Seshadri \cite{N-S,Sesh}.
If $\dim X >1$, to obtain a projective moduli space 
we have to consider also torsion free sheaves, and this was 
done by Gieseker, Maruyama and Simpson \cite{Gi,Ma,Si}.
Ramanathan \cite{Ra1,Ra2,Ra3} defined a notion of stability
for principal $G$-bundles, and constructed the projective
moduli space of semistable principal bundles on a curve.

We equivalently reformulate in terms of 
filtrations of the associated adjoint
bundle of (Killing) orthogonal algebras
the Ramanathan's notion of (semi)stability, which is 
essentially of slope type 
(negativity of the degree of some associated line bundles), 
so when we generalize
principal bundles to higher dimension by allowing 
their adjoints to be
torsion free sheaves we are able to 
just switch degrees by Hilbert
polynomials as definition of (semi)stability.
We then construct a projective coarse moduli space of
such semistable principal $G$-sheaves.
Our construction 
proceeds by reductions to intermediate groups, as in 
\cite{Ra3}, 
although starting the chain higher, namely in a
moduli of semistable tensors (as constructed  in \cite{G-S1}).
In performing
these reductions we have switched the technique,
in particular studying the non-abelian 
\'etale cohomology sets
with values in the groups involved,
which provides a simpler proof also in Ramanathan's case 
$\dim X=1$.
However, for the proof of properness we have been able
to just generalize the idea of \cite{Ra3}.

In order to make more precise these notions and results,
let $G'=[G,G]$ be the commutator subgroup, and let
$\fg=\fz\oplus \fgp$ be the Lie algebra of $G$, where 
$\fg'$ is the semisimple part and $\fz$ is the center.
As notion of principal $G$-sheaf, it seems natural 
to consider a rational principal 
$G$-bundle $P$, i.e. a principal $G$-bundle on an open
set $U$ with $\codim X\setminus U \geq 2$, and
a torsion free extension of the form $\fz_X\oplus E$,
to the whole of $X$, of the vector bundle
$P(\fg)=P(\fz\oplus \fgp)=\fz_U\oplus P(\fgp)$
associated to $P$ by the adjoint representation of $G$ in $\fg$.
This clearly amounts to the following

\begin{definition}
\label{numero}
A principal $G$-sheaf $\SP$ over $X$ is a triple $\SP=(P,E,\psi)$
consisting of a torsion free sheaf $E$ on $X$, 
a principal $G$-bundle $P$ on the open set $U_E$ where
$E$ is locally free, and an isomorphism of vector bundles
$$
\psi:P(\fgp)\stackrel{\isom}{\too} E|_{U_E} \, .
$$
\end{definition}

Recall that the algebra structure of $\fgp$ given by the
Lie bracket provides $\fgp$ an orthogonal (Killing) structure, i.e.
$\kappa:\fgp\otimes \fgp\to \CC$ inducing an isomorphism
$\fgp\isom \fgp^\vee$. Correspondingly, the adjoint vector bundle
$P(\fgp)$ on $U$ has a Lie algebra structure
$P(\fgp)\otimes P(\fgp)\to P(\fgp)$ and an orthogonal
structure, i.e. $\kappa:P(\fgp)\otimes P(\fgp)\to \SO_U$
inducing an isomorphism $P(\fgp)\isom P(\fgp)^\vee$.
In lemma \ref{lemmapepsi}
it is shown that the Lie algebra structure 
uniquely extends to a homomorphism
$$
[,]: E\otimes E \too E^{\vee\vee} \, ,
$$
where we have to take $E^{\vee\vee}$ in the target
because an extension $E\otimes E\to E$ does not
always exist
(so the above definition of a principal $G$-sheaf
is equivalent to the one given in our announcement
of results \cite{G-S2}). Analogously, the 
Killing form extends uniquely to
$$
\kappa:E\otimes E \too \SO_X
$$
inducing an inclusion $E\inj E^\vee$.
This form assigns an orthogonal 
$F^\perp=\ker(E\inj E^\vee \surj F^\vee)$
to each subsheaf $F\subseteq E$.

\begin{definition}
\label{filtrations1}
An orthogonal algebra filtration of $E$ is a filtration 
\begin{equation}
\label{eqfilt1}
0 \subsetneq E_{-l} \subseteq E_{-l+1} \subseteq \cdots
\subseteq E_l=E
\end{equation}
with 
$$
(1)\quad E_i^\perp = E_{-i-1}^{}
\quad \text{and} \quad
(2)\quad [E_i,E_i] \subseteq E_{i+j}^{\quad\vee\vee} 
$$
for all $i$, $j$.
\end{definition}

We will see that, if  
$U$ is an open set with $\codim X\setminus U\geq 2$ such that
$E|_U$ is locally free,
a reduction 
of structure group of
the principal bundle $P|_{U}$ to a parabolic subgroup
$Q$ together with a dominant character of $Q$ produces a
filtration of $E$, and the filtrations arising in
this way are precisely the orthogonal algebra filtrations
of $E$
(lemmas \ref{bijection} and \ref{selforth}).
We define the Hilbert polynomial $P_{E_\bullet}$ of
a filtration $E_\bullet\subseteq E$ as
$$
P_{E_\bullet}=\sum ( r P_{E_i} - r_{i}P_E )
$$
where $P_{E}$, $r$, $P_{E_i}$, $r_i$ always denote
the Hilbert polynomials with respect to $\SO_X(1)$ 
and ranks of $E$ and $E_i$.
If $P$ is a polynomial, we write $P\prec 0$ 
if $P(m)<0$ for $m\gg 0$, and analogously for ``$\preceq$'' and
``$\leq$''.
We also use the usual convention: whenever ``(semi)stable'' and
``$(\preceq)$''  appear in a sentence, two statements
should be read: one with ``semistable'' and ``$\preceq$'' and another 
with ``stable'' and ``$\prec$''.

\begin{definition}
\label{stab1}
(See equivalent definition in lemma \ref{stabilitylie}) 
A principal $G$-sheaf $\SP=(P,E,\psi)$ is said to be 
(semi)stable if 
all orthogonal algebra filtrations 
$E_\bullet\subset E$ have
$$
P_{E_\bullet} (\preceq) 0
$$
\end{definition}

In proposition \ref{samestable} we prove that this
is equivalent to the condition that 
the associated tensor 
$$
(E,\phi: E\otimes E\otimes \wedge^{r-1}E \too \SO_X)
$$
is (semi)stable (in the sense of \cite{G-S1}).

To grasp the meaning of this definition, recall that suppressing
conditions (1) and (2) in definitions \ref{filtrations1} 
and \ref{stab1}
amounts to the (semi)stability of $E$ as a torsion free sheaf,
while just requiring condition (1) amounts to the
(semi)stability of $E$ as an orthogonal sheaf
(cfr. \cite{G-S2}). Now, demanding (1) and (2) is
having into account both the orthogonal and the algebra structure
of the sheaf $E$, i.e. considering its (semi)stability
as orthogonal algebra. By corollary \ref{stabilitylie},
this definition coincides with the one given in the
announcement of results \cite{G-S2}.

Replacing the Hilbert polynomials $P_E$ and $P_{E_i}$ by degrees we obtain
the notion of \textit{slope-(semi)stability}, which
in section \ref{sec5} will be 
shown to be 
\textit{equivalent
to the Ramanathan's notion of (semi)stability} 
\cite{Ra2,Ra3} 
\textit{of the rational principal
$G$-bundle $P$} (this has been written at the end
just to avoid interruption of the main argument of the article, and 
in fact we refer sometimes to section \ref{sec5} 
as a sort of appendix). Clearly 
$$
\text{slope-stable $\Longrightarrow$ stable
$\Longrightarrow$ semistable $\Longrightarrow$ 
slope-semistable}
$$

Since $G/G'\isom
\CC^{*q}$, given a principal $G$-sheaf, the principal
bundle $P(G/G')$ obtained by extension of structure group
provides $q$ line bundles on $U$, and since 
$\codim X\setminus U \geq 2$, 
these line bundles extend uniquely to line bundles on $X$.
Let $d_1,\ldots,d_q\in H^2(X,\CC)$ be their Chern classes. 
The rank $r$ of $E$ is clearly the dimension of $\fgp$.
Let $c_i$ be the Chern classes of $E$.

\begin{definition}[Numerical invariants]
\label{ginvariants}
We call the 
data $\tau=(d_1,\ldots,d_q,c_i)$ 
the numerical invariants of the principal
$G$-sheaf $(P,E,\psi)$. 
\end{definition}

\begin{definition}[Family of principal $G$-sheaves]
A family of (semi)stable principal $G$-sheaves 
parametrized by a complex scheme $S$
is a triple 
$(P_S,E_S,\psi_S)$ where $E_S$ is a torsion free sheaf
on $X\times S$, flat over $S$, 
$P_S$ is a principal $G$-bundle on 
the open set $U_{E_S}$ where $E_S$ is locally free, 
and $\psi:P_S(\fgp)\to E_S|_{U_{E_S}}$ is an isomorphism of
vector bundles.

Furthermore, it is asked that for
all closed points $s\in S$ the corresponding principal $G$-sheaf
is (semi)stable with numerical invariants $\tau$.
\end{definition}

An isomorphism between two such families $(P_S,E_S,\psi_S)$ and
$(P'_S,E'_S,\psi'_S)$ is a pair 
$$
(\beta:P^{}_S\stackrel{\isom}\too P'_S ,\gamma:E^{}_S \stackrel{\isom}\too E'_S)
$$
such that the following diagram is commutative
$$
\xymatrix{
{P_S(\fgp)} \ar[r]^{\psi} \ar[d]_{\beta(\fgp)} & 
{E_S|_{U_{E_S}}} \ar[d]^{\gamma|_{U_{E_S}}} \\
{P'_S(\fgp)} \ar[r]^{\psi'} & {E'_S|_{U_{E_S}}} 
}
$$
where $\beta(\fgp)$ is the isomorphism of vector bundles induced by
$\beta$.
Given an $S$-family $\SP_S=(P_S,E_S,\psi_S)$ and a morphism
$f:S'\to S$, the pullback is defined as
$\wt{f}^* \SP_S=(\wt{f}^*P_S ,\overline{f}^* E_S,\wt{f}^*\psi_S)$,
where $\overline{f}=\id_X\times f:X\times S  \to 
X\times S'$ and
$\wt{f}=i^*(\overline{f}):U_{\overline{f}^*E_S}\to U_{E_S}$,
denoting $i:U_{E_S}\to X\times S$ the inclusion of the open
set where $E_S$ is locally free.

\begin{definition}
The functor $\wt{F}^\tau_G$ is the sheafification of the 
functor
$$
F^\tau_G: (\Sch/\CC) \too (\Sets) 
$$
sending a complex scheme $S$, locally of finite type, 
to the set of isomorphism classes of families
of semistable principal 
$G$-sheaves with numerical invariants $\tau$, and it is
defined on morphisms as pullback.
\end{definition}

Let $\SP=(P,E,\psi)$ be a semistable principal 
$G$-sheaf on $X$. An orthogonal algebra
filtration $E_\bullet$ of $E$ which is 
\textit{admissible}, i.e. having 
$P_{E_\bullet} = 0$,
provides a reduction $P^Q$ 
of $P|_U$ to a parabolic subgroup $Q\subset G$
(lemma \ref{bijection}) on the open set
$U$ where it is a bundle filtration.
Let $Q\surj L$ be its Levi
quotient, and  $L\inj Q\subset G$ a splitting. 
We call
the semistable principal $G$-sheaf
$$
\big(
P^Q(Q\surj L\inj G),
\oplus E_i/E_{i-1}, 
\psi'
\big)
$$
the associated \textit{admissible deformation} of $\SP$,
where $\psi'$ is the natural isomorphism between
$P^Q(Q\surj L\inj G)(\fgp)$ and $\oplus E_i/E_{i-1}|_U$.
This principal $G$-sheaf is semistable.
If we iterate this process, it stops after a finite number
of steps, i.e. 
a semistable $G$-sheaf ${\rm grad}\,\SP$ (only depending on $\SP$) 
is obtained such that all its admissible deformations are
isomorphic to itself (cfr. proposition \ref{propsequiv}).

\begin{definition}
\label{sequiv}
Two semistable $G$-sheaves $\SP$
and $\SP'$ are said S-equivalent if 
${\rm grad}\,\SP\isom{\rm grad}\,\SP'$.
\end{definition}
 
When $\dim X=1$ this is just Ramanathan's notion of
S-equivalence of semistable principal $G$-bundles.
Our main result generalizes Ramanathan's \cite{Ra3} to 
arbitrary dimension:

\begin{theorem}
For a polarized projective variety $X$ there is a coarse projective
moduli space of S-equivalence classes of semistable $G$-sheaves on $X$
with fixed numerical invariants.
\end{theorem}

Principal $\GLR$-sheaves
are not objects equivalent to 
torsion free sheaves of rank $R$, but only in the case of bundles.
As we show in section \ref{secexample},
even in this case,
the (semi)stability of both objects do not coincide.
The philosophy is that, just as Gieseker changed 
in the theory of stable vector bundles both
the objects (torsion free sheaves instead of vector bundles)
and the condition of (semi)stability (involving Hilbert
polynomials instead of degrees) 
in order to make $\dim X$ a parameter of the
theory, it is now needed to change again the objects (principal 
sheaves) and the condition of (semi)stability 
(as that of the adjoint sheaf of orthogonal algebras) 
in order to make the group $G$
a parameter of the theory (such variations of the 
conditions of stability and semistability are in 
both generalizations very slight, 
as these are implied by slope stability and imply
slope semistability, and the slope 
conditions do not vary).
The deep reason is that what we intend to do
is not generalizing the notion of 
vector bundle of rank $R$ (which
was the task of Gieseker and Maruyama), but
that of principal $\GLR$-bundle, and although 
both notions happen to be
extensionally the same, i.e. happen to define 
equivalent objects, they are essentially different.
This subtle fact is recognized by the
very sensitive condition of existence of a moduli
space, i.e. by stability.

The results of this article where announced in \cite{G-S2}.
There is independent work by Hyeon \cite{Hy}, who constructs,
for higher dimensional varieties,
the moduli space of principal bundles whose associated adjoint
is a Mumford stable vector bundle, 
using the techniques of Ramanathan \cite{Ra3},
and also by Schmitt \cite{Sch} who
chooses a faithful representation of $G$ in 
order to obtain and compactify a moduli space of principal
$G$-bundles.

\medskip
\noindent\textbf{Acknowledgments.}
We would like to thank M.S. Narasimhan for suggesting this
problem in a conversation with the first author in ICTP (Trieste)
in August 1999 and for discussions. 
We would also like to thank J.M. Ancochea, O. Campoamor, 
N. Fakhruddin, 
S. Ilangovan, J. M. Marco,
V. Mehta, A. Nair, 
N. Nitsure, S. Ramanan, T.N. Venkataramana and A. Vistoli for 
comments and fruitful discussions.

The authors are members of VBAC (Vector Bundles on Algebraic Curves),
which is partially supported by EAGER (EC FP5 Contract no.
HPRN-CT-2000-00099) and by EDGE (EC FP5 Contract
no. HPRN-CT-2000-00101).
T.G. was supported by a postdoctoral fellowship of
Ministerio de Educaci\'on y Cultura (Spain), and wants to thank the
Tata Institute of Fundamental Research (Mumbai, India), where
this work was done while he was a postdoctoral student.
I.S. wants to thank the very warm hospitality of the members of 
the Institute during his visit to Mumbai.

\section*{Preliminaries}

\noindent\textbf{Notation.}
We denote by $(\Sch/\CC)$ the category of schemes over
$\Spec \CC$, locally of finite type. All schemes considered
will belong to this category. 
If $f:Y \to Y'$ is a morphism,
we denote $\overline{f}=\id_X\times f:X\times Y \to X\times Y'$.
If $E_S$ is a coherent sheaf on $X\times S$, we denote
$E_S(m):= E_S \otimes p^*_X\SO_X(m)$.
An open set $U\subseteq Y$ of a scheme $Y$ will be called
\textit{big} if $\codim Y\setminus U \geq 2$.
Recall that in the \'etale topology, 
a cover of a scheme $U$ is a finite collection 
of morphisms $\{f_i:U_i\to U\}_{i\in I}$ such that
each $f_i$ is \'etale, and $U$ is the
union of the images of the $f_i$.
\medskip

Given a principal $G$-bundle $P\to Y$ and a left 
action $\sigma$ of $G$ in
a scheme $F$, we denote
$$
P(\sigma,F) \;:=\; P\times_G F \;=\; (P\times F)/G,
$$
the associated fiber bundle.
If the action $\sigma$ is clear from the context, we will 
write $P(F)$. 
If $\rho:G\to H$ is a group homomorphism, let $\sigma$ be
the action of $G$ on $H$ defined by left multiplication
$h\mapsto \rho(g)h$. Then the associated fiber bundle
is a principal $G'$-bundle, and it is denoted
$\rho_*P$.
If $\sigma$ is a character of the group,
we will denote by $P(\sigma)$ the corresponding line bundle.

Let $\rho:H \to G$ be a homomorphism of groups, and let $P$ be a
principal $G$-bundle on a scheme $Y$. 
A reduction of structure group of $P$ 
to $H$ is a pair $(P^H,\zeta)$, where $P^H$ is 
a principal $H$-bundle on $Y$ and $\zeta$ is an isomorphism
between $\rho_* P^H$ and $P$. 
Two reductions $(P^{H}_T,\zeta^{}_T)$
and $(Q^{H}_T,\theta^{}_T)$ are isomorphic if there is an isomorphism
$\alpha$ giving a commutative diagram
\begin{equation}
\label{isopair}
\xymatrix{{P^H_T} \ar[d]^{\alpha}_{\cong} \\ 
{Q^{H}_T}  }
\qquad
\xymatrix{
  {\rho_{*} P^H_T} \ar[r]^{\zeta^{}_T}
\ar[d]^{\rho_{*}{\alpha}} & {P^{}_T} \ar@{=}[d] \\
 {\rho_{*} Q^{H}_T} \ar[r]^{\theta^{}_T} & {P^{}_T} 
}
\end{equation}

Let $p:Y\to S$ be a morphism of schemes, and let $P_S$ be a
principal $G$-bundle on the scheme $Y$. Define the functor
of families of reductions
\begin{eqnarray*}
\label{reduction2}
\Gamma(\rho,P_S): ({\rm Sch}/S) & \too & ({\rm Sets})\\
(t:T \too S) & \longmapsto & \big\{ (P^{H}_T,\zeta^{}_T) 
\big\}/\text{isomorphism}
\end{eqnarray*}
where $(P^{H}_T,\zeta^{}_T)$ is a reduction of structure
group of $P_T:=P_S\times_{S} T$ to $H$.

If $\rho$ is injective, then $\Gamma(\rho,P_S)$ is a sheaf,
and it is in fact representable by a scheme $S'\to S$,
locally of finite type \cite[lemma 4.8.1]{Ra3}.
If $\rho$ is not injective, this functor is not necessarily a
sheaf, so we denote by 
$\wt\Gamma(\rho,P_S)$ its sheafification 
with respect to the \topol topology on $(\Sch/S)$.

\begin{lemma}
\label{zeroes}
Let $Y$ be a scheme, and let $f:\SK\to \SF$ be a 
homomorphism of sheaves on $X\times Y$. Assume that
$\SF$ is flat over $Y$. 
Then there is a unique closed subscheme $Z$
satisfying the following universal property:
given a Cartesian diagram
$$
\xymatrix{
{X\times S} \ar[r]^-{\overline{h}} \ar[d]^{{p_S}}& 
{X\times Y} \ar[d]^{p}
\\
{S} \ar[r]^-{h} & {Y} 
}
$$
$\overline{h}^* f=0$ if and only if $h$ factors through $Z$.

\end{lemma}

\begin{proof}
Uniqueness is clear.
To show existence, assume that $\SO_X(1)$ is very ample 
(taking a multiple if necessary).
Recall that if $\SG$ is a coherent sheaf on $X\times Y$, we denote 
$\SG(m)=\SG\otimes p_X^*\SO_X^{}(m)$.
Since $\SF$ is $Y$-flat, taking $m'$ 
large enough, $p_*^{} \SF(m')$ is locally 
free. The question is local on $Y$,
so we can
assume, shrinking $Y$ if necessary, 
that $Y=\Spec A$ and $p_*^{} \SF(m')$ is given by 
a free $A$-module.
Now, since $Y$ is affine, the homomorphism
$$
p^{}_*f(m'): p_*^{} \SK(m') \too p_*^{} \SF(m')
$$
of sheaves on $Y$ is equivalent to a homomorphism of 
$A$-modules
$$
M \stackrel{(f_1,\ldots,f_n)}{\too} A\oplus \cdots\oplus A
$$
The zero locus of $f_i$ is defined by the ideal $I_i\subset A$
image of $f_i$, thus the zero scheme of $(f_1,\ldots,f_n)$ 
is given by the ideal $I=\sum I_i$,
hence $Z'_{m'}$ is a closed subscheme.

Since $\SO_X(1)$ is very ample, if $m''>m'$
we have an injection $p_*^{} \SF(m')\inj p_*^{} \SF(m'')$
(and analogously for $\SK$), hence $Z_{m''}\subset Z_{m'}$,
and since $Y$ is noetherian, there exists $N'$ such
that, if $m'>N'$,
we get a scheme $Z$ independent of $m'$.

To check the universal property first we will show that if
$\overline{h}^*f=0$ then $h$ factors through $Z$.  Since the question
is local on $S$, we can take $S=\Spec(B)$, $Y=\Spec(A)$, 
and the morphism $h$ is locally
given by a ring homomorphism $A\to B$.  Since $\SF$ is flat over $Y$,
for $m'$ large enough the natural homomorphism $\alpha:h^*p_*\SF(m')\to
{p_S}_*\overline{h}^*\SF(m')$ (defined as in \cite[Th. III 9.3.1]{Ha})
is an isomorphism.  Indeed, for $m'$ sufficiently large,
$H^i(X,\SF_y(m'))=0$ and $H^i(X,\overline{h}^* (\SF(m'))_s)=0$ for all
points $y\in Y$, $s\in S$ and $i>0$, and since $\SF$ is flat, this
implies that $h^*p_*\SF(m')$ and ${p_S}_*\overline{h}^*\SF(m')$ are
locally free. Then, to prove that the homomorphism $\alpha$ 
is an isomorphism,
it is enough to check it at the fiber of every $s\in S$, but this
follows from \cite[Th. III 12.11]{Ha} or \cite[II \S 5 Cor. 3]{Mu2}.

Hence the commutativity of
the diagram
$$
\xymatrix{
{{p_S}_*\overline{h}^*\SK(m')} 
\ar[rrr]^{{p_S}_*\overline{h}^*f(m')=0}&&&
{{p_S}_*\overline{h}^*\SF(m')}  \\
{h^*p_*\SK(m')} \ar[rrr]^{h^*p_*f(m')} \ar[u] &&&
{h^*p_*\SF(m')} \ar[u]_{\isom}
}
$$
implies that $h^*p_*f(m')=0$. This means that 
for all $i$,  in the
diagram
$$
\xymatrix{
{M} \ar[r]^{f_i}\ar[d] & {A} \ar[r]\ar[d] & {A/I_i}\ar[d] \ar[r] & {0}\\
{M\otimes_A B} \ar[r]^-{f_i\otimes B} 
& {B} \ar[r] & {A/I_i\otimes_A B}\ar[r] &{0}
}
$$
it is $f_i\otimes B=0$.
Hence the image $I_i$ of
$f_i$ is in the kernel $J$ of $A\to B$. Therefore $I\subset J$,
hence $A\to B$ factors through $A\to A/I$, which means that
$h:S\to Y$ factors through $Z$.

Now we show that if we take $S=Z$ and $h:Z\inj Y$ the inclusion,
then $\overline{h}^* f=0$.
By definition of $Z$ we have $h^*p_* f(m')=0$ for any $m'$ with 
$m'>N'$.
Showing that $\overline{h}^* f=0$
is equivalent to showing that
$$
\overline{h}^* f(m')  : \overline{h}^* \SK(m') \too \overline{h}^* \SF(m')
$$ 
is zero for some $m'$.
Take $m'$ large enough so that 
$\ev:p^*p_*\SK(m')\to \SK(m')$ is surjective.
By the right exactness of $\overline{h}^*$ the homomorphism 
$\overline{h}^*\ev$ is still surjective. The commutative
diagram
$$
\xymatrix{
{\overline{h}^* \SK(m')} \ar[rrr]^{\overline{h}^* f(m')} &&& 
{\overline{h}^* \SF(m')}\\ 
{\overline{h}^*p^*_{}p_*  \SK(m')} 
\ar[rrr]^{\overline{h}^*p^*_{}p_*  f(m')}
\ar@{>>}[u]^{\overline{h}^*\ev} &&& 
{\overline{h}^*p^*_{}p_* \SF(m')} \ar[u]\\
{p_S^{*}h^*_{}p_*  \SK(m')} 
\ar[rrr]^{p_S^{*}h^*_{}p_*  f(m')=0} \ar@{=}[u]&&& 
{p_S^{*}h^*_{}p_* \SF(m')} \ar@{=}[u]
}
$$
implies $\overline{h}^*f(m')=0$, hence $\overline{h}^*f=0$.

\end{proof}

The following two lemmas and corollary will help to relate
the three main objects that will be introduced in this section.

\begin{lemma}
\label{lemma1}
Let $E$ and $F$ be coherent sheaves on a scheme $Y$, and $L$ a locally
free sheaf on $Y$. There is a
natural isomorphism
$$
\Hom (E\otimes F, L) \cong 
\Hom (E,\SH om(F,L)) \cong
\Hom (E,F^\vee\otimes L)
$$
\end{lemma}

\begin{proof}
Given a homomorphism $\varphi:E\otimes F \to L$, 
we define $\psi:E\to {\SH om}(F,L)$ by sending 
a local section $e$ of $E$ to $\varphi(e,\cdot)$.
Conversely, to a homomorphism $\psi:E\to {\SH om}(F,L)$, we
associate the homomorphism 
$$
E\otimes F \stackrel{\psi\otimes F}{\too} 
{\SH om}(F,L) \otimes F \too L
$$
where the second map is the natural pairing. It is easy to check
that both constructions are inverse to each other.
Finally, since $L$ is locally free, ${\SHom}(F,L)=F^\vee\otimes L$.
\end{proof} 

\begin{lemma}
\label{lemma2}
Let $E$ be a torsion free sheaf of rank $r$ on a scheme $Y$. Then
there is a canonical isomorphism
$$
E^{\vee\vee} \cong (\bigwedge{}^{r-1} E)^\vee \otimes \det E \,.
$$
\end{lemma}

\begin{proof}
This isomorphism is obvious if we restrict to the open set $U_E$ where
$E$ is locally free. Since both sheaves are reflexive
and $\codim X\setminus
U_{E}\geq 2$, it uniquely extends to an isomorphism on the whole of
$X$ (\cite[Prop. 1.6(iii)]{Ha2}).
\end{proof}

Combining lemmas \ref{lemma1} and \ref{lemma2} we obtain the
following

\begin{corollary}
\label{cor12}
Let $E$ be a torsion free sheaf of rank $r$ on a scheme $Y$, and $L$ a
line bundle on $Y$. Then, giving 
a homomorphism
$$
\varphi: E\otimes E \too E^{\vee\vee}\otimes L
$$
is equivalent to giving a homomorphism
$$
\phi: E\otimes E \otimes  E^{\otimes r-1}\;=E^{\otimes r+1} \too \det E
\otimes L
$$
which is skew-symmetric in the last $r-1$ entries, i.e.
which induces a homomorphism on $E\otimes E \otimes \bigwedge{}^{r-1} E$.
\end{corollary}

Now we introduce the progressively richer concepts of tensor, 
$\fgp$-sheaf, and principal $G$-sheaf, defining them relative to a
scheme $S$. As usual, if no mention to the base scheme $S$
is made,
it will be understood $S=\Spec \CC$.
For each of these three concepts we give compatible
notions of (semi)stability, leading in each case to a projective coarse
moduli space.

\begin{definition}[Tensor]
\label{deftensor}
A family of tensors parametrized by a scheme $S$ is a
triple $(E_S,\phi_S,N)$ 
consisting of a torsion free sheaf $E_S$ on 
$X\times S$, flat over $S$, which restricts to a torsion free sheaf
with trivial determinant and fixed Hilbert polynomial $P$
on every slice $X\times s$, 
a line bundle $N$ on $S$ and
a homomorphism $\phi_S$  
\begin{equation}
\phi^{}_S:E^{}_S{}^{\otimes a}
\, \too \,  p_S^*N,
\end{equation}

A tensor is called a Lie tensor if $a=r+1$ for $r$ the rank of $E_S$, 
and 
\begin{enumerate}
\item $\phi_S$ is skew-symmetric in the last $r-1$ entries, i.e.
it induces a homomorphism on $E_S\otimes E_S \otimes \bigwedge{}^{r-1}
E_S$,

\item the homomorphism 
$\varphi_S:E_S\otimes E_S\to E^{\vee\vee}_S\otimes N$ 
associated to $\phi_S$ by corollary \ref{cor12}
is antisymmetric,

\item $\varphi_S$ satisfies the Jacobi identity
\end{enumerate}
\end{definition}

To give a precise definition of the Jacobi identity,
first define a homomorphism
$$
[[\cdot,\cdot],\cdot]: E_S\otimes E_S \otimes E_S 
\stackrel{\varphi_S \otimes E_S}{\too}  
E_S^{\vee\vee}  \otimes N \otimes E_S
\stackrel{E_S^{\vee\vee}\otimes N \otimes \varphi_S}{\too} 
E_S^{\vee\vee} \otimes E_S^{\vee}
\otimes E_S^{\vee\vee} \otimes N^2
\too E_S^{\vee\vee} \otimes N^2
$$
where the last map comes from the natural pairing of the first two
factors. Then define
\begin{eqnarray}
\label{jacobin}
J: E_S\otimes E_S \otimes E_S  &\too & E_S^{\vee\vee} \otimes N^2 \\
(u,v,w) &\longmapsto & [[u,v],w]+[[v,w],u]+[[w,u],v] \nonumber
\end{eqnarray}
and require $J=0$.

An isomorphism between two
families of tensors $(E_S,\phi_S,N)$ and
$(E'_S,\phi'_S,N')$ parametrized by $S$
is a pair of isomorphisms
$\alpha:E_S\to E'_S$ and 
$\beta:N\to N'$ such that the induced diagram
$$
\xymatrix{
{E^{}_S{}^{\otimes r+1}} 
\ar[rr]^{\phi_S} \ar[d]^{\alpha^{\otimes r+1}} 
& & { p_S^*N} \ar[d]^{  \beta }\\
{E'_S{}^{\otimes r+1}} \ar[rr]^{\phi'_S} 
& &  { p_S^*N'}
}
$$
commutes.
In particular, $(E,\phi)$ and $(E,\lambda\phi)$ are isomorphic
for $\lambda\in \CC^*$.
Given an $S$-family of tensors $(E_S,\phi_S,N)$
and a morphism $F:S'\to S$, the pullback is the 
$S'$-family defined
as $(\wt{f}^* E_S, \wt{f}^* \phi_S, f^* N_S)$.

Since we will work with GIT (Geometric invariant theory, \cite{Mu1}), 
the notion of filtration $E_\bullet$ of
a sheaf is going to be essential for us. By this we always
understand a $\ZZ$-indexed filtration
$$
\ldots \subseteq E_{i-1}\subseteq E_{i}
\subseteq E_{i+1}\subseteq \ldots
$$
starting with $0$ and ending with $E$. Of course, only a finite
number of inclusions can be strict
$$
0 \subsetneq E_{\lambda_1} \subsetneq E_{\lambda_2} 
\subsetneq \;\cdots\; \subsetneq E_{\lambda_t} \subsetneq
E_{\lambda_{t+1}}=E \qquad \lambda_1<\cdots<\lambda_{t+1}
$$
where we have deleted, from $0$ onward, all the non-strict 
inclusions. 
Reciprocally,
from $E_{\lambda_\bullet}$ we recover $E_\bullet$ by 
defining $E_m=E_{\lambda_{i(m)}}$, 
\textit{where $i(m)$ is the maximum 
index with $\lambda_{i(m)}\leq m$.}

\begin{definition}[Balanced filtration]
\label{defbalancedfiltration}
A filtration $E_{\bullet}\subseteq E$ 
of a torsion free sheaf $E$ is called balanced 
if $\sum i \rk{E^i}=0$ for $E^i=E_i/E_{i-1}$.
In terms of $E_{\lambda_\bullet}$,
this is $\sum_{i=1}^{t+1} \lambda_i \rk (E^{\lambda_i})=0$
for $E^{\lambda_i}=E_{\lambda_i}/E_{\lambda_{i-1}}$.
\end{definition}

\begin{remark}
\label{gitsheaves}
\textup{
The notion of balanced filtration appeared naturally
in the Gieseker-Maruyama construction of the moduli space of
(semi)stable sheaves, the condition of (semi)stability
of a sheaf $E$ being that all balanced filtrations of
$E$ have negative (nonpositive) Hilbert polynomial.
In this case the condition ``balanced'' can be 
suppressed, since $P_{E_\bullet}=P_{E_{\bullet+l}}$
for any shift $l$ in the indexing
(and furthermore it is enough to consider filtrations
of one element, i.e. just
subsheaves).}
\end{remark}


Let $\SI_a=\{1,\ldots,t+1\}^{\times a}$ be the set of all
multi-indexes $I=(i_1,\ldots,i_a)$ of cardinality $a$. 
Define
\begin{equation}
\label{muE}	
\mu_{\rm tens}(\phi,E_{\lambda_\bullet})=\min_{I\in \SI_a} \big\{
\lambda_{i_1}+\dots+\lambda_{i_a}: \,
\phi|_{E_{\lambda_{i_1}}\otimes\cdots \otimes E_{\lambda_{i_a}}}\neq 0
  \big\}
\end{equation}

\begin{definition}[Stability for tensors]
\label{stabilitytensor}
Let $\delta$ be a polynomial of degree at most $n -1$ and positive
leading coefficient. 
We say that $(E,\phi)$ is $\delta$-(semi)stable if 
$\phi$ is not identically zero and
for all balanced
filtrations $E_{\lambda_\bullet}$ of $E$, it is 
\begin{equation}
\label{stabtensor}
\Big(\sum_{i=1}^t (\lambda_{i+1}-\lambda_i)
\big( r P_{E_{\lambda_i}} -r_{\lambda_i} P \big)\Big) +
\mu_{\rm tens}(\phi,E_{\lambda_\bullet}) \, \delta \;(\preceq)\; 0
\end{equation}
\end{definition}

In this definition, \textit{it suffices to consider saturated
filtrations}, 
which means that
the sheaves $E^i=E_i/E_{i-1}$ (i.e. the
$E^{\lambda_i}=E_{\lambda_i}/E_{\lambda_{i-1}}$) 
are all torsion free, 
and thus with $\rk(E_{\lambda_i})<\rk(E_{\lambda_{i+1}})$ for all $i$.
There is a coarse moduli space of semistable
tensors \cite{G-S1}.


Now we go to our second main concept, that of a
$\fgp$-sheaf.
A family of Lie algebra 
sheaves 
parametrized by $S$ is a pair
$$
\big(E_S,\; \varphi_S:E_S\otimes E_S\to E_S^{\vee\vee}),
$$
where $E_S$ is a torsion free sheaf on $X\times S$, flat over $S$,
and $\varphi_S$ is a homomorphism such that 
\begin{enumerate}
\item $\varphi_S:E_S\otimes E_S\too E^{\vee\vee}_S$ 
is antisymmetric

\item $\varphi_S$ satisfies the Jacobi identity
\end{enumerate}

The precise definition of the Jacobi identity is
as in definition \ref{deftensor}, but with 
$\SO_{X\times S}$ instead of $N$.
An isomorphism between two families 
is an isomorphism 
$\alpha:E^{}_S \to E'_S$ 
with
$$
\xymatrix{
{E^{}_S\otimes E^{}_S} \ar[r]^{\varphi^{}_S} \ar[d]^{\alpha\otimes \alpha} &  
{E_S^{\vee\vee} } \ar[d]^{\alpha^{\vee\vee}}\\
{E'_S\otimes E'_S} \ar[r]^{\varphi'_S} &  
{E_S^{\prime\vee\vee} } 
}
$$

Note that since conditions (1) and (2) above 
are closed, it is not enough
to check that they are satisfied for all closed points of $S$, 
because $S$ could be nonreduced.


\begin{definition}
\label{defkilling}
The Killing form $\kappa_S$ associated to a Lie algebra sheaf
$(E_S,\varphi_S)$
is the composition
$$
{E_S\otimes E_S} 
 \stackrel{[,]\otimes[,]}{\too}
{(E_S^\vee\otimes E_S^{\vee\vee})\otimes (E_S^\vee\otimes
E_S^{\vee\vee})}
\too 
{(E_S^\vee\otimes E_S^{\vee\vee})} \too
\SO_{X\times S}
$$
\end{definition}

If the Lie algebra is semisimple, in the sense that the
induced homomorphism $E_S^{\vee\vee} \to E_S^\vee$ is
an isomorphism, the fiber of
$E_S$ over a point $(x,s)\in X\times S$ where $E_S$ is locally free
has the structure of a semisimple Lie algebra,
which, because of the rigidity of semisimple Lie algebras,
must be constant on connected components of $S$. This
justifies the following

\begin{definition}[$\fgp$-sheaf]
\label{deffgpsheaf}
A family of $\fgp$-sheaves is a family of Lie algebra
sheaves where the Lie algebra associated to each
connected component of the parameter space 
$S$ is $\fgp$.
\end{definition}


The following is the sheaf version of the well known
notion of Lie algebra filtration (see \cite{J} for instance,
recalled in section \ref{sec5}).

\begin{definition}[Algebra filtration]
A filtration $E_{\bullet}\subseteq E$ 
of a Lie algebra sheaf $(E,[,])$ is called an algebra filtration
if for all $i$, $j$,
$$
[E_{i}, E_{j}] \subseteq E_{i+j}^{\quad\vee\vee}.
$$
In terms of $E_{\lambda_\bullet}$, this is
$$
[E_{\lambda_i}, E_{\lambda_j}] \subseteq E_{\lambda_{k-1}}^{\quad\vee\vee}
$$
for all $\lambda_i$, $\lambda_j$, $\lambda_k$ with
$\lambda_i+\lambda_j<\lambda_k$,
\end{definition}

\begin{definition}
\label{stabilityliealgebra}
A $\fgp$-sheaf is (semi)stable if
for all balanced algebra filtrations $E_{\bullet}$
it is
$$
\sum_{i=1}^{t} 
\big( r P_{E_{i}} -r_{i} P_E \big)
\;(\preceq)\; 0 
$$
or, in terms of $E_{\lambda_\bullet}$,
\begin{equation}
\label{stablie}
\sum_{i=1}^{t} (\lambda_{i+1}-\lambda_i)
\big( r P_{E_{\lambda_i}} -r_{\lambda_i} P_E \big)
\;(\preceq)\; 0 
\end{equation}
\end{definition}

As usual, in definition \ref{stabilityliealgebra} it is enough to
consider saturated filtrations. 

\begin{remark}
\label{remstar}
\textup{
We will see in corollary \ref{selforth} that for an algebra
filtration of a $\fgp$-sheaf, the fact of being balanced 
and saturated is equivalent to being orthogonal,
i.e. $E_{-i-1}^{}=E^\perp_i=
\ker(E\inj E^{\vee\vee} \stackrel{\kappa}{\isom} E^\vee
\to E_i^\vee)$. Thus, in the previous definition we can
change ``balanced algebra filtration'' by 
``orthogonal algebra filtration''.} 
\end{remark}

\begin{remark}
\textup{
Observe that the condition ``balanced'' cannot be suppressed in this
case, is it was in remark \ref{gitsheaves}, because
a shifted filtration $E_{\bullet+l}$ of an algebra filtration is no
longer an algebra filtration.}
\end{remark}


\begin{construction}[Correspondence between 
Lie algebra sheaves and Lie tensors]
\label{constr}
\textup{
Consider a family of Lie tensors
$$
( F_S, \phi_S: F^{}_S{}^{\otimes r+1} \too p^*_S N, N)
$$
Corollary \ref{cor12} gives
$$
( F_S, \phi'_S: F_S\otimes F_S \to
F^{\vee\vee}_S\otimes(\det F_S)^{-1}\otimes p^*_S N, N)
$$
If we 
tensor  $\phi'_S$ with $(\det F_S)^2\otimes p^*_S N^{-2}$ and
define $E_S=F_S\otimes \det F_S \otimes p^*_S N^{-1}$,
we obtain a Lie algebra sheaf
\begin{equation}
\label{rliesheaves}
( E_S, \varphi_S: E_S\otimes E_S \to E^{\vee\vee}_S)
\end{equation}
such that for all $s\in S$ and $x\in U_{E_s}$, 
$\varphi_s(x)$ is a Lie algebra structure on the
fiber $E_s(x)$. }

\textup{
Conversely, given a Lie algebra sheaf as in 
(\ref{rliesheaves}),
corollary \ref{cor12} gives a homomorphism
$$
\phi_S^{}: E^{}_S{}^{\otimes r+1} \too \det E_S^{} \,.
$$
Recall that we always assume that $E_S=p^* L$,
where $L$ is a line bundle on the base scheme $S$,
hence this is a Lie tensor.}
\end{construction}

If $S=\Spec \CC$, this gives a bijection of isomorphism
classes, but not for arbitrary $S$, 
because $E_S$ is not in general isomorphic
to $F_S$. They are only locally isomorphic, in
the sense that we can cover $S$ with open sets $S_i$ (where
the line bundles $L$ and $N$ are trivial), so that 
the objects restricted to $S_i$ are isomorphic, and
this suffices to provide an
\textit{isomorphism between the
sheafified functors.
We will show that, for $\fgp$-sheaves,
its (semi)stability is equivalent to that of the
corresponding tensor, hence there is a projective moduli space
of $\fgp$-sheaves.}
This is the key initial point of this article, allowing
us to use in section \ref{sec1}
the results in \cite{G-S1} to construct
the moduli space of $\fgp$-sheaves.


Recall now, from the introduction, the notion 
of a principal $G$-sheaf $\SP=(P_S,E_S,\psi_S)$ for a reductive
connected group $G$ and its notion of (semi)stability. 
Let $\fgp$ be the semisimple
part of its Lie algebra. We associate now to
$\SP$ a $\fgp$-sheaf $(E_S,\varphi_S)$ by the following

\begin{lemma}
\label{lemmapepsi}
let $\SU=U_{E_S}$ be the open set where $E_S$ is locally free. 
The homomorphism
$\varphi_{\SU}:E_S|_{\SU}\otimes E_S|_{\SU}\to E_S|_{\SU}$, given by the Lie
algebra structure of $P_S(\fgp)$ and the isomorphism $\psi_S$,
extends uniquely to a homomorphism
$$
\varphi_S: E_S\otimes E_S \too E_S^{\vee\vee},
$$
\end{lemma}

\begin{proof}
The homomorphism 
$\varphi_{\SU}^{}$ can be seen as a section of
$$
E_S|_\SU^\vee \otimes E_S|_\SU^\vee \otimes E_S|_\SU
\isom (E_S|_\SU \otimes E_S|_\SU \otimes E_S|_\SU^\vee)^\vee
\isom (E_S \otimes E_S \otimes E_S^\vee)^\vee|_\SU^{},
$$ 
Since $(E_S \otimes E_S \otimes E_S^\vee)^\vee$ is a reflexive sheaf
on $X\times S$, this section extends uniquely to 
an element of
$$
\Gamma(X\times S,(E_S \otimes E_S \otimes E_S^\vee)^\vee)
= \Hom(E_S \otimes E_S \otimes E_S^\vee, \SO_{X\times S})
= \Hom(E_S \otimes E_S, E_S^{\vee\vee}),
$$
where the two equalities follow from corollary \ref{cor12},
and this element is the extended homomorphism $\varphi_S$.
\end{proof}


The following  corollary of remark \ref{remstar}
provides an equivalent definition
of (semi)stability

\begin{corollary}
\label{stabilitylie}
A principal $G$-sheaf $\SP=(P,E,\psi)$ is 
(semi)stable (definition \ref{stab1}) if and only if
the associated $\fgp$-sheaf $(E,\varphi)$ 
is (semi)stable (definition \ref{stabilityliealgebra}).
\end{corollary}


\begin{remark}
\label{fgpaut}
\textup{
Lemma \ref{lemmapepsi} implies  
that there is a natural bijection between the isomorphism
classes of families of $\fgp$-sheaves and those of principal
$\Aut(\fgp)$-sheaves.}
\end{remark}


\begin{lemma}
\label{lemma0}
Let $G$ be a connected reductive algebraic group. Let $P$ be a 
principal $G$-bundle on $X$ and let $E=P(\fgp)$ be 
the vector bundle associated to $P$ by the 
adjoint representation of $G$ on the semisimple part of 
its Lie algebra $\fgp$.
Then $\det E\cong \SO_X$.
\end{lemma}

\begin{proof}
We have $\autlie (\fgp) \subset \ogp$, where the orthogonal structure
on $\fgp$ is given by its nondegenerate Killing form. Note that $P(\fgp)$ is
obtained by extension of structure group using the composition
$$
\rho: G \too \autlie (\fgp) \inj \ogp \inj \glgp.
$$
Since $G$ is connected, the image of $G$ in $\ogp$ lies in the 
connected component of identity, i.e. in $\sogp$. Hence
$P(\fgp)$ admits a reduction of structure group to $\sogp$,
and thus $\det P(\fgp)\cong \SO_X$.
\end{proof}


We end this section by extending to principal sheaves some well 
known definitions and
properties of principal bundles and by recalling some notions
of GIT \cite{Mu1}.
Let $m:H\times R \to R$ be an action of  
an algebraic group $H$ on a scheme $R$. 
Let $p_R:H\times R \to R$ be the projection to the second factor.

\begin{definition}[Universal family]
Let $\SP_R$ be a family of principal $G$-sheaves parametrized by
$R$. Assume there is a lifting of the action of $H$ to $\SP_R$,
i.e. there is an isomorphism 
$$
\Lambda: \overline{m}^*\SP_R^{}
 \stackrel{\isom}{\too} \overline{p}_R^*\SP_R^{}
$$
Assume that
\begin{enumerate}
\item Given a family $\SP_S$ parametrized by $S$ and a closed 
point $s\in S$,
there is an open \topol neighborhood $i:S_0\inj S$ of $s$ and a 
morphism $t:S_0\to R$ such that $\overline{i}^*\SP_S \;\isom \; 
\overline{t}^*\SP_R$.

\item Given two morphisms $t_1, \, t_2:S\to R$ and an isomorphism
$\beta:\overline{t_2}^*\SP \to \overline{t_1}^*\SP$, 
there is a unique $h:S\to H$
such that $t_2=h[t_1]$ and $\overline{(h\times t_1)}^*
\Lambda=\beta$.
\end{enumerate}

Then we say that $\SP_R$ is a universal family with group $H$ for
the functor $\wt{F}^\tau_G$.
\end{definition}

\begin{definition}[Universal space]
Let $F:(\Sch/\CC)\to (\Sets)$ be a functor. 
Let $\underline{R}/\underline{H}$ be the sheaf on $(\Sch/\CC)$
associated to the presheaf $S \mapsto \Mor(S,R)/\Mor(S,H)$.
We say that
$R$ is a universal space with group $H$ for the functor
$F$ if the sheaf $F$ is isomorphic to $\underline{R}/\underline{H}$.
\end{definition}

The difference between these two notions can be understood
as follows. 
Given a stack $\SM:(\Sch/\CC)\to ({\rm Groupoids})$, we
denote by $\overline{\SM}:(\Sch/\CC)\to (\Sets)$ the
functor associated by taking the set of isomorphism classes
of each groupoid. 
Let $[R/H]$ be the quotient stack and let $\SF$ be
the stack of semistable principal $G$-sheaves.
Then $R$ is a universal space with group $H$ if 
$\overline{[R/H]}\isom \overline{\SF}$, whereas it is
a universal family if $[R/H]\isom \SF$, i.e. if the
isomorphism holds at the level of stacks, without
taking isomorphism classes.

\begin{definition}[Categorical quotient]
A morphism $f:R\to Y$ of schemes 
is a categorical quotient for an action 
of an algebraic group $H$ 
on $R$ if
\begin{enumerate}
\item 
It is $H$-equivariant when we provide $Y$ with the trivial action.

\item
If $f':R\too Y'$ is another morphism satisfying (1), then there
is a unique morphism $g:Y \to Y'$ such that $f'=g\circ f$.
\end{enumerate}
\end{definition}

\begin{definition}[Good quotient]
A morphism $f:R\to Y$ of schemes is a good quotient for an 
action of an algebraic group $H$ 
on $R$ if
\begin{enumerate}
\item $f$ is surjective, affine and $H$-equivariant, 
when we provide $Y$ with the trivial action.

\item $f_*(\SO^H_R)=\SO_Y^{}$, where $\SO^H_R$ is the sheaf of
$H$-invariant functions on $R$.

\item If $Z$ is a closed $H$-invariant subset of $R$, then $p(Z)$ is
closed in $Y$. Furthermore, if $Z_1$ and $Z_2$ are two 
closed $H$-invariant subsets of $R$ with $Z_1\cap Z_2=\emptyset$, 
then $f(Z_1)\cap f(Z_2)=\emptyset$
\end{enumerate}
\end{definition}

\begin{definition}[Geometric quotient]
A geometric quotient $f:R\to Y$ is a good quotient 
such that $f(x_1)=f(x_2)$ if and only if
the orbit of $x_1$ is equal to the orbit of $x_2$.
\end{definition}

Clearly, geometric quotients are good quotients, and 
these are categorical quotients.
Assume that $R$ is projective, and the action of $H$ on $R$ has
a linearization on an ample line bundle $\SO_R(1)$.
A closed point $y\in R$ is called GIT-semistable if, for some
$m>0$, there is an $H$-invariant section $s$ of $\SO_R(m)$
such that $s(y)\neq 0$. If, moreover, every orbit of 
$H$ in $R_s=\{x\in R| s(x)\neq 0\}$ is closed and of
the same dimension as $H$, then $y$ is called 
a GIT-stable point.
We will use the following characterization in 
\cite{Mu1} of GIT-(semi)stability: let 
$\lambda:\CC^*\to H$ be a one-parameter subgroup,
and $y\in R$. Then $\lim_{t\to 0}\lambda(t)\cdot y=y_0$
exists, and $y_0$ is fixed by $\lambda$. Let
$t\mapsto t^a$ be the character by which $\lambda$
acts on the fiber of $\SO_R(1)$. Defining 
$\mu(y,\lambda)=a$, Mumford proves that 
$y$ is GIT-(semi)stable if and
only if, for all one-parameter subgroups, 
it is $\mu(y,\lambda)(\leq)0$.

\begin{proposition}
Let $R^{ss}$ (respectively $R^{s}$) be the open subset of
GIT-semistable points (respectively GIT-stable). 
Then there is a good quotient $R^{ss}\to R\gitq H$,
and the restriction $R^{s}\to R^{s}\gitq H$
is a geometric quotient.
Furthermore, $R\gitq H$ is projective and $R^{s}\gitq H$
is an open subset.
\end{proposition}

\begin{definition}
\label{corepresents}
A scheme $Y$ corepresents a functor $F:(\Sch/\CC)
\to (\Sets)$ if
\begin{enumerate}

\item There exists a natural transformation $f:F \to \underline{Y}$
(where $\underline{Y}={\rm Mor}(\cdot,Y)$ 
is the functor of points represented by $Y$).

\item For every scheme $Y'$ and natural transformation 
$f':F \to \underline{Y'}$, there exists a unique 
$g:\underline{Y}\to \underline{Y'}$ 
such that $f'$ factors through $f$.

\end{enumerate}
\end{definition}

\begin{remark}
\label{corep}
\textup{
Let $R$ be a universal space with group $H$ for $F$, and
let $f:R\to Y$ be a categorical quotient. It follows from
the definitions that $Y$ corepresents
$F$.
}
\end{remark}

\begin{proposition}
\label{410}
Let $\SP_R=(P_R,E_R,\psi_R)$ 
be a universal family with group $H$ for the functor 
$\wt{F}^\tau_\group$. Let $\rho:\groupp \to \group$ be a homomorphism 
of groups,
such that the center $Z_\groupp$ of $\groupp$ is applied 
to the center $Z_\group$ of $\group$ 
and the induced homomorphism
$$
\Lie(\groupp/Z_\groupp)\too \Lie(\group/Z_\group)
$$
is an isomorphism.
Assume that the functor $\wt\Gamma(\rho,P_R)$ is represented
by a scheme $M$. Then
\begin{enumerate} 
\item
There is a natural action of $H$ on $M$, making
it a universal space with group $M$ for the functor
$\wt{F}^\tau_\groupp$.

\item
Moreover, if $\rho$ is injective (so that $\Gamma(\rho,P_R)$ 
itself is representable by $M$), then the action of $H$ lifts
to the family $\SP_M$ given by $\Gamma(\rho,P_R)$, and
then $\SP_M$ becomes a universal family with group $H$ for
the functor $\wt{F}^\tau_\groupp$
\end{enumerate}
\end{proposition}

\begin{proof}
Analogous to \cite[lemma 4.10]{Ra3}.
\end{proof}

\section{Construction of $R$ and $R_1$}
\label{sec1}

Given a principal $G$-bundle, we obtain a pair 
$(E,\varphi:E\otimes E\to E)$, where $E=P(\fgp)$ is the
vector bundle associated to the adjoint representation of
$G$ on the semisimple part $\fgp$ of the Lie algebra of $G$,
and $\varphi$ is given by the Lie algebra structure.
To obtain a projective moduli space we have to allow $E$ 
to become a torsion free sheaf. For technical reasons, 
when $E$ is not locally free, we make $\varphi$ take
values in $E^{\vee\vee}$.

The first step to construct the moduli space is the construction
of a scheme parameterizing semistable based $\fgp$-sheaves, 
i.e. triples
$(q:V\otimes \SO_X(-m)\surj E, E,
\varphi:E\otimes E\to E^{\vee\vee})$, where $V$ is a fixed vector
space, $m$ is a suitable large integer depending only on the 
numerical invariants, and $(E,\varphi)$ is
a semistable $\fgp$-sheaf. 

We have already seen that a $\fgp$-sheaf can be described
as a \tensor\ in the sense of \cite{G-S1}, where a notion
of (semi)stability for \tensors\ is given, depending on a polynomial
$\delta$ of degree at most $n -1$ and positive leading coefficient.
\textit{
In this article we will always assume that $\delta$ has degree 
$n-1$.}
Now we will prove, after some lemmas,
that the (semi)stability of the $\fgp$-sheaf coincides with
the $\delta$-(semi)stability of the corresponding \tensor
(in particular for the tensors associated to $\fgp$-sheaves,
its $\delta$-(semi)stability does not depend on $\delta$,
as long as $\deg(\delta)=n-1$),
so that we can apply the results of \cite{G-S1}, and the
moduli space of semistable $\fgp$-sheaves is a subscheme
of the moduli space of $\delta$-semistable tensors.

Given a $\fgp$-sheaf $(E,\varphi)$ and a balanced 
filtration
$E_{\lambda_\bullet}$, define
\begin{eqnarray}
\label{muvarphi}
\mu(\varphi,E_{\lambda_\bullet}) &=& 
\min \big\{
\lambda_i + \lambda_j - \lambda_k \;:\; 0\neq 
\varphi:E_{\lambda_i}\otimes E_{\lambda_j} \too 
E^{\vee\vee}/E^{\quad\vee\vee}_{\lambda_{k-1}} \big\}\\
\nonumber
& =& 
\min \big\{
\lambda_i + \lambda_j - \lambda_k \;:\;  
[E_{\lambda_i}, E_{\lambda_j}] \nsubseteq 
E^{\quad\vee\vee}_{\lambda_{k-1}} \big\}
\end{eqnarray}

\begin{lemma}
\label{samemu}
\textup{
If $(E,\phi)$ is the associated tensor, then 
$\mu(\varphi,E_{\lambda_\bullet})$ in
(\ref{muvarphi}) is equal
to $\mu_{\rm tens}(\phi,E_{\lambda_\bullet})$ in
(\ref{stabilitytensor}).}
\end{lemma}

\begin{proof}
For a general $x\in X$ let $e_1,\ldots,e_r$ be a basis 
adapted to
the flag $E_{\lambda_\bullet}(x)$, thus giving a splitting
$E(x)=\oplus E^{\lambda_i}(x)$. 
Writing $r^{\lambda_i}= \dim E^{\lambda_i}(x)$, 
\begin{eqnarray*}
&& \mu_{\rm tens}(\phi,E_{\lambda_\bullet}) =  
\\&&  =\min \big\{
\lambda_i+\lambda_j + \lambda_1 r^{\lambda_1} + \cdots +
\lambda_k (r^{\lambda_k} -1)+
\cdots + \lambda_{t+1} r^{\lambda_{t+1}} \; : 
\\&& \rule{1cm}{0cm} e_1 {}_\wedge 
e_2 {}_\wedge \ldots {}_\wedge e_{k'-1} 
{}_\wedge \varphi_x(e_{i'}\otimes e_{j'})
{}_\wedge e_{k'+1} {}_\wedge \ldots {}_\wedge e_r\neq 0 
\;\; \text{for some} 
\\&& \rule{2cm}{0cm}
e_{i'}\in E^{\lambda_i}(x), \, e_{j'}\in E^{\lambda_j}(x), \,
1\leq k'\leq r \big\} =
\\&& 
= \min \big\{
\lambda_i+\lambda_j-\lambda_k \; : \; 
\varphi_x(E^{\lambda_i}(x),E^{\lambda_j}(x)) \nsubseteq
E_{\lambda_{k-1}}(x) 
\\&& \rule{1cm}{0cm} \text{and}\;
\varphi_x(E^{\lambda_i}(x),E^{\lambda_j}(x)) \subseteq
E_{\lambda_{k}}(x) \big\} =  \\&& 
= \min \big\{
\lambda_i+\lambda_j-\lambda_k \; : \; 
[E_{\lambda_i},E_{\lambda_j}] \nsubseteq
E_{\lambda_{k-1}}^{\quad\vee\vee}
\\&& \rule{1cm}{0cm} \text{and}\;
[E_{\lambda_i},E_{\lambda_j}] \subseteq
E_{\lambda_{k}}^{\quad\vee\vee}   \} =  \\&&
= \mu(\varphi,E_{\lambda_\bullet})
\end{eqnarray*}
\end{proof}

We will need the following result, due to Ramanathan
\cite[lemma 5.5.1]{Ra3}, whose proof we recall for 
convenience of the reader.

\begin{lemma}
\label{ra551}
Let $W$ be a vector space, and let $p\in
\PP(W^\vee\otimes W^\vee \otimes W)$ be the point
corresponding to a  Lie algebra structure on $W$.
If the Lie algebra is semisimple, this point is GIT-semistable
for the natural action of $\slw$ and linearization in  
$\SO(1)$ on $\PP(W^\vee\otimes W^\vee \otimes W)$ 
\end{lemma}

\begin{proof}
Define the $\slw$-equivariant homomorphism
\begin{eqnarray*}
g:\,(W^\vee\otimes W^\vee \otimes W)=\Hom(W,\End W) &\too& 
(W\otimes W)^\vee\\
f & \mapsto & g(f)(\cdot\otimes \cdot)=\tr (f(\cdot)\circ f(\cdot))
\end{eqnarray*}
Choose an arbitrary linear space isomorphism between $W$ and $W^\vee$.
This gives an isomorphism $(W\otimes W)^\vee\isom\End(W)$. 
Define the determinant map $\det:(W\otimes W)^\vee\isom\End(W)\to
\CC$. Then $\det\circ g$ is an $\slw$-equivariant polynomial on
$W^\vee\otimes W^\vee \otimes W$ and it is nonzero when evaluated on 
the point $f$ corresponding to a semisimple Lie algebra, because it is
the determinant of the Killing form. Hence this point is GIT-semistable.
\end{proof}

\begin{lemma}
\label{mu0}
If $\varphi$ is a $\fgp$-sheaf,
then $\mu(\varphi,E_{\lambda_\bullet})\leq 0$ for any
balanced filtration $E_{\lambda_\bullet}$,
and $\mu(\varphi,E_{\lambda_\bullet})=0$ if and only if
it is an algebra filtration.
\end{lemma}

\begin{proof}
Since $E^{\vee\vee}$ is torsion free, the formula (\ref{muvarphi}) is
equivalent to
\begin{equation}
\label{equivmuvarphi}
\mu(\varphi,E_{\lambda_\bullet}) =
\min \big\{
\lambda_i + \lambda_j - \lambda_k \;:\;  
[E_{\lambda_i}(x), E_{\lambda_j}(x)] \nsubseteq 
E^{\quad\vee\vee}_{\lambda_{k-1}}(x) \big\}
\end{equation}
where $x$ is a general point of $X$, so that 
$E_{\lambda_\bullet}$ is a vector bundle filtration near $x$.
Fixing a Lie algebra isomorphism between the fiber $E(x)$ and $\fgp$, 
the filtration $E_{\lambda_\bullet}$ induces a filtration on 
$\fgp$. Consider a vector space splitting $\fgp=\oplus \fgp^{\lambda_i}$
of this filtration and a basis $e_l$ of $\fgp$ such that
$e_l\in \fgp^{i(l)}$, in order to define a monoparametric subgroup
of ${\rm SL}(\fgp)$ given by $e_l \mapsto t^{\lambda_{i(l)}} e_l$
for all $t\in \CC^*$ (cfr. notation introduced for definition
\ref{defbalancedfiltration}). The Lie algebra structure on $\fgp$ gives
a point $\langle \varphi_\fgp\rangle \in\PP(\fgp^\vee\otimes \fgp^\vee \otimes \fgp)$.
Let $a_{lm}^n$ be the homogeneous coordinates of this point,
i.e. $[e_l,e_m]=\sum_n a_{lm}^n e_n$. The monoparametric subgroup
acts as $t^{\lambda_{i(l)}+\lambda_{i(m)}-\lambda_{i(n)}} a_{lm}^n$
on the coordinates $a_{lm}^n$.
Hence (\ref{equivmuvarphi}) is equivalent to 
$$
\mu(\varphi,E_{\lambda_\bullet}) = 
\min \big\{
\lambda_{i(l)} + \lambda_{i(m)} - \lambda_{i(n)} \;:\; 
a^n_{lm}\neq 0 \big\}
$$
By lemma \ref{ra551}, the point $\varphi_\fgp$ 
is semistable under the ${\rm SL}(\fgp)$ action because it
corresponds to a semisimple Lie algebra, hence 
$\mu(\varphi,E_{\lambda_\bullet})\leq 0$.

Now assume that $\mu(\varphi,E_{\lambda_\bullet})= 0$. Then
it follows from (\ref{muvarphi}) that
$$
[E_{\lambda_i},E_{\lambda_j}]\subseteq E^{\quad\vee\vee}_{\lambda_{k-1}}
$$
for all $\lambda_i$, $\lambda_j$, $\lambda_k$ with
$\lambda_i+\lambda_j-\lambda_k<0$, i.e. $E_{\lambda_\bullet}$ is
an algebra filtration of $E$.

Conversely, if $E_{\lambda_\bullet}$ is an
algebra filtration of $E$, then $\mu(\varphi,E_{\lambda_\bullet})= 0$,
because if $\mu(\varphi,E_{\lambda_\bullet})< 0$, then for some
triple $(\lambda_i,\lambda_j,\lambda_k)$ with
$\lambda_i+\lambda_j<\lambda_k$ it is 
$[E_{\lambda_i},E_{\lambda_j}] \nsubseteq
E^{\quad\vee\vee}_{\lambda_{k-1}}$,
contradicting that $E_{\lambda_\bullet}$ is an
algebra filtration.

\end{proof}

\begin{lemma}
\label{slopess}
Let $(E,\varphi:E\otimes E \to E^{\vee\vee})$ 
be a $\fgp$-sheaf, and let
$(E,\phi:E^{\otimes r+1}\to \SO_X)$ 
be the associated Lie tensor. Assume
that one of the following conditions is satisfied
\begin{enumerate}
\item $(E,\varphi)$ is a semistable $\fgp$-sheaf (definition 
\ref{stabilitylie})

\item $(E,\phi)$ is a $\delta$-semistable tensor (definition 
\ref{stabilitytensor})
\end{enumerate}
Then $E$ is a Mumford semistable sheaf.
\end{lemma}

\begin{proof}
Assume $E$ is not Mumford semistable. 
Consider its Harder-Narasimhan
filtration, i.e. the filtration
\begin{equation}
\label{HNfilt}
0 = E_0 \subsetneq E_1 \subsetneq E_2 \subsetneq \cdots 
\subsetneq E_t \subsetneq E_{t+1}=E
\end{equation}
such that $E^i=E_i/E_{i-1}$ is Mumford semistable for all $i=1,\ldots,t+1$, and
\begin{equation}
\label{decreasing}
\mu_{\max}(E):=\mu(E^1) > \mu(E^2) > \cdots > \mu(E^{t+1})
=:\mu_{\min}(E),
\end{equation}
where $\mu(F):=\deg(F)/\rk(F)$ denotes the slope of a sheaf $F$.
Define
\begin{equation}
\label{weights}
\lambda_i=-r! \mu(E^i)
\end{equation}
(the  factor $r!$ is used to make sure that $\lambda_i$ is integer).
Changing the indexes $i$ by $\lambda_i$, the Harder-Narasimhan
filtration becomes
$$
0  \subsetneq E_{\lambda_1} \subsetneq 
E_{\lambda_2} \subsetneq \cdots \subsetneq E_{\lambda_t} \subsetneq
E_{\lambda_{t+1}}=E
$$
Since $\deg(E)=0$ (by lemma \ref{lemma0}), it follows that this
filtration is balanced (definition \ref{defbalancedfiltration}).
Now we will check that it is an algebra filtration. Given a triple
$(\lambda_i,\lambda_j,\lambda_k)$, with
$\lambda_i+\lambda_j<\lambda_k$,
we have to show that 
$$
[E_{\lambda_i},E_{\lambda_j}] \subset E_{\lambda_{k-1}}^{\quad\vee\vee}.
$$
Let $k'$ be the minimum integer for which
$$
[E_{\lambda_i},E_{\lambda_j}] \subset E_{\lambda_{k'-1}}^{\quad\vee\vee}.
$$
We have to show that $k'\leq k$.
By definition of $k'$, the following composition is nonzero
$$
E_{\lambda_i}\otimes E_{\lambda_j} \stackrel{[\cdot,\cdot]}{\too}
E_{\lambda_{k'-1}}^{\quad\vee\vee} \too 
E_{\lambda_{k'-1}}^{\quad\vee\vee}/E_{\lambda_{k'-2}}^{\quad\vee\vee}
$$
It is well known that, if a homomorphism $F_1\to F_2$ between two torsion
free sheaves is nonzero, then $\mu_{\min}(F_1)\leq \mu_{\max}(F_2)$,
hence 
\begin{equation}
\label{mus}
\mu_{\min}(E_{\lambda_i}\otimes E_{\lambda_j}) \leq
\mu_{\max}(E_{\lambda_{k'-1}}^{\quad\vee\vee}/
E_{\lambda_{k'-2}}^{\quad\vee\vee})
\end{equation}
Using (\ref{weights}) and the fact that 
$\mu_{\min}(E_{\lambda_1}\otimes E_{\lambda_2})=
\mu_{\min}(E_{\lambda_1})+\mu_{\min}(E_{\lambda_2})$ 
\cite[Prop. 2.9]{A-B}), 
the left hand side is 
$$
\mu_{\min}(E_{\lambda_i}\otimes E_{\lambda_j})=
\frac{-1}{r!}(\lambda_i+ \lambda_j)
$$
Since the quotient
$E_{\lambda_{k'-1}}^{\quad\vee\vee}/E_{\lambda_{k'-2}}^{\quad\vee\vee}$ 
is Mumford semistable, the right hand side is
$$
\mu_{\max}(E_{\lambda_{k'-1}}^{\quad\vee\vee}/E_{\lambda_{k'-2}}^{\quad\vee\vee})
=\mu(E_{\lambda_{k'-1}}^{\quad\vee\vee}/E_{\lambda_{k'-2}}^{\quad\vee\vee})
=\frac{-1}{r!}\lambda_{k'-1}
$$
Hence the inequality (\ref{mus}) becomes
$$
\lambda_i+ \lambda_j \geq \lambda_{k'-1},
$$
and then $\lambda_{k'-1}<\lambda_k$, hence $k'\leq k$, and we
conclude that $E_{\lambda_\bullet}$ is a balanced algebra
filtration. 

If we plot the points 
$(r_{\lambda_i},d_{\lambda_i})=(\rk E_{\lambda_i}, \deg E_{\lambda_i})$, 
$1\leq i\leq t+1$ in the plane $\ZZ\oplus \ZZ$ we get a polygon, called 
the Harder-Narasimhan polygon. Condition (\ref{decreasing})
means that this polygon is (strictly) convex. Since $d=0$ (and
$d_{\lambda_1}>0$), this implies that $d_{\lambda_i}>0$ for $1\leq i
\leq t$, and then 
\begin{equation}
\label{polyleading}
\sum_{i=1}^{t} r!\big(\mu(E^i)-\mu(E^{i+1})\big)(r d_{\lambda_i} - r_{\lambda_i} d) >0.
\end{equation}
Therefore
\begin{equation}
\label{poly}
\sum_{i=1}^{t} (\lambda_{i+1}-\lambda_i)
\big( r P_{E_{\lambda_i}} -r_{\lambda_i} P_E \big)\succ 0
\end{equation}
because the leading coefficient of (\ref{poly}) is
(\ref{polyleading}), and then $(E,\varphi)$ is not
semistable as a $\fgp$-sheaf.
Hence, if $(E,\varphi)$ is semistable, then 
$E$ is Mumford semistable.

Now, since the Harder-Narasimhan filtration 
(\ref{HNfilt}) of $E$ is
an algebra filtration, it is, by lemma \ref{mu0}, 
$\mu(\varphi,E_{\lambda_\bullet})=0$.
Now, lemma \ref{samemu} implies 
$\mu_{\rm tens}(\phi,E_{\lambda_\bullet})=0$,
hence 
$$
\sum_{i=1}^{t} (\lambda_{i+1}-\lambda_i)
\big( r P_{E_{\lambda_i}} -r_{\lambda_i} P_E \big) +
\mu(\phi,E_{\lambda_\bullet})
=
\sum_{i=1}^{t} (\lambda_{i+1}-\lambda_i)
\big( r P_{E_{\lambda_i}} -r_{\lambda_i} P_E \big)\succ 0
$$
and then 
$(E,\phi)$ is not $\delta$-semistable as a tensor.
Hence, if $(E,\phi)$ is $\delta$-semistable, it follows
that $E$ is Mumford semistable.

\end{proof}

\begin{proposition}
\label{samestable}
Let $(E,\varphi:E\otimes E\to E^{\vee\vee})$ be a $\fgp$-sheaf and 
let $(E,\phi:E^{\otimes r+1}\to \SO_X)$ be the associated \tensor.
The following conditions are equivalent
\begin{enumerate}
\item $(E,\phi)$ is a $\delta$-(semi)stable tensor 

\item $(E,\varphi)$ is a (semi)stable $\fgp$-sheaf

\end{enumerate}
\end{proposition}

\begin{proof}
Assume that $(E,\phi)$ is $\delta$-semistable. By lemma
\ref{slopess}, the sheaf $E$ is Mumford semistable.
Let  $E_{\lambda_\bullet}$ be a balanced algebra
filtration. Then
$\mu_{\rm tens}(\phi,E_{\lambda_\bullet})=\mu(\varphi,E_{\lambda_\bullet})=0$ 
(lemmas \ref{samemu} and \ref{mu0}),
hence inequality (\ref{stabtensor}) in definition 
\ref{stabilitytensor} becomes (\ref{stablie}) in definition
\ref{stabilityliealgebra}.

Conversely, assume that the $\fgp$-sheaf $(E,\varphi)$ 
is (semi)stable,
thus $E$ is again Mumford semistable, 
and consider a balanced filtration $E_{\lambda_\bullet}$ of
$E$. We have to show that
(\ref{stabtensor}) is satisfied.
If the filtration is an algebra
filtration, then $\mu(\varphi,E_{\lambda_\bullet})=0$ by lemma
\ref{mu0}, hence (\ref{stabtensor}) holds. 
If it is not an algebra
filtration, then $\mu(\varphi,E_{\lambda_\bullet})<0$ (again by
lemma \ref{mu0}). Since $E$ is Mumford semistable, 
it is $r d_{\lambda_i} -r_{\lambda_i} d\leq 0$ for all $i$.
Denote by $\tau/(n-1)!$ the coefficient of $t^{n-1}$ in
$\delta$. It is $\tau>0$ because $\deg \delta=n-1$.
Then the leading coefficient of the polynomial of 
(\ref{stabtensor}) becomes
$$
\Big(\sum_{i=1}^{t} (\lambda_{i+1}-\lambda_i)
\big( r d_{\lambda_i} -r_{\lambda_i} d \big)\Big)
+ \tau \mu(\varphi,E_{\lambda_\bullet})
< 0 ,
$$
and thus (\ref{stabtensor}) holds.

\end{proof}

Now, let us recall briefly how the moduli space
of tensors was constructed in \cite{G-S1}.
Start with a $\delta$-semistable tensor
$$
\phi:F ^{\otimes a} \too \SO_X
$$
with $\rk F=r$, Hilbert polynomial $P_F=P$ and $\det F\isom \SO_X$.
Let $m$ be a large integer (depending only on the polarization and 
numerical invariants of $F$)
and an isomorphism $g$  between 
$H^0(F(m))$
and a fixed
vector space $V$ of dimension $h^0(F(m))$.
This gives a quotient
$$
q: V\otimes \SO_X(-m)\too F
$$
and hence a point in  the Hilbert scheme $\SH$ 
of quotients of $V\otimes \SO_X(-m)$
with Hilbert polynomial $P$. 
Let $l>m$ be an integer, 
and $W=H^0(\SO_X(l-m))$. The  quotient $q$ 
induces homomorphisms
$$
\begin{array}{rccl}
q\,\,: & V\otimes \SO_X(l-m)&\surj& F(l) \\
q'\,:& V\otimes W&\to& H^0(F(l)) \\
q'':& \bigwedge{}^{P(l)}(V\otimes W) &\to& 
\bigwedge{}^{P(l)} H^0(F(l))\;\cong \; \CC
\end{array}
$$
If $l$ is large enough, these homomorphisms are surjective, and they 
give Grothendieck's embedding
$$
\SH \;\inj\; \PP\big(\bigwedge{}^{P(l)}(V^\vee\otimes W^\vee)\big).
$$
and hence a very 
ample line bundle $\SO_\SH(1)$ on $\SH$ (depending on
$m$ and $l$).
The isomorphism $g:V \stackrel{\isom}{\to} H^0(F(m))$ and
$\varphi$ induces a linear map
$$
\Phi: V^{\otimes a} \;\too\; 
H^0(F(m)^{\otimes a}) \;\too\;
H^0(\SO_X(am)) \; =:\; B,
$$
and hence the tensor $\varphi$ and the isomorphism
$g$ give a point in
$$
\PP\big(\bigwedge{}^{P(l)}(V^\vee\otimes W^\vee)\big) \times
\PP\big((V^{\otimes a})^\vee \otimes B\big) 
\;=\;
\PP \times \PP'
$$
Let $Z$ be the closure of the points associated to 
$\delta$-semistable tensors.
We give $Z$ a polarization $\SO_Z(1)$, by restricting
a polarization $\SO_{\PP\times\PP'}(b,b')$, where the ratio between
$b$ and $b'$ depends on the polynomial $\delta$ and
the integers $m$ and $l$
$$
\frac{b'}{b}=\frac{P(l)\delta(m)-\delta(l)P(m)}{P(m)-a\delta(m)}
$$
There is a tautological family of tensors 
parametrized by $Z$
\begin{equation}
\label{tautmoduli1}
\phi_Z^{}: F_Z^{\otimes r+1}\too p^*_{\PP'} \SO_{\PP'}^{}(1),
\end{equation}

The scheme $Z$ has an open
dense set $Z^{ss}$ representing the sheafification of the functor
\begin{equation}
\label{functorF}
F^b:({\rm Sch}/\CC) \too ({\rm Sets})
\end{equation}
associating to a scheme $S$ the set of
equivalence classes of families of $\delta$-semistable ``based'' 
\tensors\ 
$$
\big(q_S^{}: V\otimes \SO_{X\times S}(-m)^{} {\to} F_S^{} ,\;
F_S^{},\; \phi_S^{}:F_S^{\otimes a}\to p^*_S N,\;  N \big)
$$
where $q_S$ is a surjection inducing an isomorphism 
$$
g_S^{}=p^{}_{S*}(q_S^{}(m)):V\otimes \SO_S^{} \to p^{}_{S*}(F_S^{}(m))
$$
and $(F_S,\phi_S, N)$ is a family of
$\delta$-semistable tensors (definition \ref{stabilitytensor})
with fixed rank $r$, Hilbert polynomial $P$ and trivial 
determinant. In particular, 
\begin{equation}
\label{pullbackL}
\det(F_S)\cong p_S^* L,
\end{equation}
where $L$ is a line bundle on $S$.
From now on, we will assume $a=r+1$, where $r$ is 
the rank of $F$.

\begin{proposition}
\label{r}
There is a closed subscheme $R$ of $Z^{ss}$ 
representing the sheafification 
$\wt{F}^b_{\rm Lie}$ of the subfunctor
of (\ref{functorF}) 
\begin{eqnarray}
\label{functorFlie}
F^b_{\rm Lie}: ({\rm Sch}/\CC) & \too &({\rm Sets})\\
\nonumber
 S &\longmapsto & F^b_{\rm Lie}(S)\subset F^b(S)
\end{eqnarray}
where $F^b_{\rm Lie}(S)\subset F^b(S)$ is the subset of 
families of based $\delta$-semistable Lie tensors.

A point of the closure $\overline{R}$ of $R$ in $Z$ 
is GIT-(semi)stable with respect to the natural $\slv$-action
and linearization on $\SO_{\overline R}(1)=\SO_Z(1)|_{\overline R}$
(see \cite{G-S1}) if and only if
the corresponding \tensor\  is $\delta$-(semi)stable 
and $q$ induces an isomorphism $V\cong H^0(E(m))$. 
In particular
the open subset of semistable points of $\overline{R}$ is $R$.
\end{proposition}

\begin{proof}
Let $(q_{Z^{ss}}^{}, F^{}_{Z^{ss}}, 
\phi^{}_{Z^{ss}}:F_{Z^{ss}}^{\otimes r+1}\to p^*_{Z^{ss}} N, 
N)$ 
be the tautological family on $Z^{ss}$ coming from 
(\ref{tautmoduli1}). 
For each pair $(i,j)$ with $1\leq i<j\leq r+1$, let
$$
\sigma^{}_{ij}(\phi^{}_{Z^{ss}}): F_{Z^{ss}}^{\otimes r+1}
\too N
$$ 
be the homomorphism obtained from $\phi_{Z^{ss}}$ by interchanging
the factors $i$ and $j$. Let $Z_{ij}\subset Z^{ss}$ be the 
zero subscheme defined by 
$\phi^{}_{Z^{ss}} +\sigma^{}_{ij}(\phi^{}_{Z^{ss}})$,
using lemma \ref{zeroes}. Finally, define
$$
Z_{\rm skew}= \bigcap_{3\leq i<j\leq r+1} Z_{ij}.
$$
{}From the universal property of $Z_{ij}$ (lemma
\ref{zeroes}) it follows that, for a family satisfying
condition (1) of definition \ref{deftensor}, 
the classifying morphism into $Z^{ss}$ factors 
through $Z_{\sk}$. Furthermore, the restriction of the tautological
family to $Z_{\sk}$ satisfies condition (1), hence by corollary
\ref{cor12} we have a family parametrized by $Z_{\sk}$
\begin{equation}
\label{universal}
(q_{Z_{\sk}}^{} , F_{Z_{\sk}}^{}, \varphi_{Z_{\sk}}^{}: 
F_{Z_{\sk}}^{}\otimes F_{Z_{\sk}}^{} \too F_{Z_{\sk}}^{\quad\vee\vee}\otimes
p^*_{Z_\sk} N, N)
\end{equation}
The closed subscheme (``antisymmetric locus'') 
$Z_{\asym}\subset Z_{\sk}$ 
is defined as the zero subscheme
of $\varphi_{Z_{\sk}} + \sigma_{12}(\varphi_{Z_{\sk}})$ given
by lemma \ref{zeroes}. It follows that if a family
satisfies conditions (1) and (2) of definition \ref{deftensor}, 
then the classifying morphism factors
through $Z_{\asym}$, and furthermore the restriction of the 
tautological family to $Z_{\asym}$ satisfies conditions (1) and (2).

Let $J$ be the homomorphism defined as in (\ref{jacobin}), using 
the tautological family parametrized by $Z_{\asym}$.
Note that this homomorphism is zero if and only if the associated
homomorphism (lemma \ref{lemma1})
$$
J':F_{Z_{\asym}}^{}\otimes F_{Z_{\asym}}^{}\otimes F_{Z_{\asym}}^{} 
\otimes F_{Z_{\asym}}^\vee \too p^*_{Z_{\asym}} N^2
$$
is zero. Finally, define the closed subscheme $R\subset Z_{\asym}$ as 
the zero subscheme of $J'$ given in lemma \ref{zeroes}.
It follows that if a family satisfies conditions (1) to (3)
of definition \ref{deftensor},
then the classifying morphism will factor
through $R$, and furthermore the restriction of the 
tautological family to $R$ satisfies conditions (1) to (3).

The criteria for stability follows from \cite{G-S1}.
\end{proof}

Recall that a $\fgp$-sheaf is (semi)stable if and only if
the associated Lie tensor is $\delta$-semistable
(proposition \ref{samestable}).
\begin{proposition}
\label{r1}
There is a subscheme $R_1\subset R$ representing 
the sheafification $\wt{F}^{b}_{\fgp}$ of the subfunctor of (\ref{functorFlie})
\begin{eqnarray}
\label{functorFliegp}
F^{b}_{\fgp}: ({\rm Sch}/\CC) & \too &({\rm Sets})\\
\nonumber
 S &\longmapsto & F^{b}_{\fgp}(S)\subset F^b_{\rm Lie}(S)
\end{eqnarray}
where $F^{b}_{\fgp}(S)\subset F^{b}_{\rm Lie}(S)$ 
is the subset of $S$-families of based $\delta$-semistable 
Lie tensors such that the homomorphism associated by
construction \ref{constr} provides a family of based 
semistable $\fgp$-sheaves
with fixed numerical invariants $\tau$.

Furthermore, $R_1$ is a union of connected components of $R$,
hence the inclusion $R_1 \inj R$ is proper.
\end{proposition}

\begin{proof}
Consider the tautological family parametrized by $R$
$$
(q_R^{}, F_R^{}, \phi_R^{}: F_R^{\otimes r+1} \too p^*_R N, N)
$$
and the associated family obtained as in construction \ref{constr}
\begin{equation}
\label{rr1}
( q_R^{}, E_R^{}, \varphi_R^{}: E_R^{}\otimes E_R^{} \to E^{\vee\vee}_R)
\end{equation}
Let $\kappa$ be the Killing form (definition \ref{defkilling})
$$
\kappa: E_R\otimes E_R \too \SO_{X\times R}.
$$
This induces a homomorphism $\det \kappa':\det E_R \to \det E_R^\vee $.
Recall from (\ref{pullbackL}) that 
$\det(F_R)$ is the pullback of a line
bundle from $R$, hence the same holds for $\det(E_R)$, 
and then $\det \kappa'$ is constant along the fibers
of $\pi:X\times R\to R$. Hence $\det \kappa'$ is nonzero on an open set
of the form $X\times W$, where $W\subset R$ is an open set.

A point $(q,E,\varphi)\in R$ belongs to $W$ if and only if
for all $x\in U_E$ the Lie algebra
$(E(x),\varphi(x))$ is semisimple, because the Killing form is
nondegenerate if and only if the Lie algebra is semisimple.

Now we show that the open set $W$ is in fact equal to $R$.
Let $(q,E,\varphi:E\otimes E \to E^{\vee\vee})$ be a 
based algebra sheaf corresponding to a point in $R\setminus W$.
Then its Killing form $\kappa:E\otimes E\to \SO_X$ is degenerate.
Let $E_1$ be the kernel of the homomorphism induced by $\kappa$
$$
0 \too E_1 \too E \too E^\vee.
$$
By lemma \ref{slopess}, $E$ is Mumford semistable,
thus $E^\vee$ is Mumford semistable, and being
both of degree $0$, the sheaf 
$E_1$ is also of degree $0$ and Mumford semistable. 
Note that $E_1$ is a solvable ideal of $E$, i.e. 
the fibers of $E_1$ are solvable ideals
of the fibers of $E$ (over points where both sheaves are locally free)
\cite[proof of Th. 2.1 in chp VI]{Se2}. 
Since $E_1\otimes E_1$ (modulo torsion) and $E_1^{\vee\vee}$ are 
Mumford
semistable of degree zero, the image $E_2'=[E_1,E_1]$ of the
Lie bracket
homomorphism $\varphi:E_1\otimes E_1 \to E_1^{\vee\vee}$, 
is a Mumford semistable subsheaf of $E^{\vee\vee}$
of degree zero. Define $E_2=E_2'\cap E$. It
is a Mumford semistable subsheaf of $E$ of degree zero.
Similarly $E_3'=[F_2,F_2]$, $E_3$, etc... are all Mumford 
semistable
sheaves of degree zero. Since $E_1$ is solvable, we
arrive eventually to a non-zero sheaf $E'$ of degree zero, 
which is an abelian ideal of $E$.

We claim that the balanced filtration 
$E_{\lambda_1}= E'\subset E_{\lambda_2} = E$
with $\lambda_1=\rk{E'}-r$ and $\lambda_2=\rk E'$
contradicts the $\delta$-semistability of
the tensor $(E,\phi)$ associated to $(E,\phi)$ by
construction \ref{constr}.

To prove this we  
need to calculate 
$\mu_\tens(\phi,E_{\lambda_\bullet})$ (cfr. (\ref{muE})).
By lemma \ref{samemu} this is equal to
$\mu(\varphi,E_{\lambda_\bullet})$ 
(cfr. (\ref{muvarphi})). We need to estimate which triples 
$(i,j,k)$ are relevant to calculate the minimum, 
i.e. for which triples it is $[E_{\lambda_i}, E_{\lambda_j}] \nsubseteq 
E^{\quad\vee\vee}_{\lambda_{k-1}}$. 
Since $E'$ is abelian, $[E',E']=0$, so $(1,1,k)$ is not relevant.
Since $E'$ is an ideal, we have $[E',E]\subset E'{}^{\vee\vee}$. 
If $E'$ is 
in the center, then this bracket is zero, hence $(1,2,k)$ is 
not relevant.
On the other hand, if $E'$ is not in the center, then 
$[E',E]\neq 0$, hence $(1,2,1)$ is relevant, 
and the associated weight is $\lambda_1+\lambda_2-\lambda_1=
\rk(E')>0$. Since $E$ is not abelian, it is $[E,E]\neq 0$.
There are two possibilities: if $[E,E]\subset E'{}^{\vee\vee}$, then 
$(2,2,1)$ is relevant and
$\lambda_2+\lambda_2-\lambda_1=\rk(E')+\rk(E)>0$.
Otherwise $(2,2,2)$ is relevant, and 
$\lambda_2+\lambda_2-\lambda_2=
\rk(E')>0$. 
Summing up, we obtain
$$
\mu(\varphi, E_{\lambda_\bullet})>0.
$$
Since $\deg(E')=\deg(E)=0$, the leading coefficient of 
$$
\big(r P_{E'} - \rk(E') P_E\big) + \mu(\varphi, E_{\lambda_\bullet}) \delta
$$
is positive,
hence $(E,\phi)$ is not $\delta$-semistable (and by proposition
\ref{samestable}, $(E,\varphi)$ is not semistable), contradicting 
the assumption, so we have proved that $W=R$.

Now assume that we have two based $\fgp$-sheaves 
$(q,E,\varphi)$ and $(q',E',\varphi')$ belonging to the same connected
component of $R$, and $x\in U_E$,  $x'\in U_{E'}$. Then we have
$$
(E(x),\varphi(x))\cong (E'(x'),\varphi'(x'))
$$
as Lie algebras, because of the well known 
rigidity of semisimple Lie algebras 
(see \cite{Ri}, for instance).
Hence $R_1$ is the union of the connected components
of $R$ with $(E(x),\varphi(x))\cong \fgp$.

\end{proof}

We will denote by $\SE_{R_1}$ the tautological family of
$\fgp$-sheaves parametrized
by $R_1$ obtained by restricting (\ref{rr1})
\begin{equation}
\label{family1}
\SE_{R_1}=(E_{R_1}, \varphi_{R_1})
\end{equation}

Giving a family of (semi)stable $\fgp$-sheaves is equivalent
to giving a family of (semi)stable principal $\autfgp$-sheaves. 
By lemma \ref{stabilitylie}, 
the (semi)stability conditions for a $\fgp$-sheaf and
the corresponding principal $\autfgp$-sheaf coincide, hence 
$(E_{R_1},\varphi_{R_1})$ can be seen as 
a family of semistable
principal $\autfgp$-sheaves.

Recall that $\SH$ is
the Hilbert scheme classifying quotients $V\otimes \SO_X(-m)\to F$
(of fixed Chern classes), 
$\PP'=\PP\big((V^{\otimes r+1})^\vee \otimes 
H^0(\SO_X((r+1)m))\big)$ and, by construction \ref{constr}, we have 
$E_{R_1}=F_{R_1}\otimes \det F_{R_1} \otimes p^* \SO_{\PP'}(-1),$
where $F_{R_1}$ is the restriction of (\ref{tautmoduli1}) to $R_1$, and 
$p$ is 
$$
p:R_1 \inj \PP \times \PP' \to \PP'
$$
Let $\tau:V\otimes \SO_\glv \to V\otimes \SO_\glv$ be the 
universal automorphism. Let $\pi_{\glv}$, 
$\pi_{R_1}$ be the projections
to the two factors of $\glv\times R_1$.
The group $\glv$ acts on $R_1$, and this action 
lifts to $F_{R_1}$ (\cite[\S 4.3 pg. 90]{H-L}) 
and $p^* \SO_{\PP'}(1)$, giving 
isomorphisms $(\Lambda,\SB)$
\begin{equation}
\label{tenact}
\xymatrix{
{V\otimes \SO_{X\times R_1}(-m)} 
\ar[d]_{\pi_{\glv}^*\tau} \ar@{>>}[r]^-{\overline{\sigma}^*q^{}_{R_1}}  
& {\overline{\sigma}^*F^{}_{R_1}} \ar[d]^{\Lambda}_{\isom}\\
{V\otimes \SO_{X\times R_1}(-m)} 
\ar@{>>}[r]^-{\overline{\pi^{}_{R_1}}^*q^{}_{R_1}}  
& {\overline{\pi^{}_{R_1}}^*F^{}_{R_1}} 
}
\qquad
\xymatrix{
{\overline{\sigma}^* F_{R_1}^{\otimes r+1}}
\ar[d]^{\Lambda{}^{\otimes r+1}}_{\isom}
\ar[r]^-{\overline{\sigma}^*\phi^{}_{R_1}}&
 {\overline{\sigma}^* N} 
\ar[d]^{\SB}_{\isom} \\
{\overline{\pi^{}_{R_1}}^* F_{R_1}^{\otimes r+1}}  \ar[r]^
-{\overline{\pi^{}_{R_1}}^*\phi^{}_{R_1}}&
 {\overline{\pi^{}_{R_1}}^* N}
}
\end{equation}
between the pullbacks of the family of Lie tensors
$(F_{R_1},\phi_{R_1})$ by the action $\sigma:\glv\times R_1\to R_1$
and the projection $\pi^{}_{R_1}$ to the second factor. 

Since this action lifts to $F_{R_1}$ and $p^* \SO_{\PP'}(1)$, 
it also lifts to
$E_{R_1}$. An element $\lambda$ in the center of $\glv$
acts trivially on $R_1$, hence the action $\sigma$
factors through an action
action $m:\pglv\times R_1\to R_1$ of $\pglv$ on $R_1$. 
The element $\lambda$ acts as multiplication by $\lambda$
on $F_{R_1}$ and as multiplication by 
$\lambda^{-r-1}$ on $\SO_{\PP'}(-1)$, hence
it acts trivially on $E_{R_1}$. Therefore the action 
of $\glv$ on
$E_{R_1}$ factors through $\pglv$
\begin{equation}
\label{fgpact}
\xymatrix{
{\overline{m}^*E_{R_1}} \ar[d]^{\Lambda}_{\isom} \\
{\overline{p^{}_{R_1}}^*E_{R_1}}
}
\qquad
\xymatrix{
{\overline{m}^*E_{R_1}}\otimes {\overline{m}^*E_{R_1}}
\ar[rr]^-{\overline{m}^*\varphi_{R_1}^{}} 
\ar[d]^{\Lambda\otimes \Lambda} &&
{\overline{m}^*E_{R_1}^{\vee\vee}}
\ar[d]^{\Lambda^{\vee\vee}} \\
{\overline{p^{}_{R_1}}^*E_{R_1}}\otimes{\overline{p}_2^*E_{R_1}}
\ar[rr]^-{\overline{p^{}_{R_1}}^*\varphi_{R_1}^{}} &&
{\overline{p^{}_{R_1}}^*E_{R_1}^{\vee\vee}}
}
\end{equation}
where 
$p^{}_{R_1}$ is the projection of $\pglv\times R_1$ to 
the second factor.
\textit{
This gives a lift of the $\pglv$ action
on $R_1$ to the family $\SE_{R_1}$.
}

\begin{proposition}
\label{redr1}
With this action,
$(E_{R_1},\varphi_{R_1})$ 
becomes a 
universal family with group $\pglv$ for the
functor $\wt F^\tau_{\autlie(\fgp)}$ (cfr. remark \ref{fgpaut}).
\end{proposition}

\begin{proof}
Let $(E_S,\varphi_S)$ be a family of semistable $\fgp$-sheaves.
Shrink $S$ if necessary, so that $\det E_S\isom \SO_{X\times S}$. Using
this isomorphism and construction \ref{constr} we obtain
a family of $\delta$-semistable Lie tensors
$(E_S^{}, \phi_S^{}:E_S^{\otimes r+1}\to \SO_{X\times S}^{})$.
By proposition \ref{r1}, after shrinking $S$ if necessary,
there is a morphism
$f:S\to R_1$ such that the pullback 
$(\overline{f}^* E,\overline{f}^* \phi)$ of the 
family of Lie tensors 
parametrized by $R_1$ is isomorphic to $(E_S,\phi_S)$,
hence the families of $\fgp$-sheaves associated by
construction \ref{constr} to both of them are isomorphic. 

Now we are going to check the second condition in the definition of
universal family.
Let $t_1,t_2:S\to R_1$ be two morphisms, and let 
$\alpha:E_2\to E_1$
be an isomorphism between the two pullbacks $(E_1,\varphi_1)$
and $(E_2,\varphi_2)$ of $\SE_{R_1}$ under $t_1$ and $t_2$.
We have to find a morphism $h:S\to \pglv$ such that
$t_2=h[t_1]$ and $\overline{(h\times t_1)}^*
\Lambda=\alpha$. Since the question is local on $S$,
we may shrink $S$ when needed along the proof.

By pulling back the family $(F_{R_1},\phi_{R_1})$,
these morphisms also give two families of 
semistable
based Lie tensors
$(q_1, F_1,\phi_1)$ and $(q_2, F_2,\phi_2)$.
By definition of $E_{R_1}$, we have 
$E_i=F_i\otimes \det F_i \otimes N_i$, $i=1,2$.
After eventually shrinking $S$, there are isomorphisms $a_i:\det F_i \otimes
N_i\to \SO_{X\times S}$. Define $\alpha'$ by
$$
\xymatrix{
{E_2} \ar@{=}[r] \ar[d]_{\alpha}& {F_2\otimes \det F_2 \otimes N_2} 
\ar[r]^-{F_2\otimes a_2} 
& {F_2} \ar[d]^{\alpha'} \\
{E_1} \ar@{=}[r] & {F_1\otimes \det F_1 \otimes N_1} \ar[r]^-{F_1\otimes a_1}&
{F_1}
}
$$
and hence $\alpha=\alpha'\otimes (a_1^{-1}\circ a_2^{})$.
Given an isomorphism $\beta:N_2^{-1}\to N_1^{-1}$, we obtain
an isomorphism
$$
\alpha' \otimes \det\alpha' \otimes \beta:
E_2=F_2\otimes \det F_2 \otimes N_2 \too
E_1=F_1\otimes \det F_1 \otimes N_1 \, .
$$
Choose $\beta$ so that $\alpha'\otimes (a_1^{-1}\circ a_2^{})=
\alpha' \otimes \det\alpha' \otimes \beta$. Since 
$\alpha=\alpha' \otimes \det\alpha' \otimes \beta$, the
commutativity of 
$$
\xymatrix{
{E_2\otimes E_2} \ar[r]^-{\varphi^{}_2} \ar[d]_{\alpha\otimes \alpha} 
&{E_2^{\vee\vee}} \ar[d]^{\alpha^{\vee\vee}} \\
{E_1\otimes E_1} \ar[r]^-{\varphi^{}_1} & {E_1^{\vee\vee}}
}
$$
implies the commutativity of
$$
\xymatrix{
{F_2^{\otimes r+1}} \ar[r]^{\phi^{}_2} \ar[d]_{\alpha^{\prime\otimes r+1}} 
&{N_2} \ar[d]^{\beta} \\
{F_1^{\otimes r+1}} \ar[r]^{\phi^{}_1} &
{N_1}
}
$$
and hence the pair $(\alpha',\beta)$ gives an isomorphism
between $(F_1,\phi_1)$ and $(F_2,\phi_2)$.
Using the based Lie tensors $(q_1,F_1,\phi_1)$  
and $(q_2,F_2,\phi_2)$, let
$g_i=p^{}_{S*}(q_i(m))$, $i=1,2$, and define
the isomorphism $h'$
$$
\xymatrix{
{V\otimes \SO_S} \ar[d]_{h'} \ar[r]^{g^{}_2} & 
{p^{}_{S *}(F_2(m))}   \ar[d]^{\isom}_{p^{}_{S*}(\alpha'(m))} \\
{V\otimes \SO_S} \ar[r]^{g^{}_1} & {p^{}_{S *}(F_1(m))}  
}
$$
This isomorphism can be seen as a morphism $h':S \to \glv$. 
By construction, it is $t_2=h'[t_1]$, and $(\alpha',\beta)$
is the pullback of the isomorphism (\ref{tenact}) by $h'\times t_1$. 
Denote by $h:S \to \pglv$ the composition with the projection 
to $\pglv$.
Then we have $t_2=h[t_1]$, and $\alpha$ is the pullback
of the left arrow in (\ref{fgpact}) by $\overline{h\times t_1}$.

Finally, we have to check that these two properties determine
$h$ uniquely. Let $h_1, h_2:S\to \pglv$ be two such morphisms.
Define $h=h_1^{} h_2^{-1}$. Then $h[t_1]=t_1$, and the pullback
$\overline{h\times t_1}^* \Lambda$
is the identity automorphism.
Replacing  $S$ by an \'etale cover, we can lift $h$ to a morphism
$h':S\to \glv$, and this induces an automorphism 
$\alpha'=\overline{h\times t_1}^* \Lambda '$ of 
$F_{S}=\overline{t_1}^*F_{R_1}$
\begin{equation}
\label{prevp}
\xymatrix{
{V\otimes \SO_{X\times S}} 
\ar[d]_{h'} \ar@{>>}[r]^-{\overline{t_1}^*q^{}_{R_1}}  
& {F_{S}(m)} \ar[d]^{\alpha'}\\
{V\otimes \SO_{X\times S}} \ar@{>>}[r]^-{\overline{t_1}^*q^{}_{R_1}}  
& {F_{S}(m)} 
}
\end{equation}
Applying $p^{}_{S*}$ to (\ref{prevp}), we obtain
$$
\xymatrix{
{V\otimes \SO_{S}} 
\ar[d]_{p^{}_{S*} h'} \ar[r]^-{H^0(q^{}_1)}_-{\isom}  
& {p^{}_{S*} F_{S}(m)} \ar[d]^{p^{}_{S*} \alpha'}\\
{V\otimes \SO_{S}} \ar[r]^-{H^0(q^{}_1)}_-{\isom} 
& {p^{}_{S*} F_{S}(m)} 
}
$$
Since $\overline{h\times t_1}^* \Lambda=\id$, the automorphism
$\alpha'$ is a family of homotethies, i.e. $p^{}_{S*} \alpha'$ 
can be seen as a morphism $S\to \CC^*$, and considering the
previous diagram, $p^{}_{S*} h'$ can also be seen as
a morphism from $S$ to $\CC^*$, the center of $\glv$, 
hence $h$ is the identity morphism from $S$ to $\pglv$.

\end{proof}

\section{Construction of $R_2$}

Recall that 
\textit{all schemes considered are locally of finite
type over $\Spec \CC$.} 
In this section and the following we are going to make use
of the category of complex analytic spaces.
For a scheme $Y$, we denote by
$Y^\an$ the associated complex analytic space 
(\cite[XII]{SGA1}, \cite[App. B]{Ha}),
and given a morphism $f$ in the category
of schemes, we denote by $f^\an$ the corresponding
morphism in the category of analytic spaces.
Recall that the
underlying set of $Y^\an$ is the set of closed points
of $Y$, and it is endowed with the analytic topology.

\begin{lemma}
\label{generalpos}
Let $S$ be a scheme (not necessarily smooth).
Let $\SZ\subset X\times S$ be a closed subscheme with
$\codim_{\RR}(\SZ_s^\an,X^\an\times s)\geq m$ for all closed points 
$s\in S$, and $\SU\subset X\times S$ its complement.
Let $M$ be a real manifold with $\dim_\RR(M)\leq m-1$ and compact
boundary, and let
$$
f=(f_X,f_S):M\too X^\an\times S^\an
$$ 
be a continuous map such that the image of the 
boundary lies in $\SU^\an$. Then $f$ can be modified
by a homotopy, relative to its boundary to a continuous 
map $\wt f$ whose image lies in $\SU^\an$.
\end{lemma}

\begin{proof}
Consider the cartesian product (in the category of topological
spaces and continuous maps)
$$
\xymatrix{
{\SZ_{M}^\an} \ar[r] \ar[d] & {\SZ^\an} \ar[d]\\
{M} \ar[r]^{f_{S}} & {S^\an}
}
$$
The map $f$ factors as
$$
\xymatrix{
{M} \ar@(ur,ul)[rrrr]^{f}
\ar[rr]^-{(f_X,\id)} && {X^\an \times M} \ar[rr]^-{(\id,f_S)} && 
{X^\an\times S^\an}\\
&& {\SZ_M^\an} \ar@{^{(}->}[u] \ar[rr] && {\SZ^\an} \ar@{^{(}->}[u] 
}
$$
By hypothesis $\codim_{\RR}(\SZ_s^\an,X^\an\times s)\geq m$
for all $s\in S$, so $\codim_\RR(\SZ_M^\an,X^\an\times M)\geq m$,
and since $\dim_\RR(M)\leq m-1$
and $X^\an\times M$ is smooth, we can modify $(f_X,\id)$
homotopically, relative to its boundary,
to a map $\wt{f}_1$ whose image does not intersect 
$\SZ_M^\an$.
Then the image of $\wt f=(\id,f_S)\circ\wt{f}_1$ lies 
in $\SU$.
\end{proof}

\begin{lemma}
\label{extension}
For a scheme $S$, let $\SZ\subset X\times S$ be a closed subscheme
such that $\codim_\RR(\SZ^\an_s,X^\an\times s)\geq 4$
for all $s\in S$.
Let $\SU$ be the complement of $\SZ$, and
let $x\in \SU\subset X\times S$ be a closed point.
Then the inclusion  $i^\an:\SU^\an\inj X^\an\times S^\an$
induces an isomorphism of topological fundamental groups
$$
\pi_1(i^\an,x):\pi_1(\SU^\an,x)\stackrel{\cong}{\too} 
\pi_1( X^\an\times S^\an,x).
$$
\end{lemma}

\begin{proof}
To check that $\pi_1(i^\an)$ is injective, let $f:\sone\to \SU^\an$ be a 
continuous based loop (i.e. 
a continuous map from the unit interval $[0,1]$
sending $0$ and $1$ to the base point $x$)
mapping to zero in $\pi_1( X^\an\times S^\an,x)$.
So there is a continuous map $g$ fitting into a
commutative diagram
$$
\xymatrix{
{\sone} \ar[r]^f \ar@{^{(}->}[d] & {\SU^\an} \ar@{^{(}->}[r] & 
{X^\an\times S^\an} \ar@{=}[d]\\
{\DD} \ar[rr]^{g}  && 
{X^\an\times S^\an}
}
$$
where $\DD$ denotes the unit disk (whose boundary is $\sone$).
By lemma \ref{generalpos} we can change $g$ by a homotopy
relative to its boundary to a map whose 
image is in $\SU^\an$, hence
$[f]\in \pi_1(\SU^\an,x)$ is zero.

To check that $\pi_1(i^\an)$ is surjective, let 
$$
f: \sone \too X^\an\times S^\an
$$
be a continuous based loop. 
Applying lemma \ref{generalpos} we can change
$f$, by a homotopy relative to the endpoints of the
interval, to a based loop in 
$\SU^\an$.
\end{proof}

\begin{corollary}
With the same notation and hypothesis 
as in lemma \ref{extension},
the inclusion $i$ induces an isomorphism of algebraic
fundamental groups
$$
\pi^{\rm alg}(i,x):
\pi^{\rm alg}(\SU,x)\stackrel{\cong}{\too} \pi^{\rm alg}( X\times S,x).
$$
\end{corollary}

\begin{proof}
The algebraic fundamental group is canonically 
isomorphic to the completion
of the topological fundamental group
with respect to the topology of finite index subgroups
(cfr.  \cite[XII Cor. 5.2]{SGA1}), hence the result
follows from lemma \ref{extension}.
\end{proof}

The monomorphism $\rho_2: G/Z \inj \autlie(\fgp)$ is 
the inclusion of the 
connected component of the 
identity of $\autlie(\fgp)$. Thus
$F=\autlie(\fgp)/(G/Z)$ is a finite group.

Recall that the tautological family (\ref{family1})
parametrized by $R_1$ is denoted  
$$
\SE_{R_1}=( E_{R_1}, \varphi_{R_1})
$$
Let $\SU_{R_1}\subset X\times R_1$ be the open set where $E_{R_1}$ is
locally free. Then $\SE_{R_1}$ gives a principal $\autlie(\fgp)$-bundle
$P_{R_1}$ on $\SU_{R_1}$. Consider the functor $\Gamma(\rho_2,P_{R_1})$
of reductions defined as in (\ref{reduction2}).

\begin{proposition}
The functor $\Gamma(\rho_2,P_{R_1})$ is represented 
by a scheme $R_2\to R_1$ which 
is \'etale and finite
over $R_1$,
so there is a tautological family parametrized
by $R_2$
\begin{equation}
\label{family2}
(q_{R_2}^{}, P_{R_2}^{G/Z}, E_{R_2}^{}, \psi_{R_2}^{} )
\end{equation}
\end{proposition}

\begin{proof}
The set of isomorphism classes of $S$-families of 
$\rho_2$-reductions is bijective to
the set
\begin{equation}
\label{hom1}
\Mor_{\SU_S}(\SU_S , P_{S}(F))
\end{equation}
of sections of the pulled back principal $F$-bundle $P_{S}(F)\to
\SU_S$.  Since $F$ is a finite group, giving the principal $F$-bundle
$p:P_{R_1}(F) \to \SU_{R_1}$ is equivalent to giving a representation
of the algebraic fundamental group $\pi^{\rm alg}(\SU_{R_1},x)$ in $F$
(\cite[V \S 7]{SGA1}).  By lemma \ref{extension} this fundamental
group is isomorphic to $\pi^{\rm alg}(X\times R_1,x)$, so
there is a unique principal $F$-bundle $\overline{P_{R_1}(F)}$ on
$X\times R_1$ whose restriction to $\SU_{R_1}$ is isomorphic to
$P_{R_1}(F)$.  We claim that the set (\ref{hom1}) is bijective to
\begin{equation}
\label{hom2}
\Mor_{X\times S}(X\times S , \overline{P_{R_1}(F)}_{S}).
\end{equation}
Indeed, an element of the set (\ref{hom1}) corresponds to 
a trivialization of the principal bundle $P_{S}(F) \to \SU_{S}$.
If this is trivial, then the principal bundle  
$\overline{P_{R_1}(F)}_{S}\to X\times S$ will also be trivial, 
and trivializations of the former are in bijection with
trivializations of the later, and these correspond to
elements of (\ref{hom2}), thus proving the claim.

Finally, the morphism 
$X\times R_1 \to R_1$ is projective and
faithfully flat, $\overline{P_{R_1}(F)} \to X\times R_1$
is an \'etale and surjective, 
and $\overline{P_{R_1}(F)} \to R_1$ is projective. It
follows from \cite[lemma 4.14.1]{Ra3} that the
functor $\Gamma(\rho_2,P_{R_1})$ is representable
by a scheme $R_2 \to R_1$ which is \'etale and finite
over $R_1$. 
\end{proof}

From proposition, together with proposition \ref{410}, 
we obtain the following

\begin{corollary}
The family $\SP_{R_2}^{}=(P^{G/Z}_{R_2},E_{R_2}^{}, \psi_{R_2}^{})$
is a universal family with group $\pglv$ for the functor
$\wt{F}^{\tau}_{G/Z}$.
\end{corollary}

Recall $G'=[G,G]$ denotes the commutator subgroup. Clearly 
$G/G'\cong \CC^{*q}$,
and giving a principal $G/G'$-bundle is equivalent to 
giving $q$ line bundles.
Note that $G/Z\times G/G'=G/Z'$, where $Z'$ is the center of $G'$.
Denote the projection in the first factor by
$$
\rho'_{2}:G/Z'\to G/Z.
$$
Let $d_1,\ldots,d_q$ be $q$ fixed elements of $H^2(X,\CC)$. 
Define
$$
R'_2=J^{d_1}(X)\times \cdots \times J^{d_q}(X)\times R_2,
$$
where $J^{d_i}(X)$ is the Jacobian parameterizing line bundles
on $X$ with first Chern class equal to $d_i\in H^2(X,\CC)$. 
Using a Poincar\'e line bundle on $J^{d_i}(X)\times X$, 
we construct a tautological family parametrized by $R'_2$
\begin{equation}
\label{family2p}
(q_{R'_2}^{}, P_{R'_2}^{G/Z'}, E_{R'_2}^{}, \psi_{R'_2}^{})
\end{equation}
where the principal $G/Z'$-bundle $P_{R'_2}^{G/Z'}$ is the
product of the pullback of 
the principal $G/Z$-bundle $P_{R_2}^{G/Z}$ of 
the family (\ref{family2}), and the principal
$\CC^*$-bundles associated to line bundles
on $X\times R'_2$
pulled back from Poincar\'e line bundles on $X\times J^{d_i}$.

\begin{lemma}
\label{redr2p}
The scheme $R'_2$ over $R_2$ 
represents the functor $\Gamma(\rho'_2,P_{R_2})$.
\end{lemma}

\begin{proof}
It follows easily from the construction of $R'_2$.
\end{proof}

There is a lift of the trivial $\CC^*$ action on the Jacobian $J(X)$
to the Poincar\'e bundle, providing it with a structure of a 
universal family with group
$\CC^*$. Using this action, we obtain from 
lemma \ref{redr2p} and 
proposition \ref{410} the following
 
\begin{corollary}
\label{repr2p}
There is a natural action of $G/G'\times \pglv$ on
the family of principal $G/Z'$-sheaves
$\SP^{G/Z'}_{R'_2}
=(P^{G/Z'}_{R'_2}, E_{R'_2}^{}, \psi_{R'_2}^{})$,
providing it with a structure of 
universal family with group $G/G'\times \pglv$
for the functor $\wt{F}^\tau_{G/Z'}$.
\end{corollary}

\section{Construction of $R_3$}

Let $Z'$ be the center of the commutator subgroup $G'=[G,G]$. 
It is a finite abelian group. Consider the exact sequence of groups
\begin{equation}
\label{exse}
1 \too Z' \too G \stackrel{\rho_3}{\too} G/Z' \too 0.
\end{equation}
Recall that the family (\ref{family2p}) parametrized by $R'_2$ 
provides a principal $G/Z'$-bundle
\begin{equation}
\label{pringz}
P^{G/Z'}_{R'_2}\too\SU_{R'_2}^{}\subset X\times R'_2,
\end{equation}
where $\SU_{R'_2}$ is the open set where the torsion free sheaf 
$E_{R'_2}$ of (\ref{family2p})
is locally free.

We first recall some facts about nonabelian cohomology.
For a scheme $Y$ and a group $H$, we denote by $\underline{H}$
the trivial \'etale sheaf on $Y$ with fiber $H$.
Given a morphism $p:Y \to S$, we
define $R^ip_*(\underline{H})$ the \'etale sheaf on $S$ 
associated to the presheaf
$$
(u:U\to S) \longmapsto \check H^i_\et(Y_U,\underline{H}),
$$
where $\check H^i_\et$ denotes the Cech cohomology set with 
respect to the \'etale topology, and $Y_U=Y\times_S U$
For a finite abelian group $F$, let
$H^i(Y^\an;F)$ be the singular cohomology of
$Y^\an$ with coefficients in $F$.
We will need the following comparison

\begin{theorem}
\label{thmcomparison}
Let $F$ be a finite abelian group, and $Y$ a scheme, locally 
of finite type. Then there is a canonical isomorphism
$$
\check H^i_\et(Y,\underline{F})
\;\isom\;
H^i(Y^\an;F)
$$
\end{theorem}

\begin{proof}
Follows from \cite[XVI Th. 4.1]{SGA4} 
($\check H^i_\et(Y,\underline{F})\isom 
H^i_{\rm cl}(Y^\an;\underline{F})$) and the fact that \'etale cohomology
can be calculated using Cech cohomology.
\end{proof}

\begin{lemma}
\label{hi}
Let $p:\SU_{R'_2}\to R'_2$ be the projection to $R'_2$. 
Then, for $i\leq 2$,
$$
R^i p_* \underline{Z'} = \underline{H^i(X^\an;Z')},
$$
i.e. $R^ip_* \underline{Z'}$ is the constant sheaf with fiber
$H^i(X^\an;Z')$, the singular cohomology group of $X^\an$ with
coefficients in $Z'$.
\end{lemma}

\begin{proof}
Let $U\to R'_2{}$ be an \'etale open set of $R'_2{}$,
and let $\SU_U=\SU_{R'_2}\times_{R'_2} U$.
The isomorphism of the homotopy groups in lemma \ref{extension}
provides an isomorphism of the singular homology groups
$$
H_1(\SU_U{}^\an;\ZZ) \stackrel{\isom}{\too} H_1(X^\an\times U^\an;\ZZ)
$$
Now we will show that 
$$
H_2(\SU_U{}^\an;\ZZ) \too H_2(X^\an\times U^\an;\ZZ)
$$
is an isomorphism.
To check that it is injective, consider a class 
$\alpha$ in $H_2(\SU_U{}^\an;\ZZ)$ which maps
to zero. This class is represented by a sum
with integer coefficients $\sum n_i f_i$,
where $f_i:M^2_i\to \SU_U{}^\an$ are continuous
maps with $M^2_i$ a polyhedron of real dimension
2. Since it maps to zero, there is a 3-dimensional
singular chain $\beta$ in $X^\an \times U^\an$,
represented by a sum with integer coefficients
$\sum m_j g_j$, where the $g_j:M^3_j\to X^\an \times U^\an$ are 
continuous maps with $M^3_j$ a polyhedron of real dimension
3, and we can assume that the boundary of $M^3_j$ is
mapped to the union of the images of $f_i$. In particular,
the image of this boundary is in $\SU_U{}^\an$.

By lemma \ref{generalpos}, each map $g_j$ can be changed by a homotopy,
relative to its boundary, to a
map $\wt g_j$ whose image lies in $\SU_U{}^\an$.
Then $\sum m_j \wt g_j$ is a cycle in $\SU_U{}^\an$
whose boundary is $\sum n_i f_i$, 
hence $\alpha$ is already zero in $H_2(\SU_U{}^\an;\ZZ)$.

To check surjectivity, note that a singular cocycle
in $X^\an\times U^\an$ can be represented by a sum
$\sum n_i f_i$ where, for each $i$,
$$
f_i:M^2_i\too X^\an\times U^\an
$$
is a continuous map from  $M^2_i$, a
closed manifold with real dimension 2
with a triangulation. By lemma \ref{generalpos}
the map $f_i$ can modified by a homotopy to a map $\wt f_i$
whose image lies in $\SU_U{}^\an$.
This modification does not change the homology class,
so this proves surjectivity.

The inclusion $j:\SU_U{}^\an \inj X^\an \times U^\an$ induces
an isomorphism 
$$
j^*:H^i(X^\an \times U^\an;Z')
 \stackrel{\cong}\too 
H^i(\SU_U{}^\an;Z')
$$
for $i=1$ or $2$. 
Indeed, denoting $\SU=\SU_U{}^\an$ and 
$\SM=X^\an \times U^\an$, 
the inclusion induces a commutative diagram
$$
\xymatrix{
{0} \ar[r] & {\Ext^1(H_{i-1}(\SM;\ZZ),Z')} \ar[r] \ar[d]^{\isom}&
{H^i(\SM;Z')} \ar[r] \ar[d]^{j^*}& {\Hom(H_i(\SM;\ZZ),Z')}
\ar[r]\ar[d]^{\isom} & {0}\\
{0} \ar[r] & {\Ext^1(H_{i-1}(\SU;\ZZ),Z')} \ar[r] &
{H^i(\SU;Z')} \ar[r] & {\Hom(H_i(\SU;\ZZ),Z')}
\ar[r] & {0}
}
$$
where the exact rows are given by the universal coefficient theorem
for singular cohomology (\cite[Ch. 5 \S 5]{Sp}), and then
$j^*$ is an isomorphism by the 5-lemma.

By theorem \ref{thmcomparison}, 
\'etale cohomology coincides with singular cohomology, hence
taking sheafification we obtain
$$
R^ip_* \underline{Z'} 
\stackrel{\cong}\too  \underline{H^i(X;Z')}.
$$
\end{proof}

Given a scheme $Y$, the exact sequence (\ref{exse}) 
gives an exact sequence of pointed sets
\cite{Se1}, \cite{F-M} 
$$
\check H^1_\et(Y,\underline{G}) \too 
\check H^1_\et(Y,\underline{G/Z'}) \too 
\check H^2_\et(Y,\underline{Z'})
$$
where the distinguished element for each set corresponds
to the trivial cocycle (and exactness means that the 
inverse image of the distinguished element of the last set
is equal to the image of the first map).

This exact sequence implies that, if $p:Y\to S$ is a morphism of
schemes, there is an exact sequence of sheaves of sets on $S$
\begin{equation}
\label{secrel}
R^1p_*\underline{G} \too R^1p_*\underline{G/Z'} \too
R^2p_*\underline{Z'},
\end{equation}
which can be thought of as the relative version of the previous 
sequence.

\begin{lemma}
\label{red1}
Assume there is a reduction $(P^G,\zeta)$
to $G$ of an algebraic principal $G/Z'$-bundle $P$ on a scheme $Y$.
Then the set of algebraic isomorphism classes of reductions
is an $H^1(Y^\an;Z')$-torsor.
\end{lemma}

\begin{proof}
Recall that this means that $H^1(Y^\an;Z')$ 
acts simply transitively on this set, i.e.
it is a principal $H^1(Y^\an;Z')$-bundle over a point, 
and hence, for each reduction
$(P^G,\zeta)$, there is a natural bijection
between $H^1(Y^\an;Z')$ and 
the set of isomorphism classes of reductions,
sending the zero element
of $H^1(Y^\an;Z')$ 
to $(P^G,\zeta)$.

Since $Z'$ is discrete abelian, $H^1(Y^\an;Z')=
\check H^1_\et(Y,\underline{Z'})$
(theorem \ref{thmcomparison}).
The action of this group on the set of reductions is
defined as follows. Let $(P^G,\zeta)$ be an analytic reduction,
and $\alpha\in \check H^1_\et(Y,\underline{Z'})$. 
Let $\{g_{ij}\}$ be a $\underline{G}$-cocycle 
representing the isomorphism class of $P^G$, and let $\{z_{ij}\}$
be a cocycle representing $\alpha$.
Then  $\{g_{ij}z_{ij}\}$ defines a 
principal $G$-bundle $\hat P^{G}$
and, using $\zeta$, an isomorphism $\hat\zeta:\rho_{3*}(\hat P^{G})\isom
P$. The action is
$$
(P^G,\zeta)\cdot \alpha = (\hat P^G,\hat\zeta).
$$ 
It is easy to check that this is well defined on
the set of isomorphism classes of reductions, and 
the action is simply transitively.
\end{proof}

\begin{remark}
\textup{
In the previous proof we have used the fact that $Z'$ is 
in the center of $G$.
In general the set of reductions is bijective to a cohomology
set with twisted coefficients.}
\end{remark}

The relative version of this bijection is as follows. Assume
that we have a morphism of schemes $p:Y\to S$, and a principal 
$G/Z'$-bundle $P_S$ on $Y$

\begin{lemma}
\label{lemmaredsheaf}
Let $P_S^G$ be a principal $G$-bundle on $Y$, with 
$\rho_{3*}P_S^G\isom P_S$. Then, for all 
\'etale open sets
$U\to S$
\begin{equation}
\label{redsheaf}
\wt\Gamma(\rho_3,P_S)(U)= R^1p_* \underline{Z'}(U)
\end{equation}
where $\wt\Gamma(\rho_3,P_S)$ is the 
sheaf of reductions defined in the preliminaries.
\end{lemma}

\begin{proof}
Lemma \ref{red1} gives a bijection (depending only on $P^G_S$)
$$
\Gamma(\rho_3,P_S)(U) = \check H^1_\et(Y_U, \underline{Z'})
$$
Sheafifying sides, we obtain the result.
\end{proof}

\begin{proposition}
\label{r3p}
The functor $\wt\Gamma(\rho_3,P^{G/Z'}_{R'_2})$
is representable by a scheme $R'_3$ \'etale and finite over $R'_2$.
\end{proposition}


\begin{proof}
The strategy of the proof is as follows.
First we see that the subscheme $\hat R'_2\subset R'_2$
corresponding to principal bundles that admit a reduction 
of structure group to $G$ is a union of connected components of $R'_2$.
Then we show that the functor $\wt\Gamma(\rho_3,P^{G/Z'}_{R'_2})$
is a principal space over $\hat R'_2$,
and the structure group of this principal space is 
the finite group $H^1(X^\an;Z')$, hence affine,
and then it follows from descent theory that the functor is
representable \cite{FGA}.

The principal $G/Z'$-bundle $P^{G/Z'}_{R'_2}\to \SU_{R'_2}$ 
(cfr. \ref{pringz})
gives a section $\sigma'$ of $R^1p_* \underline{G/Z'}$ over $R'_2$,
and using (\ref{secrel}) we obtain
a section of $R^2p_*\underline{Z'}$. 
The principal $G/Z'$-bundle corresponding to a point in $R'_2$ can be lifted
to $G$ if and only if this section is zero at this point.
By lemma \ref{hi} this sheaf is constant with finite fiber, hence
the section is locally constant and it vanishes in a subscheme
$\hat R'_2 \subset R'_2$, which is a union of certain 
connected components of $R'_2$.

By exactness of the sequence (\ref{secrel}), we can cover $\hat R'_2$
with open sets $U_i$ (in the \topol topology) such that the
section $\sigma'|_{U_i}$ of $R^1p_* \underline{G/Z'}$ over $U_i$
lifts to a section $\sigma_i$ of $R^1p_* \underline{G}$. Refining
the cover $U_i$ if necessary, we can assume that 
$$
\sigma_i\in H^1(\SU_{U_i},\underline{G}).
$$
This means that there are principal $G$-bundles $P^G_i \to \SU_{U_i}$
such that $\rho_{3*} P^G_i \isom P^{G/Z'}_{U_i}$.
The action of $H^1(X^\an;Z')$ described in the proof of 
lemma \ref{red1} gives an action $\Theta$ on the functor
of reductions $\wt \Gamma (\rho_3,P^{G/Z'}_{R'_2})$.
By lemma \ref{lemmaredsheaf}, after restricting to $U_i$
we have an equality of functors
$$
\wt\Gamma(\rho^{}_3,P^{G/Z'}_{U_i})=R^1p_*\underline{Z'}|_{U_i} :
(\Sch/U_i) \too (\Sets)
$$
By lemma \ref{hi}, $R^1p_*\underline{Z'}$ is the sheaf
of sections of $R'_2\times H^1(X^\an;Z') \to R'_2$,
and then 
$\wt\Gamma(\rho_3,P^{G/Z'}_{U_i})$ is represented by 
the scheme $U_i \times H^1(X^\an;Z')$,
the action $\Theta$ becoming just multiplication on the 
right. Hence the functor $\wt\Gamma(\rho_3,P^{G/Z'}_{R'_2})$
is a principal space with group $H^1(X^\an;Z')$.
Since this group is affine,
by descent theory it follows that it
is represented by a 
principal $H^1(X^\an;Z')$-bundle over $\hat R'_2$, and 
the result follows.
\end{proof}

Let $R_3^{}\subset R'_3$ be the union of components corresponding to principal
$G$-sheaves with fixed numerical invariants $\tau$. 
The morphism $R_3 \to R'_2$ is also finite.
Then, proposition \ref{r3p} together with corollary \ref{repr2p}
and proposition \ref{410} conclude


\begin{corollary}
\label{repr3}
The scheme $R_3$ is a universal space with group $\pglv$ for
the functor $\wt{F}_G$.
\end{corollary}

Note that we have used the fact that the action of $G/G'$ on
$R_2$ is trivial.

\section{Construction of a quotient}
\label{sec4}

Let $H$ be a reductive algebraic group acting on two schemes $T$ and $S$.
We will use the following (\cite[lemma 5.1]{Ra3})

\begin{lemma}[Ramanathan]
\label{ramanathan}
If $f:T\to S$ is an affine $H$-equivariant morphism and 
$p:S\to \hat{S}$ is a good quotient for the action of $H$, 
then there is
a good quotient $q:T \to \hat{T}$ by $H$, and
the induced morphism $\hat{f}:\hat{T}\to \hat{S}$
is affine.

Furthermore, if $f$ is finite, then $\hat{f}$ is finite.
When $f$ is finite and $p:S\to \hat{S}$ is a geometric
quotient, then $q:T \to \hat{T}$ is also a geometric quotient.
\end{lemma}

\begin{theorem}
There is a projective scheme $\mathfrak{M}^\tau_G$ corepresenting the 
functor $\wt{F}_G$ of families of semistable principal $G$-sheaves
with numerical invariants $\tau$. 
There is an open subscheme $\mathfrak{M}^{\tau,s}_G$ whose closed
points are in bijection with isomorphism classes of stable
principal $G$-sheaves.
\end{theorem}

\begin{proof}
We use the notation of proposition \ref{r}.
Using geometric invariant theory,
it is proved in \cite{G-S1} that there is a good quotient
for the action of $\slv$ on the scheme $R$ of based
$\delta$-semistable Lie tensors
$$
p^{}_R: R \too \overline{R}\gitq \slv,
$$
where $\overline{R}$ is the closure of $R$ defined in
proposition \ref{r},
and $\overline{R}\gitq \slv$ is a projective scheme,
and that it is a geometric quotient on the open subscheme
$R^s$ of based $\delta$-semistable Lie tensors.
By proposition \ref{r1}, the inclusion 
of based semistable $\fgp$-sheaves 
$R_1 \inj R$ is proper, hence
the restriction of $p_R$
$$
p^{}_{R_1}: R_1 \too R_1/\slv = \mathfrak{M}_1,
$$
is also a good quotient onto a projective scheme, and
it is a geometric quotient on the open set $R_1^s$ 
corresponding to based stable $\fgp$-sheaves.
Since the center of $\slv$ acts trivially on $R_1$, this is
also a quotient by $\pglv$.

For the scheme $R_3$ of based semistable principal $G$-sheaves,
i.e. pairs $(q,\SP)$ where 
$\SP=(P,E,\psi)$ is a semistable principal $G$-sheaf and
$q:V\otimes \SO_X(-m)\surj E$ is a surjection inducing
an isomorphism $V\isom H^0(E(m))$, 
the following composition is a finite morphism
$$
f: R_3 \too R'_2 =
J^{\underline{d}}\times R_2
\too J^{\underline{d}} \times R_1,
$$
where $J^{\underline{d}}=J^{d_1}(X)\times \cdots \times J^{d_q}(X)$.
Let $\pglv$ act trivially on $J^{\underline{d}}$. Then
$$
p:J^{\underline{d}}\times R_1 \too J^{\underline{d}}\times R_1/\slv
$$
is a good quotient by $\pglv$, whose restriction to 
$J^{\underline{d}}\times R_1^s$ is a geometric quotient.
Therefore, by lemma \ref{ramanathan}, there exists a good 
quotient by $\pglv$
$$
q:R_3 \too \mathfrak{M}^\tau_G
$$
which is a geometric quotient on the subscheme $R_3^s$
of based stable principal $G$-sheaves.
Furthermore, the induced morphism 
$\overline{f}: \mathfrak{M}^\tau_G \to 
J^{\underline{d}} \times \mathfrak{M}_1$ is finite,
hence $\mathfrak{M}^\tau_G$ is projective.

By corollary \ref{repr3}, the scheme $R_3$ is a universal space with group 
$\pglv$ for the functor $\wt{F}_G$, hence, by remark
\ref{corep}, the projective scheme $\mathfrak{M}^\tau_G$
corepresents the functor $\wt{F}_G$.

The last statement follows also from Ramanathan's lemma, because
$f$ is finite.
\end{proof}

\textit{Two semistable principal sheaves are called GIT-equivalent if they
correspond to the same point in the moduli space.} 
Now we will show that this amounts to the notion
of S-equivalence given in the introduction 
(definition \ref{sequiv}).

Let $\SP=(P,E,\psi)$ be a semistable principal
sheaf. 
If it is not stable, let $E_\bullet$, or
$E_{\lambda_\bullet}$ be 
an admissible filtration, 
i.e. a balanced algebra filtration
with 
\begin{equation}
\label{admissible}
\sum_{i\in \ZZ}
\big( r P_{E_{i}} -r_{i} P_E \big)
\; = \;
\sum_{i=1}^{t} (\lambda_{i+1}-\lambda_i)
\big( r P_{E_{\lambda_i}} -r_{\lambda_i} P_E \big)
\; = \; 0 \, .
\end{equation}
Let $U'$ be the open subset of $X$ where it is 
a vector bundle filtration. By lemma \ref{bijection} this 
bundle filtration amounts
to a reduction
$P^Q$ of $P|_{U'}$ to a parabolic subgroup $Q\subset G$
together with a character $\chi$ of the Lie algebra of $Q$.
Let $Q\surj L$ be its Levi
quotient, and  $L\inj Q\subset G$ a splitting. 
In the introduction we called the principal $G$-sheaf
$$
\big(
P^Q(Q\surj L\inj G),
\oplus E^i,
\psi'
\big) 
$$
the 
\textit{admissible deformation of $\SP$ associated to
$E_{\bullet}$},
whose associated $\fgp$-sheaf is
$\oplus [,]^{i,j}:E^i \otimes E^j \to E^{i+j\,\vee\vee}$.

\begin{proposition}
\label{propsequiv}
Any admissible deformation of a semistable principal $G$-sheaf $\SP$ is semistable.
After a finite number of admissible deformations, a principal 
$G$-sheaf is obtained such that 
any further admissible deformation is
isomorphic to itself. This principal $G$-sheaf depends only on $\SP$,
and we denote it $\grad\SP$ (and $\grad\SP:=\SP$ if $\SP$ is stable).

Two principal sheaves $\SP$ and $\SP'$ are GIT-equivalent if and only if 
they are S-equivalent in the sense that
$\grad\SP \isom \grad\SP'$.
\end{proposition}

\begin{proof}
Let $z\in R_3$ and let $\overline{\slv\cdot z}$ be the closure of
its orbit. It is a union of orbits, and
by definition of good quotient, it has a unique closed orbit
$B_3(z)$, which is characterized as the unique orbit
in $\overline{\slv\cdot z}$ with minimal dimension.
Thus, two points $z$ and $z'$ in $R_3$ are GIT-equivalent 
(i.e. mapped to the same point in the moduli space)
if and only if $B_3(z)=B_3(z')$.

\medskip
\noindent\textbf{Claim.}
If $\slv\cdot z$  is not closed, then there exists a one-parameter
subgroup $\lambda$ of $\slv$
with $\mu(f(z),\lambda)=0$ such that the limit $z_0=\lim_{t\to 0}
\lambda (t) \cdot z$ is in $\overline{\slv \cdot z}\setminus \slv
\cdot z$. 

Indeed, recall that we have
a finite $\slv$ equivariant morphism
$$
R_3 \stackrel{f}{\too}  J^{\underline{d}} \times R_1
\,\subset\,  J^{\underline{d}} \times \overline{R_1} 
$$
where $\overline{R_1}$ is the closure of $R_1$ in 
the projective variety $\overline{R}$ defined in
proposition \ref{r}.
Note that $J^{\underline{d}} \times R_1$ is the
open subscheme of semistable points of the
projective variety $J^{\underline{d}} \times \overline{R_1}$.
Since  $z$ is not in $B_3(z)$, the point $f(z)$ is not in 
$B(f(z))$
(the closed orbit in the closure of $\slv\cdot f(z)
\subset J^{\underline{d}} \times R_1$), 
because the morphism $f$ sends orbits to orbits and $\dim(f(\slv\cdot z))=
\dim(\slv\cdot f(z))$, since $f$ is equivariant and finite.
By \cite[lemma 1.25]{Si}, there is a one parameter subgroup
$\lambda$ of $\slv$ such that 
$\overline{f(z)}:=\lim_{t\to 0} \lambda (t) \cdot f(z)\in B(f(z))$.
Since $f(z)$ is semistable, $\mu(f(z),\lambda)\leq 0$.
If this inequality were strict, then 
$\mu( \overline{f(z)} ,\lambda^{-1})>0$, which is 
impossible because $\overline{f(z)}$ is a semistable
point. Therefore $\mu(f(z),\lambda)= 0$.
Since $f$ is proper, $\lim_{t\to 0} \lambda (t) \cdot z $ exists,
and furthermore it belongs to $B(z)
\subset \overline{\slv\cdot z}\setminus \slv\cdot z$,
thus proving our claim.
\smallskip

For any one-parameter subgroup with
$\mu(f(z),\lambda)=0$, 
$\lim_{t\to 0} \lambda (t) \cdot f(z)$
exists and is semistable \cite[Prop. 2.14]{G-S}, and since $f$ is proper, 
$\lim_{t\to 0} \lambda (t) \cdot z$
also exists in $R_3$.

\medskip
\noindent\textbf{Claim.}
There is a bijection between one-parameter subgroups
of $\slv$ with $\mu(f(z),\lambda)=0$ on the one side,
and admissible 
($P_{E_{\lambda_\bullet}}=0$)
saturated balanced algebra filtrations
$E_{\lambda_\bullet}$ of $E$
together with a splitting
of the induced filtration  $H^0(E_{\lambda_{\bullet}}(m))$ in $V$
on the other side.

Indeed, in \cite{G-S1} we established a bijection between
one-parameter subgroups of $\slv$ with $\mu(f(z),\lambda)=0$ and 
balanced filtrations with
$$
P_{E_{\lambda_\bullet}} + \mu_{\rm tens}(E_{\lambda_\bullet},\phi)\delta =0
$$
where $(E,\phi)$ is the tensor corresponding to the
point $f(z)$. Therefore, the $\delta$-semistability of this tensor 
implies that the filtration
$E_{\lambda_\bullet}$ is saturated (since the left hand side
of the former  equality is bigger for the saturation).
The leading coefficient is
\begin{equation}
\label{equal0}
\sum_{i=1}^t (\lambda_{i+1}-\lambda_i)( \deg E_{\lambda_i} \rk E -
\rk E_{\lambda_i} \deg E ) + \mu_{\rm tens}(E_{\lambda_\bullet},\phi)\tau=0
\end{equation}
By lemma \ref{lemma0}, $\deg E=0$. Lemma \ref{slopess} implies
$\deg E_{\lambda_i}\leq 0$, and recall $\tau>0$.  Therefore lemmas
\ref{samemu} and \ref{mu0} imply 
$\mu_{\rm tens}(E_{\lambda_\bullet},\phi)=
\mu(E_{\lambda_\bullet},\varphi)\leq 0$. Since we have
equality in  (\ref{equal0}), it must
be $\mu(E_{\lambda_\bullet},\varphi)= 0$.
Hence, by lemma \ref{mu0}, the filtration $E_{\lambda_\bullet}$
is an algebra filtration, thus proving the claim.
\medskip

Now, let $\SP=(P,E,\psi)$ be a semistable principal $G$-sheaf,
choose a 
quotient $q:V\otimes\SO_X(-m)\to E$, and
let $z\in R_3$ be the point
corresponding to the based principal $G$-sheaf $(q,\SP)$. 
Let $\lambda:\CC^* \to \slv$ be the one-parameter subgroup associated
to an admissible saturated algebra filtration. 
The action
of $\lambda$ on the point $z$ define a morphism
$\CC^*\to R_3$ that extends to
$$
h:T=\CC \too R_3 \, ,
$$
with $h(t)=\lambda(t)\cdot z$ for $t\neq 0$
and $h(0)=\lim_{t\to 0} \lambda(t)\cdot z=z_0$.  
In the rest of this section we shall show that
the point $z_0$
corresponds to 
the associated admissible deformation. 
Then it will follow that the limit 
$z_0$ fails to be
in the orbit of $z$ if and only if the associated 
admissible deformation
fails to be isomorphic to $\SP$.

If $z_0$ is not in the orbit of $z$,
since $\slv \cdot z_0 \subset \overline{\slv
\cdot z}\setminus \slv \cdot z$, it is $\dim \slv \cdot z_0 < 
\dim \slv \cdot z$, so if we iterate this process
(with $z_0$ and another one-parameter
subgroup as before) we get a sequence 
of points $z_0$, $z^\prime_0$, $z^{\prime\prime}_0$,...
that must stop giving a point in $B(z)$. 
Hence, the principal $G$-sheaf $\grad \SP$, up to isomorphism,
depends only on $\SP$, because there is only one closed orbit in
$\overline{\slv\cdot z}$.

To finish the proof of the proposition
it only remains to show that
the point $z_0$
corresponds to 
the associated admissible deformation.
This will be done constructing a based family 
$(q_T,\SP_T)=(q_T,P_T,E_T,\psi_T)$ 
such that $(q_t,\SP_t)$ corresponds
to the point $h(t)\in R_3$ when $t\neq 0$
and $\SP_0$ is the associated admissible deformation. 
Since $R_3$ is separated,
it will follow that $(q_0,\SP_0)=(q_0,P_0,E_0,\psi_0)$ 
corresponds to $z_0$.

First we define a based family $(q_T,E_T,\varphi_T)$
of $\fgp$-sheaves.
For any $n\in \ZZ$,
define $E_n=E_{\lambda_{i(n)}}$, where 
(recall from Preliminaries before definition \ref{gitsheaves})
$i(n)$ is the 
maximum index with $\lambda_{i(n)}\leq n$. 
Let $-N$ be a negative
integer such that $E_n=0$ for $n\leq -N$, and write $V_n=H^0(E_n(m))$.
Borrowing the formalism from \cite[\S 4.4]{H-L},
define
$$
E_T=\bigoplus_n E_n\otimes t^n \subset 
E\otimes_\CC t^{-N}\CC[t] \subset
E\otimes_\CC \CC[t,t^{-1}]
$$
$$
\begin{array}{rcccc}
q_T: V\otimes\SO_X(-m)\otimes \CC[t] & \stackrel{\gamma}\too &
\oplus_n V_n\otimes\SO_X(-m)\otimes t^n & \too  & E_T \\
 {v^n\otimes 1} & \longmapsto &
{v^n\otimes t^n} & \longmapsto  & {q(v^n)\otimes t^n}
\end{array}
$$
\begin{eqnarray*}
\varphi_T: (\oplus E_n\otimes t^n) \otimes (\oplus E_n\otimes t^n) &
\too & (\oplus E_n\otimes t^n)^{\vee\vee} \\
w_i\otimes t^i \otimes w_j\otimes t^j & \longmapsto &
[w_i,w_j]\otimes t^{i+j}
\end{eqnarray*}
where $v^n$ is a local section of $V^n\otimes \SO_X(-m)$,
and $w_i$, $w_j$ are local sections of $E_i$ and $E_j$.
Then, as in \cite[\S 4.4]{H-L}, $(q_t,E_t,\varphi_t)$ 
corresponds to $f(h(t))$ (in particular, if $t\neq 0$, then
$(E_t,\varphi_t)$ is canonically 
isomorphic to $(E,\varphi)$), and $(E_0,\varphi_0)$ is
the admissible deformation associated to 
$E_{\lambda_\bullet}$.

Now we will define the family of principal $G$-bundles
$P_T$.
The saturated balanced algebra filtration $E_{\lambda_\bullet}$
provides, by lemma \ref{bijection}, a reduction $P^Q$ of $P|_{U'}$
to a parabolic subgroup $Q$ on the open set $U'$ where
$E_{\lambda_\bullet}$ is a bundle filtration,
together with dominant character
$\chi$ of $\fq=\Lie(Q)$.
Let $Q=LU$ be a Levi decomposition of the parabolic subgroup $Q$,
and denote $\fl=\Lie(L)$, $\fu=\Lie(U)$. Let $\fh\subset \fl$ 
be a Cartan algebra.
Let $v\in \fz_\fl$ be
the element associated to $\chi$ by lemma \ref{killingdual},  
We can associate to $v$, without loss of generality,
a one-parameter subgroup 
$$
\Psi:\CC^*\to Z_L
$$
of $Z_L$, the center of the Levi factor
$L$ corresponding to $\fl$, such that $d\Psi (1)=v$. 
Indeed, on the one hand, an integer
multiple $av$ provides such a subgroup (lemma \ref{multiple}),
and on the other hand, if we replace the indexes $\lambda_i$
by $a\lambda_i$, the associated one-parameter subgroup $\lambda(t)$ 
is replaced by $\lambda(t^a)$, and $h(t)$ is replaced by $h(t^a)$,
and $v$ by $av$, but this doesn't change the limit $z_0$.

The adjoint action of $\Psi(t)$ on
any $x\in \fu$ has zero limit as $t=e^\tau\in \CC^*$ goes to zero, since
using the root decomposition 
$x=\sum_{\alpha\in R^+(\fz_\fl)} x_\alpha$ with
respect to $\fz_\fl$, this action is
$$
\Psi(t)\cdot x = \sum \Psi(t)\cdot x_\alpha = \sum e^{\tau v} \cdot x_\alpha  
= \sum e^{\tau \alpha(v)}  x_\alpha  = \sum t^{\alpha(v)} x_\alpha
$$
and the limit is zero because $\alpha\in R^+(\ft)$,
so that $\alpha(v)>0$. Therefore, since the exponential map is
$G$-equivariant with respect to the adjoint action, for any element
$u=e^x\in U$, it is
$$
\lim_{t\to 0} \Psi(t)\cdot e^x = \lim_{t\to 0} e^{\Psi(t)\cdot x }=1
$$
Thus, since $\Psi(t)$ is in the center $Z_L$ of $L$, 
the adjoint action $\Psi(t)\cdot lu = \Psi(t)^{-1}lu \Psi(t)$
on any $lu\in LU=Q$ has limit
$$
\lim_{t\to 0} \Psi(t)\cdot lu = l\, \lim_{t\to 0}  \Psi(t)\cdot u=l
$$
Let $\{g_{\alpha\beta}:U'_{\alpha\beta}\to Q\subset G\}$ be a 
1-cocycle on $U'$ describing
$P^Q|_{U'}$.
Denote by $P_T$ the principal $G$-bundle on 
$U'\times T$ described by
$$
\{\Psi(t)^{-1}g_{\alpha\beta}\Psi(t):U'_{\alpha\beta}\times T\to Q\subset  G\}
$$
Note that $\Psi$ is defined only on values $t\in\CC^*$, but the previous
observations show that this cocycle can be extended to $t=0$,
and for this special value it describes 
the principal $G$-bundle $P^Q(Q\surj L\inj G)$,
thus admitting a reduction of structure group to $L$.
Remark also that, for $t\neq 0$, there is a canonical isomorphism between
the principal $G$-bundle $P_t$ on $U'$ and 
$P|_{U'}$, hence $P_T$ extends canonically to a principal $G$-bundle
on $U_{E_T}\subset X\times T$ which we still denote $P_T$.

It remains to construct an isomorphism of vector bundles
$\psi_T:P_T(\fg')\to E_T|_{U'\times T}$.  
Let $\SW=\SO_X^{\oplus r}$, and let $\SW_n\subset \SW$ be the 
trivial subbundle defined as the direct sum of the first $\rk E_n$ 
summands, and $\SW^n=\SW_n/\SW_{n-1}$. 
Take a
covering $\{U'_\alpha\}$ 
of $U'$ with trivializations
$\psi_\alpha:\SW|_{U'_\alpha}\to  E|_{U'_\alpha}$
preserving the filtration on $E$, i.e. such that 
$\psi$ restricts to an isomorphism between 
$\SW_n|_{U'_\alpha}$ and  $E_n|_{U'_\alpha}$.
Consider the $\fg'$-sheaf isomorphism
\begin{eqnarray*}
\gamma:  \SW|_{U'_\alpha}\otimes \CC[t] 
 & \too & \oplus \SW_n|_{U'_\alpha}  \otimes t^n 
\\ v^n\otimes 1 & \longmapsto & v^n\otimes t^n
\end{eqnarray*}
where $v^n$ is a local section of $\SW^n$.
The transition functions $h_{\alpha\beta}:U'_{\alpha\beta}\to
\autfgp\subset \glgp$ of $E|_{U'}$ can be chosen to be
block-upper triangular matrices
$$
h_{\alpha\beta}= 
\left\{
\begin{array}{cccc}
M_{\lambda_1\lambda_1} & M_{\lambda_1\lambda_2} & \cdots&  M_{\lambda_1\lambda_{t+1}} \\
0                      & M_{\lambda_2\lambda_2} & \cdots&  M_{\lambda_2\lambda_{t+1}} \\
\vdots                 &  \vdots                & \ddots&  \vdots                 \\
0                      & 0                      &\cdots&M_{\lambda_{t+1}\lambda_{t+1}}
\end{array}
\right\}
$$
where $M_{\lambda_i\lambda_j}$ is a matrix of dimension 
$\rk E^{\lambda_i}\times \rk E^{\lambda_j}$.
The commutativity of the diagram
$$
\xymatrix{
{\SW |_{U'_{\alpha\beta}}\otimes \CC[t]} \ar[r]^{\gamma}_{\isom}
\ar[dd]|-{\gamma^{-1}(t) h_{\alpha\beta}\gamma(t)} & 
{\oplus \SW_n |_{U'_{\alpha\beta}}\otimes t^n} \ar@{^{(}->}[rr]
\ar[d]^{\psi\otimes \id}_{\isom} & &
{\SW|_{U'_{\alpha\beta}}\otimes t^{-N}\CC[t]} 
\ar[dd]^{h_{\alpha\beta}\otimes \id}_{\isom}
\ar[ld]^{\psi\otimes \id}_{\isom}\\
 & 
{\oplus E_n|_{U'_{\alpha\beta}} \otimes t^n} \ar@{^{(}->}[r]  & 
{E |_{U'_{\alpha\beta}}\otimes t^{-N}\CC[t]} \\
{\SW |_{U'_{\alpha\beta}}\otimes \CC[t]} \ar[r]^{\gamma}_{\isom} & 
{\oplus \SW_n |_{U'_{\alpha\beta}}\otimes t^n} \ar@{^{(}->}[rr]
\ar[u]_{\psi\otimes \id}^{\isom} & &
{\SW|_{U'_{\alpha\beta}}\otimes t^{-N}\CC[t]}
\ar[ul]_{\psi\otimes \id}^{\isom} 
}
$$
shows 
that the transition functions of $E_T|_{U'\times T}$ are 
$\gamma^{-1}(t) h_{\alpha\beta}\gamma(t):
U'_{\alpha\beta}\times T\to \autfgp\subset \glgp$, i.e.
$$
\gamma^{-1}(t) h_{\alpha\beta}\gamma(t)
=
\left\{
\begin{array}{cccc}
M_{\lambda_1\lambda_1} & M_{\lambda_1\lambda_2} t^{\lambda_2-\lambda_1} 
& \cdots&  M_{\lambda_1\lambda_{t+1}} t^{\lambda_{t+1}-\lambda_1}  \\
0                      & M_{\lambda_2\lambda_2} & \cdots&  
M_{\lambda_2\lambda_{t+1}} t^{\lambda_{t+1}-\lambda_2} \\
\vdots                 &  \ddots                &       &  \vdots                 \\
0                      & 0                      &\cdots&M_{\lambda_{t+1}\lambda_{t+1}}
\end{array}
\right\}
$$
This is
well defined for $t=0$ because 
$M_{\lambda_i\lambda_j}=0$ when $\lambda_i-\lambda_j<0$.
Since the adjoint action of $\Psi(t)$ on $h_{\alpha\beta}$
is precisely $\Psi(t)\cdot h_{\alpha\beta}=\gamma^{-1}(t)
h_{\alpha\beta}\gamma(t)$, 
we obtain an isomorphism $\psi_T:P_T(\fg')|_{U'\times T}\to
E_T|_{U'\times T}$, hence a family $\SP_T=(P_T,E_T,\psi_T)$.
Note that, for $t\neq 0$, using the canonical isomorphisms
$E_t\isom E$ and $P_t\isom P|_{U'}$, the isomorphism $\psi_t$
becomes $\psi$, hence $\psi_T$ extends to an isomorphism
$P_T(\fgp)\to E_T|_{U_{E_T}}$, which we still denote $\psi_T$.
Finally, it is easy to check that $(q_t,\SP_t)$ corresponds to $h(t)$
and $\SP_0\isom\grad \SP$.

\end{proof}

\section{Slope (semi)stability as Ramanathan (semi)stability}
\label{sec5}

In \cite{Ra2}, Ramanathan 
defines a {\em rational principal bundle} on $X$ 
as a  principal bundle $P$ over a big open set $U\subset X$, 
and gives a  notion of (semi)stability, which is a direct 
generalization of his notion of (semi)stability in \cite{Ra3}
for $\dim X=1$.

\begin{definition}[Ramanathan]
A rational principal $G$ bundle $P\to U\subset X$ is 
(semi)stable if for any reduction $P^Q$ to a parabolic subgroup $Q$
over a big open set $U'\subset U$, and for any dominant
character $\chi$ of $Q$, it is
$$
\deg P^Q(\chi) \;(\leq)\; 0.
$$
\end{definition}

Let $\SP=(P,E,\psi)$ be a principal $G$-sheaf and let $U$
be the open set where $E$ is locally free. 
We will show in this section that $\SP$ is 
slope-(semi)stable
if and only if the rational bundle $P$ is (semi)stable
in the sense of Ramanathan. In particular, 
we will obtain that, if $X$ is a curve,
our notion of (semi)stability for principal bundles
coincides with that of Ramanathan.
As mentioned in the introduction, this section
plays also the role of an appendix where we
prove some facts that have been already used.

Recall (from \cite{J}, for instance) the well known
notions of filtration and graduation of a Lie 
algebra $\fg$. An algebra filtration
$\fg_\bullet$ is a sequence 
$$
\ldots \subseteq \fg_{i-1} \subseteq \fg_{i} 
\subseteq \fg_{i+1} \ldots
$$
starting by $0$ and ending by $\fg$, such that
$$
[\fg_i,\fg_j]\subseteq \fg_{i+j}
\quad \text{for all $i$, $j\in \ZZ$}
$$
or, deleting (from $0$ onward) all nonstrict
inclusions, it is $\fg_{\lambda_\bullet}$
$$
0 \subsetneq \fg_{\lambda_1} \subsetneq \fg_{\lambda_2} 
\subsetneq \;\cdots\; \subsetneq
\fg_{\lambda_{t+1}}=\fg,
\quad (\lambda_1<\cdots<\lambda_{t+1})
$$
with
$$
[\fg_{\lambda_i},\fg_{\lambda_j}]\subseteq \fg_{\lambda_{k-1}}
\quad \text{if $\lambda_i+\lambda_j < \lambda_k$}.
$$
A graded structure $\fg^\bullet$ is a decomposition
$$
\fg= \bigoplus_{i\in \ZZ} \fg^{i} \quad
\text{with}\quad 
[\fg^i,\fg^j]\subset \fg^{i+j}
\quad \text{for all $i$, $j\in \ZZ$}
$$
or, deleting all zero summands,
$$
\fg= \bigoplus_{i=1}^{t+1} \fg^{\lambda_i}
\quad (\lambda_1<\cdots<\lambda_{t+1}).
$$
with
$$
[\fg^{\lambda_i},\fg^{\lambda_j}]\subseteq 
\left\{
\begin{array}{cl}
\fg^{\lambda_{k}}& \text{if there is $k$ with 
$\lambda_k=\lambda_i+\lambda_j$}\\
0 & \text{otherwise}
\end{array}
\right.
$$
To a graded algebra $\fg^\bullet$
it is associated a filtered algebra $\fg_\bullet$ with 
$$
\fg_i=\bigoplus_{j\leq i} \fg^j
$$
and reciprocally, to a filtered algebra $\fg_\bullet$
it is associated a graded algebra
$$
(\gr{\fg})^i= \fg_i/\fg_{i-1}
$$
with Lie algebra structure
$$
[\overline{v},\overline{w}]= [v,w] \mod \fg_{i+j-1}
$$
for $v\in \fg_i\setminus \fg_{i-1}$ and
$w\in \fg_j\setminus \fg_{j-1}$.

A graded algebra $\fg^\bullet$ is called balanced if
$\sum i \dim \fg^{i}=0$. In terms of $\fg^{\lambda_\bullet}$,
this is $\sum \lambda_i \dim \fg^{\lambda_i}=0$.
A filtered algebra is called balanced if
the associated graded algebra is so.
We start this appendix proving the following

\begin{lemma}
\label{isomfilt}
Let $\fgp_\bullet$ be a balanced algebra filtration of a 
semisimple Lie algebra $\fgp$.
There is a Lie
algebra isomorphism between $\fgp$ and the associated Lie algebra $\gr
(\fgp_{\bullet})$.
\end{lemma}

\begin{proof}
Let $W$ be the vector space underlying the Lie algebra $\fgp$.  Choose
a basis $e_l$ of $W$ adapted to the filtration
$\fgp_{\lambda_\bullet}$.  Associate the one-parameter subgroup
$\lambda(t)$ of $\glw$ expressed as ${\rm diag}(t^{\lambda_\bullet})$
in this basis. Since the filtration is balanced, this is in fact
a one-parameter subgroup of $\slw$.
The Lie algebra structure of $W$ is a point 
$v=\sum a_{lm}^n e^l\otimes e^m\otimes e_n$ 
in the linear space $W^\vee\otimes W^\vee \otimes W$.
The action of the one-parameter subgroup is  
$$
a_{lm}^n \longmapsto 
t^{\lambda_{i(l)} + \lambda_{i(m)} - \lambda_{i(n)}} a_{lm}^n,
$$
where $i(l)$ is the minimum integer for which $e_l\in \fgp_{\lambda_{i(l)}}$.
The point $\overline
v\in \PP(W^\vee\otimes W^\vee \otimes W)$ is GIT-semistable with
respect to the induced action of $\slw$ on this projective 
space and on its polarization line bundle  
$\SO_\PP(1)$ (by lemma \ref{ra551}), hence the Hilbert-Mumford criterion
implies
$$
\mu := 
\min \big\{
\lambda_{i(l)} + \lambda_{i(m)} - \lambda_{i(n)} \;:\; 
a^n_{lm}\neq 0 \big\} \leq 0
$$
Furthermore, $\mu=0$ because $\lambda_\bullet$ is an algebra
filtration. Indeed, if $\mu<0$ then for some
triple $(\lambda_i,\lambda_j,\lambda_k)$ with
$\lambda_i+\lambda_j<\lambda_k$ it would be
$[\fgp_{\lambda_i},\fgp_{\lambda_j}] \nsubseteq
\fgp_{\lambda_{k-1}}$,
contradicting the fact that $\fgp_{\lambda_\bullet}$ is 
algebra filtration.

Since $\mu=0$, the following limit exists and is nonzero
$$
v_0:= \lim_{t\to 0} \lambda(t)\cdot v \;\in 
\; W^\vee\otimes W^\vee \otimes W
$$
Since the subset of points of $W^\vee\otimes W^\vee \otimes W-\{0\}$ giving 
$W$ a Lie algebra structure is closed, the point $v_0$ itself
provides $W$ with a Lie algebra 
structure. By construction, the coordinates $b_{lm}^n$ of
$(W,v_0)$ are
$$
b_{lm}^n =
\left\{
\begin{array}{ll}
a_{lm}^n\; , & \lambda_{i(l)} + \lambda_{i(m)} - \lambda_{i(n)}=0\\
0\; , & \lambda_{i(l)} + \lambda_{i(m)} - \lambda_{i(n)}\neq 0
\end{array}
\right.
$$
In other words, 
$(W,v_0)\isom \gr (\fgp_{\lambda_\bullet})$.
Let $k(t):W\otimes W \to \CC$ 
be the Killing form of $\lambda(t)\cdot v$. Since
$\lambda(t)\in \slw$, 
$$
\det(k(t))=
\det(\lambda(t)^{-1} k(1) \lambda(t))=
\det(k(1))\neq 0 \quad \text{for all}\; t\in \CC^*,
$$
thus also for $t=0$. Since this determinant is nonzero,
$(W,v_0)$ is semisimple. By the rigidity of semisimple Lie
algebras, $(W,v_0)\isom (W,v)=\fgp$.
\end{proof}

Let $\fa$ be a 
\textit{toral algebra}
$\fa\subset \fg$, i.e. an algebra consisting of
semisimple elements, thus abelian \cite[\S 8.1]{Hum},
which is not necessarily maximal. 
Following \cite[\S 3]{B-T}, we can define 
the set $R(\fa)\subset \fa^\vee$ of $\fa$-roots 
in the following way.
For $\alpha\in \fa^\vee$, write
\begin{equation}
\label{defroot}
\fg^\alpha= \{ x\in \fg: [s,x]=\alpha(s)x, 
\, \text{for all} \; s\in \fa\}
\end{equation}
Then $R(\fa)=\{\alpha\in \fa^\vee\setminus 0: \fg^\alpha\neq 0\}$
For $\fh$ is a maximal toral algebra (i.e. Cartan algebra) 
containing $\fa$,
$\fa$-roots can be thought of as
classes of $\fh$-roots by saying that two 
$\fh$-roots are equivalent if their restrictions
to $\fa$ are the same.
Let 
$R(\fh)=R^+(\fh)\cup R^-(\fh)$ 
be a decomposition 
into positive and negative $\fh$-roots.
If $\beta\sim \beta'\nsim 0$,
then $\beta$ is positive if and only if $\beta'$ is
positive, hence there is an induced decomposition
$R(\fa)=R^+(\fa)\cup R^-(\fa) $.
In particular, this gives a partial ordering 
among $\fa$-roots: $\alpha<\alpha'$ when
$\alpha'-\alpha$ is a sum of positive $\fa$-roots.

\begin{lemma}
\label{killingdual}
Let $\fq$ be a parabolic subalgebra of a semisimple 
Lie algebra $\fgp$ and $\chi:\fq\to \CC$ a character
of $\fq$.
Let $\fq=\fl\oplus \fu$ be a Levi
decomposition, and $\fz_\fl$
the center of the Levi subalgebra $\fl$.
Then there is an element $v\in \fz_\fl$ 
such that 
$$
\chi(\cdot)=(v,\cdot):\fq\too \CC
$$
where $(\cdot,\cdot)$ is the Killing form
of $\fgp$.
\end{lemma}

\begin{proof}
Let $\fl'=[\fl,\fl]$ be the commutator subalgebra.
The decomposition 
$\fl=\fl'\oplus \fz_\fl$
is orthogonal with respect to the Killing form 
$\kappa=(\cdot,\cdot)$ on 
$\fgp$. Indeed, since $\kappa$ is $\fgp$-invariant,
if $l_1$, $l_2\in \fl$ and $z\in \fz_\fl$, then
$$
([l_1,l_2],z)=(l_1,[l_2,z])=0.
$$
Let $\fh$ be a Cartan subalgebra
of $\fgp$ 
containing $\fz_\fl$ and
contained in $\fl$.
The given decomposition of $\fl$ induces a decomposition
$\fh=(\fl'\cap \fh) \oplus \fz_\fl$ which is also $\kappa$-orthogonal.
Let $v\in \fh$ be the element in $\fh$, $\kappa$-dual
to $\chi|_{\fh}$.
The restriction $\chi|_{\fl'\cap\fh}$ is zero because 
$\fl'$ is semisimple, hence $v\in (\fl'\cap \fh)^\perp=\fz_\fl$.
\end{proof}

For a parabolic subalgebra $\fq$ and
split $\fq=\fl\oplus \fu$, 
let $R(\fz_\fl)=R^+(\fz_\fl)\cup R^-(\fz_\fl)$ be the
decomposition such that $\fgp^\alpha\subset \fq$
when $\alpha\in R^+(\fz_\fl)$.
Recall that  a character $\chi$ of $\fq$ 
is then called dominant if $2(\chi,\alpha)/(\alpha,\alpha)$
is a nonnegative integer for all 
positive $\fa$-roots $\alpha$. We call it integer if
$(\chi,\alpha)$ is integer for all $\fa$-roots
$\alpha$.

\begin{lemma}
\label{bijection}
Let $G'$ be a semisimple group.
Let $P$ be a principal $G'$-bundle over a scheme $Y$ (not necessarily
proper). 
There is a canonical bijection between the following sets
\begin{enumerate}
\item Isomorphism classes of reductions to a
parabolic subgroup $Q$ on a big open set $U\subset Y$, 
together with an integer dominant character $\chi$
of $\fq=\Lie(Q)$.

\item Isomorphism classes of saturated balanced algebra filtrations
\begin{equation}
\label{filtE2}
0 \subsetneq E_{\lambda_1} \subsetneq E_{\lambda_2} 
\subsetneq \;\cdots\; \subsetneq E_{\lambda_t} \subsetneq
E_{\lambda_{t+1}}=E
\end{equation}
of the bundle of algebras $E=P(\fgp)$ associated to $P$ by the adjoint
representation of $G'$.
\end{enumerate}
Let $\fq=\fl\oplus \fu$ be a Levi decomposition,
and $v\in \fz_\fl$ the element 
associated by lemma \ref{killingdual} to the
character $\chi$ in (1).
The set of integers 
$\{\lambda_i\}_{i=1,\ldots,t+1}$ in (2) 
is then just the set 
$\{ \alpha(v)\}_{\alpha \in R(\fa)\cup \{0\}}$
\end{lemma}

\begin{proof}
We start with a filtration (\ref{filtE2}).
Take a point $x$ of $Y$ where the filtration is a bundle filtration.
Fix an isomorphism between the fiber of $E$ at this point and $\fgp$.
We obtain a balanced algebra filtration
$\fgp_{\lambda_\bullet}$ of $\fgp$. 
By lemma \ref{isomfilt}, the 
associated graded Lie algebra
$\gr(\fgp^{\lambda_\bullet})$ is isomorphic to $\fgp$, 
and using this isomorphism we obtain a decomposition
giving $\fgp$ the structure of a graded Lie algebra
\begin{equation}
\label{gradedlambda}
\fgp= \bigoplus_{i=1}^{t+1} \fgp^{\lambda_i},
\end{equation}
such that
\begin{equation}
\label{sub}
\fgp_{\lambda_i}= \bigoplus_{j=1}^{i} \fgp^{\lambda_j},
\end{equation}
Define a linear endomorphism of $\fgp$
\begin{eqnarray*}
f:  \bigoplus_{i=1}^{t+1} \fgp^{\lambda_i} & \too &
\bigoplus_{i=1}^{t+1} \fgp^{\lambda_i} \\
v \in \fgp^{\lambda_i}& \longmapsto & -\lambda_i v
\end{eqnarray*}
If $v_i\in \fgp^{\lambda_i}$ and $v_j\in \fgp^{\lambda_j}$,
then $[v_i,v_j]\in \fgp^{\lambda_i+ \lambda_j}$ so
$$
f([v_i,v_j])=[f(v_i),v_j]+[v_i,f(v_j)]\, ,
$$
i.e. $f$ is a derivation. Thus, since 
$\fgp$ is semisimple, 
a semisimple element $v\in \fgp$ exists 
such that $f(\cdot)=[v,\cdot]$.
Let $\fz_v$ be the center of the centralizer $\fc_v$ of $v$.
It is a toral algebra.
Consider the $\fz_v$-root decomposition
(see (\ref{defroot}) or \cite[\S 3]{B-T})
\begin{equation}
\label{gradedalpha}
\fgp = \bigoplus_{\alpha\in R(\fz_v)\cup \{0\}} \fgp^{\alpha}
\end{equation}
Note that 
$\fgp^{\alpha=0}$ is just the centralizer $\fc_v$ of $v$.
This decomposition is a refinement of (\ref{gradedlambda}). 
Since
\begin{equation}
\label{refinement}
\fgp^{\lambda_i} = \bigoplus_{ \alpha(v)=-\lambda_i} 
\fgp^{\alpha}
\end{equation}

\noindent\textbf{Claim.}
The direct summand $\fgp^{\alpha=0}$
in decomposition (\ref{gradedalpha}) 
is equal to the direct summand $\fgp^{\lambda_i=0}$ 
in decomposition (\ref{gradedlambda}).

To prove this claim, let $\fz_v$-root $\alpha$ be such that
 $\alpha(v)=0$. For
$x\in \fgp^{\alpha}$ it is  $[v,x]=\alpha(v)x=0$, i.e. $x$ is in
the centralizer $\fc_v$ of $v$. By definition, $\fz_v$ is the center
of $\fc_v$, thus $[w,x]=0$ for all $w\in \fz_v$, proving the claim.

As a consequence, for all $\fz_v$-roots $\alpha$, it is
$\alpha(v)\neq 0$, and thus
$\alpha(v)>0$ gives a set of positive $\fz_v$-roots $R^+(\fz_v)$.
Using (\ref{sub}), (\ref{refinement}) and the claim,
we obtain for $\fgp_0$ in (\ref{sub})
$$
\fgp_0 = \bigoplus_{\beta\in R^{+}(\fz_v)\cup \{0\}} \fgp^{\beta},
$$
hence $\fgp_0\subset \fgp$ is a parabolic subalgebra
(\cite[\S 4]{B-T}).
Let $U$ be the big open set where $E_{\lambda_\bullet}$ 
is a bundle filtration.
The inclusion $E_0|_U\subset E|_U$ gives a reduction of 
structure group $P^Q$ of the principal $G'$-bundle $P|_U$ to the 
parabolic subgroup $Q\subset G'$ corresponding to 
$\fgp_0\subset \fgp$, because the stabilizer (under the adjoint
action of a connected group) of a parabolic subalgebra is 
the corresponding parabolic subgroup.

Finally, the character $\chi(\cdot)=(v,\cdot)$ 
of the parabolic $\fgp_0$ is dominant, 
because $(\chi,\alpha)=\alpha(v)$ is a 
positive integer for all positive $\fz_v$-roots.

Reciprocally, assume we are given a reduction $P^Q$ 
of $P$ to a parabolic subgroup
$Q$ on a big open set $U\subset Y$ 
and a dominant character $\chi$ of $\fq=\Lie(Q)$. 
Choose a decomposition $\fq=\fl\oplus\fu$ into a Levi
and a unipotent subalgebras, and let $\fz_\fl$ be the
center of $\fl$. 
Let $v\in \fz_\fl$ be the element associated to
$\chi$ by lemma \ref{killingdual}.
Consider the $\fz_\fl$-root decomposition of 
$\fgp$
(see (\ref{defroot}) or \cite[\S 3]{B-T})
$$
\fgp = \bigoplus_{\alpha\in R(\fz_\fl)\cup \{0\}} \fgp^{\alpha}\, 
$$
By hypothesis $\alpha(v)=(\chi,\alpha)$ is an integer for
all $\fz_\fl$-roots $\alpha$. Define a filtration
$\fgp_{\lambda_\bullet}$ of $\fgp$ by
\begin{equation}
\label{gsubl}
\fgp_{\lambda_i}= \bigoplus_{-\alpha(v)\leq \lambda_i}
\fgp^{\alpha}\, .
\end{equation}
This is a balanced algebra filtration of, because 
$\dim \fgp^\alpha = \dim \fgp^{-\alpha}$
and $[\fgp^\alpha, \fgp^\beta]\subset \fgp^{\alpha+\beta}$.
Clearly $\fq\subseteq\fgp_0$, and in fact $\fq=\fgp_0$ 
because the character $\chi$ of $\fq$ is dominant.
It is also clear that $\fl\subseteq\fc_v$, the centralizer
of $v$, and since $\chi$ is dominant, it is $\fl=\fc_v$,
hence the center $\fz_\fl$ of $\fl$ is the center $\fz_v$ 
of $\fc_v$.

For adjoint action of $Q$ on $\fgp$
it is
$$
Q\cdot \fgp^\alpha \subset \bigoplus_{\beta\geq \alpha} \fgp^\beta
$$
Thus the filtration (\ref{gsubl})
is preserved by this action:
$Q\cdot \fgp_{\lambda_i} \subset \fgp_{\lambda_i}$.
Since $P$ has a reduction to $Q$ on $U\subset Y$, 
this produces a vector bundle filtration of $E|_U$, 
and it extends uniquely to a saturated filtration on $Y$
as in (\ref{filtE2}).

It is easy to check that the two constructions are inverse to each
other, and by construction $\{\lambda_i\}_{i=1,\ldots,t+1} 
= \{ \alpha(v)\}_{\alpha \in R(\fz_\fl)\cup \{0\}}$.
\end{proof}

\begin{lemma}
\label{multiple}
With the same hypothesis (and notation) 
as in lemma \ref{bijection},
there are positive integers $a$ and $b$ such that 
$av$ corresponds to a one-parameter subgroup of $Z_L$
(i.e. its differential is $av$)
and 
$b\chi$ corresponds to a character of the group $Q$.
\end{lemma}

\begin{proof}
Let $\fh$ be a Cartan algebra of $\fg'$ with
$\fz_\fz\subset \fh\subset \fl$ and let $H$ be the
maximal torus of the connected group $G$ corresponding to $\fh$.  Let
$R(\fh)$ be the set of roots with respect to $\fh$.  The element $v\in
\fz_\fl \subset \fh$ is in the coweight lattice $\ZZ(W^\vee)$, 
because any $\fh$-root 
gives an integer when
evaluated on $v$. 
Indeed, the $\fz_\fl$-roots
$\alpha:\fz_\fl\to \CC$ with respect to $\fz_\fl$ are obtained by
restricting the $\fh$-roots $\beta:\fh\to \CC$ to $\fz_\fl$, but
by hypothesis, $\alpha(v)\in \ZZ$ for all
$\alpha\in R(\fz_\fl)$.
Let $X^\vee(H)$ be the lattice of  one-parameter subgroups of $H$.
Sending an element of $X^\vee(H)$ to its
differential gives an embedding 
$X^\vee(H)\inj \ZZ(W^\vee)$
with finite quotient, hence there is an integer $a$
such that $a v$ corresponds to a 
one-parameter subgroup of $H$
which can be written as
\begin{eqnarray*}
\Psi: \CC^* & \too & Z_L\subset H \\
t=e^\tau &\longmapsto & e^{\tau a v }
\end{eqnarray*}
where $Z_L$ is the center of the Lie subgroup $L$ corresponding
to $\fl\subset \fgp$.

On the other hand, the character $\chi$ of the parabolic $\fq$ 
is dominant, and 
in particular belongs to the weight lattice $\ZZ(W)$.
Let $X(H)$ be the lattice of characters of $H$.
Sending an element of $X(H)$ to its differential defines
a lattice embedding $X(H)\inj \ZZ(W)$ with finite
quotient, hence there is an integer $b$
such that $b\chi$ corresponds to a character $\Xi\in X(H)$,
i.e. the differential of $\Xi$ is $b\chi$.

Let $\fl'=[\fl,\fl]$ be the commutator subalgebra,
and $L'=[L,L]$ the commutator subgroup.
Recall that a character of $\fq$ factors as 
$\fq\surj\fl\surj \fl/\fl'\surj \CC$,
hence $\chi$ vanishes on $\fl'$.
Thus the character $\Xi$ of $H$
vanishes on $H\cap L'$, so $\Xi$ gives a group
homomorphism $L/L'\isom H/(H\cap L')\to \CC^*$.
Composing with the quotient $Q\surj L\surj L/L'$
we obtain a character of $Q$ whose differential
is $\chi$.
\end{proof}

\begin{lemma}
\label{samedegree}
Let $P$ be a principal $G'$-bundle over a big open set $U\subset X$
with a reduction $P^Q$ to a parabolic subgroup $Q\subset G'$
on a big open set $U'\subset U$. Let
$\Xi$ be a dominant character of $Q$ and $\chi$ the associated
character of $\fq$. Assume that $(\chi,\alpha)$ is an integer
for all roots of $\fgp$. Let 
\begin{equation}
\label{filtE3}
0 \subsetneq E_{\lambda_1} \subsetneq E_{\lambda_2} 
\subsetneq \;\cdots\; \subsetneq E_{\lambda_t} \subsetneq
E_{\lambda_{t+1}}=E =P(\fgp)
\end{equation}
be the balanced algebra filtration associated to it by
lemma \ref{bijection}. Then
\begin{equation}
\label{star}
\sum_{i=1}^{t+1} (\lambda_{i+1}-\lambda_i) \deg E_{\lambda_i} 
=\deg P^Q(\Xi)
\end{equation}
where $P^Q(\Xi)$ is the line bundle associated to $P^Q$ by
the character $\Xi$.
\end{lemma}

\begin{proof}
Let $L\subset Q$ be a Levi factor of $Q$, and $Z_L$ the center of $L$.
For $\fz_\fl=\Lie(Z_L)$
consider the $\fz_\fl$-root decomposition of $\fgp$
(cfr. (\ref{defroot}))
$$
\fgp = \bigoplus_{\alpha\in R(\fz_\fl)\cup \{0\}} \fgp^{\alpha}.
$$
Let $v\in \fz_\fl$ be the element associated to $\chi$
by lemma \ref{killingdual}.
Define an order $<_v$ in the set $R(\fz_\fl)\cup\{0\}$ 
by declaring
$\alpha<_v \alpha'$ if $(\alpha-\alpha')(v)<0$.
In general, $<_v$ is not a total order, because
it can happen that $(\alpha'-\alpha)(v)=0$ even if $\alpha$ and
$\alpha'$ are different.  
Choose a refinement of this to get a total order $\prec$. 
Number all the roots (including $\alpha=0$)
by $\alpha_1\succ \alpha_2 \succ \ldots \succ \alpha_{l+1}$ in descending
order, and define
a filtration $\fgp_\bullet$
\begin{equation}
\label{rootfilt}
0 \subsetneq \fgp_{\alpha_1} \subsetneq \fgp_{\alpha_2} 
\subsetneq \;\cdots\; \subsetneq \fgp_{\alpha_l} \subsetneq
\fgp_{\alpha_{l+1}}=\fgp \, ,\quad \text{with}\;\;
\fgp_{\alpha_i} = \bigoplus_{j=1}^{i} \fgp^{\alpha_j}\, .
\end{equation}
For the adjoint action of $Q$ on $\fgp$
it is
$$
Q\cdot \fgp^\alpha \subseteq \bigoplus_{\beta\geq \alpha} \fgp^\beta
\subseteq \bigoplus_{\beta\succeq \alpha} \fgp^\beta
$$
This has two consequences: on the one hand,
there is an induced action of $Q$ on
$$
(\gr\fgp)^{\alpha_i} := {\fgp_{\alpha_{i}}}/{\fgp_{\alpha_{i-1}}}
$$
and on the other hand, $P^Q$ produces a vector bundle
filtration of $E|_{U'}$, and this extends to a saturated
filtration on $U$ 
\begin{equation}
\label{thisfilt}
0 \subsetneq E_{\alpha_1} \subsetneq E_{\alpha_2} 
\subsetneq \;\cdots\; \subsetneq E_{\alpha_l} \subsetneq
E_{\alpha_{l+1}}=E
\end{equation}
Note that, although as vector spaces both 
$\fgp^\alpha$ and $(\gr\fgp)^{\alpha_i}$ are isomorphic, 
they are not isomorphic as $Q$-modules: indeed, while 
$Q \cdot (\gr\fgp)^{\alpha_i} \subset (\gr\fgp)^{\alpha_i}$,
in general we only have $Q\cdot \fgp^\alpha \subseteq 
\bigoplus_{\beta\geq \alpha} \fgp^\beta$.

The filtration (\ref{thisfilt}) is a refinement of
(\ref{filtE3}), with
\begin{equation}
\label{refine}
E_{\lambda_i} = E_{\alpha}, \quad \alpha=\max_{\prec}\big\{\beta\in
R(\fz_\fl)\cup\{0\}: -(\chi,\alpha)=-\alpha(v) \leq \lambda_i\big\}
\end{equation}
Furthermore, $E^{\alpha_i}=E_{\alpha_i}/E_{\alpha_{i-1}}$
is isomorphic to the vector bundle associated to $P^Q$ using the
action of $Q$ on $(\gr\fgp)^\alpha$.
Since this filtration is a refinement of (\ref{filtE3}),
it is
\begin{equation}
\label{sumalpha}
\deg(E^{\lambda_i}) = 
\sum_{\alpha(v)=-\lambda_i}\deg ( E^\alpha),
\end{equation}
where $E^{\lambda_i}=E_{\lambda_{i}}/E_{\lambda_{i-1}}$.

For each $\fz_\fl$-root $\alpha$ the adjoint action
of $Q$ on $(\gr \fgp)^\alpha$ gives a character
$$
\phi_\alpha:Q \stackrel{{\rm ad}}{\too} 
\gl \big((\gr \fgp)^\alpha\big) \stackrel{\det}{\too} \CC^*
$$
Every character of a parabolic subgroup factors through
its Levi quotient $L$, and two characters are equal if
they coincide when restricted to its center $Z_L$.
We have
a commutative diagram
$$
\xymatrix{
{Q} \ar[r]^{{\rm ad}} \ar@{>>}[d] & {{\gl\big((\gr \fgp)^\alpha\big)}} \ar[r]^-{\det} & 
{\CC^*}\ar@{=}[d] \\
{L} \ar[r]^{{\rm ad}} & {{\gl\big((\gr \fgp)^\alpha\big)}} 
\ar[r]^-{\det} & {\CC^*} \ar@{=}[d]\\
{Z_L} \ar[r]^{{\rm ad}} \ar@{^{(}->}[u]& {{\gl\big((\gr \fgp)^\alpha\big)}} 
\ar[r]^-{\det} & {\CC^*} \\
}
$$
It follows that 
$$
\phi_\alpha = \overline{(\dim \fgp^\alpha)\alpha},
$$
where we denote by $\overline{(\dim \fgp^\alpha)\alpha}$
the character of $Q$ such that, after restricting to 
a character $Z_L\to \CC^*$, 
the induced Lie algebra homomorphism
$\fz_\fl\to \CC$ is
$(\dim \fgp^\alpha)\alpha$.
Hence,
\begin{equation}
\label{detalpha}
\det E^\alpha \isom P^Q\big(\overline{(\dim \fgp^\alpha)\alpha}\big).
\end{equation}

Using equation (\ref{sumalpha}), the left hand side of (\ref{star}) 
is equal to the degree of the line bundle
$$
\bigotimes_{i=1}^{t+1} 
(\det E^{\lambda_i})^{ -\lambda_i}=
\bigotimes_{\alpha\in R(\fz_\fl)\cup\{0\}} 
(\det E^\alpha)^{\alpha(v)}
$$
Using (\ref{detalpha}), this 
line bundle is equal to 
\begin{equation}
\label{prev}
P^Q\big(
\overline{\sum_{\alpha\in R(\fz_\fl)\cup\{0\}}
\alpha(v)(\dim \fgp^\alpha)\alpha}
\big)
\end{equation}

\noindent\textbf{Claim.}
$$
\sum_{\alpha\in R(\fz_\fl)\cup\{0\}}
\alpha(v) (\dim \fgp^\alpha)  \alpha
= \chi
$$

Let $w\in \fz_\fl$. Then
$$
\chi(w)\;= (v,w)\;= \;\tr [([v,\cdot])([w,\cdot])] \;=\; 
\sum_{\alpha\in R(\fz_\fl)\cup \{0\}} (\dim
\fgp^\alpha)\alpha(v)\alpha(w)\, ,
$$
and the claim follows because this holds for all $w\in \fz_L$.

Since $\overline{\chi}=\Xi$, it follows that the line bundle
(\ref{prev}) is isomorphic
to $P^Q(\Xi)$, and the lemma is proved.
\end{proof}

\begin{corollary}
A principal $G$-sheaf $\SP=(P,E,\psi)$ is slope-(semi)stable
if and only if the associated rational principal $G$-bundle
$P\to U\subset X$ is (semi)stable in the sense of Ramanathan.
\end{corollary}

\begin{proof}
Without loss of generality, we can assume that $G$ is semisimple.
Assume that $\SP$ is slope-(semi)stable. 
Consider a reduction to a parabolic subgroup $Q$ of 
$P|_{U'}\to U'\subset U$, where $U'$ is a big open set, 
and a dominant character $\Xi$ of $Q$.
This gives a dominant character $\chi$ of $\fq=\Lie(Q)$.
Let $\fq=\fq\oplus \fz$ be a Levi decomposition
and $\fz_\fl$ the center of $\fl$.
A positive integer multiple $\wt\chi=c\chi$ 
has the property that $(\wt\chi,\alpha)$ is
integer for all $\fz_\fl$-roots $\alpha$.
Consider the balanced algebra filtration $\wt E^{U'}_{\lambda_\bullet}$
associated to $\wt\chi$ by lemma \ref{bijection}. 

This filtration of $E|_{U'}$ can be extended uniquely to a
saturated filtration $\wt E_{\lambda_\bullet}$ of $E$ on $X$, namely, the
intersection $\wt E_{\lambda_i}$, inside $E^{\vee\vee}$,
of $E$ and the reflexive sheaf
$F_{\lambda_i}$ extending $\wt E^{U'}_{\lambda_\bullet}$ to $X$
(cfr. \cite[II Ex. 5.15]{Ha}).  
By lemma \ref{samedegree}, and using the
slope-(semi)stability of $\SP$ we have
$$
\deg P^Q(\Xi)=
\sum_{i=1}^{t+1} \frac{\lambda_{i+1}-\lambda_i}{c}
 \deg E_{\lambda_i} \;(\leq)\;0.
$$
This means that $P\to U\subset X$ is Ramanathan (semi)stable.

Conversely, assume that $P\to U\subset X$ is Ramanathan (semi)stable. 
Consider a balanced algebra filtration of $E$.
We may assume that this filtration
is saturated.
Let $U'\subset U \subset X$ be the big open set where this 
is a bundle filtration.
Lemma \ref{bijection} produces a reduction $P^Q$ on $U'$
of $P$ to
a parabolic subgroup and a dominant character $\chi$ of $\fq=\Lie(Q)$.
By lemma \ref{multiple}, there is a positive 
integer $b$ such that $b\chi$ corresponds to
a character $\wt\Xi$ of $Q$. Then, by lemma \ref{samedegree} 
and because of the Ramanathan (semi)stability of $P$, it is
$$
\sum_{i=1}^{t+1} (\lambda_{i+1}-\lambda_i)
\deg E_{\lambda_i}
=\frac{1}{b}\deg P^Q(\wt\Xi) \;(\leq)\; 0.
$$
i.e. $\SP$ is slope-(semi)stable.
\end{proof}

\begin{corollary}
If $X$ is a curve, our notion of (semi)stability
for principal bundles coincides with that of Ramanathan.
\end{corollary}

Let us characterize  (semi)stability in
terms of the Killing form, as announced in the introduction.
An orthogonal sheaf, relative to a scheme $S$, is a pair 
$$
(E_S, E_S\otimes E_S \too \SO_{X\times S})
$$
such that the bilinear form induced on the fibers of $E_S$ 
over closed points
$(x,s)\in X\times S$
where it is locally free, is nondegenerate. 
For instance, if $(E_S,\varphi_S)$
is a $\fgp$-sheaf, the Killing form gives 
an orthogonal structure to $E_S$.

\begin{definition}[Orthogonal filtration]
A filtration $E_{\bullet}\subseteq E$ 
of an orthogonal sheaf is said to be orthogonal
if $E_i^\perp = E_{-i-1}$ for all $i$.
In terms of $E_{\lambda_\bullet}$, if the
integers
$$
\lambda_1<\lambda_2<\cdots<\lambda_t<\lambda_{t+1}
$$ 
can be denoted 
$$
\gamma_{-l}<\gamma_{-l+1}<\cdots<\gamma_{l-1}<\gamma_l
$$
so that
$$
\gamma_{-i}=-\gamma_i, \quad\text{and}\quad E_{\gamma_i}^\perp = E_{\gamma_{-i-1}}
$$
\end{definition}

Observe that an orthogonal filtration is necessarily balanced
and saturated. These filtrations were introduced in our
former article \cite{G-S1} in order to define the (semi)stability
of an orthogonal sheaf as the condition of 
admitting no orthogonal filtration
of negative (nonpositive) Hilbert polynomial.

\begin{corollary}
\label{selforth}
Let $\SP=(P,E,\psi)$ be a principal $G$-sheaf, or just
let $(E,\varphi)$ be a $\fgp$-sheaf.  
An algebra
filtration of $E$ is balanced and saturated if and only if it is orthogonal.
Therefore, $\SP$ is (semi)stable in the sense of definition
\ref{stabilitylie} if and only if it is so in the sense of
definition \ref{stab1}.
\end{corollary}

\begin{proof}
We have seen that a balanced algebra filtration 
of $\fgp$-sheaves is
induced from a filtration of Lie algebras as in (\ref{rootfilt}).
On the other hand, for a semisimple Lie algebra we have
$$
(\fgp^{\alpha})^\perp = \bigoplus_{\beta\neq -\alpha} \fgp^{\beta}
$$
for $\alpha$, $\beta\in R(\fh)\cup \{0\}$.
The first statement follows easily from these two facts.
The second follows from the first and from the fact
that it is enough to consider saturated filtrations.
\end{proof}

\section{Comparison with Gieseker-Maruyama moduli space}
\label{secexample}

The Gieseker-Maruyama moduli space is another natural
compactification of the moduli space of principal
$\GLR$-bundles. In this section we compare this with
the moduli space of semistable 
principal $\GLR$-sheaves.
We give two examples. In the first one, we show that
our moduli space does
not coincide, in general, with the Gieseker-Maruyama 
moduli space of 
torsion free sheaves.
In the second example, we construct examples showing that, 
for principal $\GLR$-bundles, our notion of 
(semi)stability does not coincide, in general, with the 
Gieseker-Maruyama (semi)stability of the associated 
vector bundle (but recall that the slope-(semi)stability
notions do coincide).

\subsection*{Example 1}

Let  $X=\PP^2$ and $G=\gl(2)$. 
The Gieseker-Maruyama moduli space of semistable torsion
free sheaves with rank 2, $c_1=1$ and $c_2=2$ is smooth
of dimension 4.
We are going to show that the moduli space of principal
$\gl(2)$-sheaves with the corresponding numerical invariants 
has a component of dimension at least 16, hence the
two moduli spaces are different.

Let $p\in \PP^2$ be a point. 
Since $\Ext^1(\SO_{\PP^2}(1)\otimes I_p,\SO_{\PP^2})=\CC$,
there is a unique extension up to isomorphism
$$
0 \too \SO_{\PP^2} \too F \too \SO_{\PP^2}(1)\otimes I_p \too 0
$$
It is easy to show that $F$ is a slope-stable vector bundle,
hence the associated principal $\gl(2)$-bundle is also slope-stable.
Let $E$ be the vector bundle associated to the adjoint
representation on $\fgp=\sltwo$.
$$
F^\vee \otimes F = {\rm ad}(F) \oplus \SO_{\PP^2} = 
E \oplus \SO_{\PP^2}
$$

The vector bundle $E$ has rank 3, $c_1=0$ and $c_2=7$, and
furthermore it is Mumford stable. To show this, note that since
$E$ has rank 3 and zero degree, if it is not Mumford stable then  
either
it has a subline bundle of nonnegative degree (hence a nonvanishing
section) or it has a subsheaf of rank two of nonnegative degree.
In this second case, it will have a rank one quotient of nonpositive
degree. Taking the dual, this produces a subline bundle of nonnegative
degree of $E^\vee\isom E^{}$. In both cases, we conclude that
if $E$ is not Mumford stable then it has a nonvanishing section.
Now let $\xi$ be a section of 
$$
H^0(F^\vee \otimes F) = H^0(E) \oplus H^0(\SO_{\PP^2})
$$
Since $F$ is slope stable, it is simple, then $\xi$ is a 
scalar multiple of identity, hence $\xi\in H^0(\SO_{\PP^2})$,
the second summand. This shows that $E$ has no sections,
hence it is Mumford stable.

Let $\Quot(E,4)$ be the Hilbert scheme of quotients
$q:E\to T$ where $T$ is a torsion sheaf of length 4 supported
on a zero-dimensional scheme. For each $q$ define 
$E_q$ to be the kernel
$$
0 \too E_q \stackrel{i}{\too} E \stackrel{q}{\too} T \too 0
$$
This torsion free sheaf inherits an $\sltwo$-sheaf structure
$$
\varphi_q^{}: E_q^{}\otimes E_q^{}
\stackrel{i\otimes i}{\too}
E^{}\otimes E^{} \stackrel{\varphi}{\too} E^{} 
\isom E_q^{\vee\vee}
$$

If $q$ and $q'$ are two quotients corresponding to different
points, then $E_q$ and $E_{q'}$ are not isomorphic.
Indeed, if $\psi$ is an isomorphism between them,
then there is a commutative diagram
$$
\xymatrix{
{0} \ar[r] & {E_q} \ar[r]^{i} \ar[d]_{\psi} & 
{E} \ar[r]^{q} \ar[d]_{\eta} &
{T} \ar[r] \ar[d]^{\xi} & 0 \\
{0} \ar[r] & {E_{q'}} \ar[r]^{i'}  & 
{E} \ar[r]^{q'}  &
{T} \ar[r]  & 0
}
$$
where $\eta$ is induced from $\psi^{\vee\vee}$
and the isomorphisms ${i'}^{\vee\vee}$ and
${i}^{\vee\vee}$. Since $E$ is Mumford stable,
it is simple, and then $\eta=\lambda \id$, a nonzero 
multiple of identity. 
Then the following diagram is commutative
$$
\xymatrix{
{E} \ar@{->>}[r]^{q} \ar@{=}[d] & 
{T} \ar[d]_{\isom}^{\frac{1}{\lambda}\xi} \\
{E} \ar@{->>}[r]^{q'}  & 
{T} 
}
$$
and this implies that $q$ and $q'$ correspond to the
same point in  $\Quot(E,4)$.

The subscheme $\Quot^0(E,4)$ corresponding to quotients supported in 4
distinct points is smooth of dimension 16, and this construction
provides a family of $\sltwo$-sheaves parametrized by $\Quot^0(E,4)$
(with different points giving nonisomorphic $\sltwo$-sheaves). 
To construct a family of principal $\gl(2)$-sheaves we have 
to consider reductions of structure group using the homomorphisms
$$
\gl(2)\too 
\gl(2)/(\ZZ/2\ZZ)=\pgl(2)\times \CC^* \too 
\pgl(2) \too \Aut(\sltwo)
$$
Let
$P^{\Aut(\sltwo)}_q$ be the principal $\Aut(\sltwo)$-bundle on
$U_{E_q}$ associated to $(E_q,\varphi_q)$. 
Since $\sltwo$ has no
outer automorphisms, $\pgl(2)=\Aut(\sltwo)$, and then this
is a principal $\pgl(2)$-bundle.
Now we have to consider a reduction to a principal 
$\pgl(2)\times \CC^*=\gl(2)/(\ZZ/2\ZZ)$-bundle with numerical invariant
equal to 1, i.e. we have to give a line bundle on 
$U_{E_q}$ with degree 1, but
since $U_{E_q}$ is a big open set, 
$\Pic(U_{E_q})=\Pic(\PP^2)$, hence there is a unique
such line bundle: the restriction of $\SO_{\PP^2}(1)$.
Finally, the reductions of a principal 
$\gl(2)/(\ZZ/2\ZZ)$-bundle to $\gl(2)$ are 
parametrized by $\check H^1_\et(U_{E_q},\ZZ/2\ZZ)$.
From theorem \ref{thmcomparison} 
and the proof of lemma \ref{hi}, this is 
isomorphic to the singular cohomology group 
$H^1(\PP^2;\ZZ/2\ZZ)$, and this 
cohomology group is 
trivial because $\PP^2$ is simply connected. Then
there is a unique reduction to a principal 
$\gl(2)$-bundle on $U_{E_q}$. Hence each
$\sltwo$-sheaf $(E_q,\varphi_q)$ produces a unique
principal $\gl(2)$-sheaf, and this provides a
family of principal $\gl(2)$-sheaves parametrized
by a scheme of dimension 16,
with different points giving nonisomorphic objects, 
hence there is 
component of the 
moduli space of principal $\gl(2)$-sheaves of dimension
at least 16.

\subsection*{Example 2}

Let $\pi:X=\wt \PP^2\to \PP^2$ be the blow up of $\PP^2$ at
one closed point. Let $D$ be the exceptional divisor, and $R$ the
divisor class of the strict transform of a line through
the blown up point, i.e. $\pi^* \SO_{\PP^2}(1)=\SO_X(D+R)$.
Hence $K_X=\SO_X(-2D-3R)$. Let $p\in X$ be a closed 
point outside the exceptional divisor. Then
\begin{eqnarray}
\label{pos1}
H^0(\SO_X(aD+bf))\neq 0& \Leftrightarrow & 
a\geq 0 \;\text{and}\; b\geq 0 \\
\label{pos2}
H^0(\SO_X(aD+bf)\otimes I_p)\neq 0& \Leftrightarrow & 
a\geq 0 \;\text{and}\; b> 0 
\end{eqnarray}

Let $c$, $s\in \ZZ$, $L=\SO_X(sR)$, $M=\SO_X(-cD+(c+s)R)$.
The local to global spectral sequence for $\Ext$ gives
an exact sequence
$$
\Ext^1(M\otimes I_p,L) \stackrel{\alpha}{\too}
H^0(\SE xt^1(M\otimes I_p,L))=\CC
\too
H^2(M^\vee\otimes L)=0
$$
The last group is zero by Serre duality and (\ref{pos1}),
hence $\alpha$ is surjective.
The second group is $\CC$ because 
$\SE xt^1(M\otimes I_p,L)\isom\CC_p$, the skyscraper sheaf at $p$. Let 
$\eta$ be an element in the first group with $\alpha(\eta)\neq 0$,
so that the extension corresponding to $\eta$
\begin{equation}
\label{seql}
0 \too L \too F \too M\otimes I_p \too 0
\end{equation}
is locally free. 
Fix the ample line bundle $\SO_X(D+2R)$.
All degrees, stability, etc... will be with respect
to this line bundle.

\begin{lemma}
\label{mumford}
The vector bundle $F$ is Mumford strictly semistable
with slope $\mu(F)=s$, and
the only subsheaf $L'$ with $\mu(L')=\mu(F)$ is
$L$.
\end{lemma}

\begin{proof}
A calculation shows $\mu(L)=\mu(F)=s$. We can assume 
that $L'$ is a line bundle
which does not factor 
through $L$. 
Then the composition
$L' \to F \to M\otimes I_p$ is nonzero, hence
$$
H^0(L'{}^\vee\otimes M \otimes I_p)\neq 0
$$
Denote  $L'=\SO_X(aD+bR)$. 
Using (\ref{pos2}) we obtain $a\leq -c$
and $b<c+s$, and then $\mu(L')<\mu(F)$.
\end{proof}

\begin{lemma}
\label{gmstab}
The vector bundle $F$ is Gieseker-Maruyama (semi)stable
if and only if
$$
3c^2+(2s-1)c+2 \;(\leq)\; 0
$$
\end{lemma}

\begin{proof}
By lemma \ref{mumford}, it is enough to check the
subbundle $L$. A calculation shows that the
following polynomial is constant
$$
P_L(m)-\frac{P_E(m)}{2}=\frac{3c^2+(2s-1)c+2}{2}
$$
and the result follows.
\end{proof}

\begin{lemma}
\label{pbstab}
The principal ${\rm GL}(2)$-bundle associated to 
the vector bundle $E$ is (semi)stable in the sense of definition
\ref{stab1}
if and only if
$$
3 \;(\leq)\; c
$$
\end{lemma}

\begin{proof}
By lemma \ref{bijection}, all orthogonal filtrations come 
from reductions to a parabolic subgroup on a big open set $U'$.
Since $F$ is a rank 2 vector bundle, such a reduction can
be seen as an extension
\begin{equation}
\label{seqf1}
0 \too F_1 \too F \too F_2 \too 0
\end{equation}
where $F_1$ is a line bundle and $F_2$ is a rank one 
torsion free
sheaf. Indeed, this gives a reduction to a maximal
parabolic subgroup on the big open set $U'$ where $F_2$ is
locally free.

Let $E$ be the vector bundle associated to the adjoint
representation on $\fgp=\sltwo$
$$
F^\vee \otimes F = {\rm ad}(F) \oplus \SO_{X} = 
E \oplus \SO_{X}
$$
i.e. the vector bundle $E$ is the sheaf of traceless
homomorphisms $\SH om(F,F)^0$
$$
0 \too \SH om(F,F)^0 \too \SH om(F,F) \stackrel{\tr}{\too} \SO_X \too 0
$$
Since the parabolic subgroup is maximal, there is a 
unique dominant character (up to scalar), and
by lemma \ref{bijection}, this reduction (and character)
gives a saturated
balanced algebra filtration $E_\bullet$ 
(we take double dual $(\,\cdot\,)^{\vee\vee}$
in order to obtain a saturated filtration)
$$
\begin{array}{ccccc}
(F_2^\vee\otimes F_1^{})^{\vee\vee} &  \subsetneq  &
(\SH om_{F_1}(F,F)^0)^{\vee\vee}  & \subsetneq &
\SH om(F,F)  \\
\Arrowvert & & \Arrowvert & & \Arrowvert \\
{E_{-1}} &  & {E_0} & & {E_1=E}
\end{array}
$$
where $\SH om_{F_1}(F,F)^0$ denotes the sheaf of traceless 
homomorphisms preserving $F_1$, i.e. it is the kernel
of the homomorphism $\alpha$
\begin{equation}
\label{respecting}
0 \too \SH om_{F_1}(F,F)^0 \too \SH om(F,F)^0
\stackrel{\alpha}{\too} F_1^\vee\otimes F_2^{}
\end{equation}
The Hilbert polynomial of the filtration $E_\bullet$ is
$$
P_{E_\bullet} = (3P_{E_{-1}}-P_E) + (3P_{E_0}-2P_E)
$$
It has degree at most 1, and the coefficient of the
term of degree 1 can be obtained substituting the Hilbert
polynomials by degrees in the previous expression.
We calculate
\begin{eqnarray*}
\deg(E_{-1})& = &\deg(F_2^\vee\otimes F_1^{}) = \deg F_1 - \deg F_2\\
\deg(E_{0})&= &\deg(\SH om_{F_1}(F,F)^0) =\\
& =& \deg(\SH om(F,F)^0) -\deg(F_1^\vee\otimes F_2^{})
= \deg F_1 - \deg F_2 \\
\deg(E) &=& 0 \\
\end{eqnarray*}
The degree of $E_0$ can be calculated from the exact 
sequence (\ref{respecting}), because the homomorphism 
$\alpha$ is surjective where $F_2$ is locally free, and 
this is a big open set.
Then the coefficient of degree 1 in the Hilbert polynomial
$P_{E_\bullet}$ of the filtration $E_\bullet$ is
$$
3(\deg(E_{-1}) + \deg(E_{0}) - \deg(E))=6 (\deg F_1 - \deg F_2)
= 12( \deg F_1 - \frac{\deg F}{2})
$$
By lemma \ref{mumford}, $\deg F_1 - (\deg F)/2 \leq 0$, and to check
the (semi)stability of the principal bundle we can assume
that the sequence (\ref{seqf1}) is (\ref{seql}):
$$
0 \too L \too F \too M\otimes I_p \too 0
$$
Taking the dual of this sequence, we obtain the short
exact sequence
$$
0 \too M^\vee \too F^\vee \too L^\vee\otimes I_p \too 0
$$
It is easy to check that the composition
$$
\SH om(F,F)^0 \inj F^\vee\otimes F \surj 
L^\vee\otimes I_p\otimes M\otimes I_p
$$
is surjective, hence $\SH om_L(F,F)^0$ fits in an
exact sequence
$$
0 \too \SH om_{L}(F,F)^0 \too \SH om(F,F)^0
\too L^\vee\otimes I_p\otimes M\otimes I_p 
\; = \; L^\vee\otimes  M\otimes (I_Z\oplus \CC_p) \too 0
$$
where $Z$ is the ``fat point'' supported at $p$, i.e.
$I^{}_Z=I^2_p$. Then 
$$
0 \too (\SH om_{L}(F,F)^0)^{\vee\vee} \too \SH om(F,F)^0
\too  L^\vee\otimes  M\otimes I_Z \too 0
$$
On the other hand, $(F_2^\vee\otimes F_1)^{\vee\vee}=M^\vee\otimes L$.
A calculation shows
$$
P_{E_\bullet}(m)=3 (-c+3)
$$
and the result follows.
\end{proof}

Lemmas \ref{gmstab} and \ref{pbstab} show that for a 
principal $\gl(2)$-bundle, its (semi)stability 
in the sense of definition \ref{stab1} does not
coincide in general with the 
Gieseker-Maruyama (semi)stability of the
associated rank 2 vector bundle.
The following table gives 
the stability of $F$ for concrete values
of the parameters $c$ and $s$,
showing this fact.

\medskip
\begin{center}
\begin{tabular}{|c|c|c|}
\hline
$(c,s)$ & Vector bundle  & Principal bundle \\
\hline
(-1,4) & unstable        &unstable\\
(-1,3) & semistable      &unstable\\
(-1,2) & stable          & unstable\\
(3,-4) & unstable        & semistable\\
(3,-5) & stable          & semistable\\
(4,-5) & unstable        & stable\\
(4,-6) & stable          & stable\\
\hline
\end{tabular}
\end{center}


\begin{thebibliography}{EMG}

\bibitem[A-B]{A-B}{B. Anchouche and I. Biswas, }
\textit{Einstein--Hermitian connections on polystable principal
bundles over a compact K\"ahler manifold, }
Amer. Jour. Math. \textbf{123} (2001), 207--228.



\bibitem[B-T]{B-T}{A. Borel, J. Tits, } 
\textit{Groupes r\'eductifs, } 
Inst. Hautes \'Etudes Sci. Publ. Math. \textbf{27} (1965), 55--150. 

\bibitem[FGA]{FGA}{A. Grothendieck, }
\textit{Technique de descente et th\'eor\`emes d'existence en
g\'eom\'etrie alg\'ebrique, }
(Bourbaki expos\'es 190, 195, 212, 221, 232 and 236). Also in:
\textit{Fondements de la g\'eom\'etrie alg\'ebrique, }
Secr\'etariat Mathematique Paris (1962).

\bibitem[F-M]{F-M} R. Friedman and J.W. Morgan : Holomorphic Principal 
Bundles over Elliptic Curves II: The Parabolic Construction. 
J. Differential Geom. \textbf{56} (2000), 301--379. 

\bibitem[Gi]{Gi}{D. Gieseker, }
\textit{On the moduli of vector bundles on an algebraic surface, }
Ann. Math., \textbf{106} (1977), 45--60.

\bibitem[G-S]{G-S}{T. G\'omez and I. Sols, }
\textit{Stability of conic bundles, } 
Internat. J. Math., \textbf{11} (2000), 1027--1055.


\bibitem[G-S1]{G-S1}{T. G\'omez and I. Sols, }
\textit{Stable tensors and moduli space of orthogonal sheaves, }
Preprint 2001, math.AG/0103150.

\bibitem[G-S2]{G-S2}{T. G\'omez and I. Sols, }
\textit{Projective moduli space of semistable principal
sheaves for a reductive group, }
Le Matematiche \textbf{15} (2000), 437--446,
conference in honor of Silvio Greco (April 2001)

\bibitem[J]{J}{N. Jacobson, }
Lie algebras, Dover, 1979.

\bibitem[Ha]{Ha}{R. Hartshorne, }
Algebraic Geometry, Grad. Texts in Math. 52, Springer Verlag, 1977.

\bibitem[Ha2]{Ha2}{R. Hartshorne, }
\textit{Stable reflexive sheaves, }
Math. Ann. \textbf{254} (1980), 121--176.

\bibitem[Hum]{Hum}{J. E. Humphreys, }
Introduction to Lie algebras
and representation theory, Grad. Texts in Math. 9,
Springer Verlag, New York-Berlin, 1972.


\bibitem[H-L]{H-L}{D. Huybrechts and M. Lehn, }
The geometry of moduli spaces of sheaves,
Aspects of Mathematics E31, Vieweg, Braunschweig/Wiesbaden 1997.


\bibitem[Hy]{Hy}{D. Hyeon, }
\textit{Principal bundles over a projective scheme, }
Trans. Amer. Math. Soc. \textbf{354} (2002), 1899--1908.


\bibitem[Ma]{Ma}{M. Maruyama, }
\textit{Moduli of stable sheaves, I and II.} 
J. Math. Kyoto Univ. \textbf{17} (1977), 91--126. 
\textbf{18} (1978), 557--614. 


\bibitem[Mu1]{Mu1}{D. Mumford, }
Geometric invariant theory. Ergebnisse der Mathematik und ihrer
Grenzgebiete, Neue Folge, Band 34. Springer-Verlag, Berlin-New York,
1965. 

\bibitem[Mu2]{Mu2}{D. Mumford, }
Abelian varieties. Oxford University Press, Bombay, 1970.

\bibitem[N-S]{N-S}
{M.S. Narasimhan and C.S. Seshadri, }
\textit{Stable and unitary vector bundles
on a compact Riemann surface.}
Ann. of Math. (2) \textbf{82} (1965), 540--567.

\bibitem[Ra1]{Ra1}{A. Ramanathan, }
\textit{
Stable principal bundles on a compact Riemann surface.}
Math. Ann. \textbf{213} (1975), 129--152. 

\bibitem[Ra2]{Ra2}{A. Ramanathan, }
\textit{Moduli for principal bundles. Algebraic geometry 
(Proc. Summer Meeting, Univ. Copenhagen, Copenhagen, 1978), } 
pp. 527--533, Lecture
Notes in Math., \textbf{732}, Springer, Berlin, 1979. 


\bibitem[Ra3]{Ra3}{A. Ramanathan, }
\textit{Moduli for principal bundles over algebraic curves: 
I and II, }
Proc. Indian Acad. Sci. (Math. Sci.), \textbf{106} (1996), 301--328,
and 421--449. 

\bibitem[Ri]{Ri}{R. Richardson, }
\textit{Compact real forms of a complex semi-simple Lie
algebra, } 
J. Differential Geometry \textbf{2} (1968), 411--419. 


\bibitem[Sch]{Sch}{A. Schmitt, }
\textit{Singular principal bundles over higher-dimensional manifolds and their
moduli spaces, }
Internat. Math. Res. Notices
\textbf{23} (2002), 1183--1210.
 

\bibitem[Se1]{Se1}{J. P. Serre, }
\textit{Espaces fibr\'es alg\'ebriques, }
S\'eminaire C. Chevalley, ENS (1958)

\bibitem[Se2]{Se2}{J. P. Serre, }
\textit{Lie algebras and Lie groups, 1964 Lectures given
at Harvard University}, (1965), W. A. Benjamin Inc.

\bibitem[Sesh]{Sesh}{C.S. Seshadri, }
\textit{Space of unitary vector bundles on a compact Riemann surface.} 
Ann. of Math. (2) \textbf{85} (1967), 303--336.  

\bibitem[SGA1]{SGA1}{A. Grothendieck and M. Raynaud, }
\textit{Rev\^etements \'etales et groupe fondemental, }
Lecture Notes in Math. \textbf{224}, Springer-Verlag 
(1971).

\bibitem[SGA4]{SGA4}{M. Artin, A. Grothendieck, J.L. Verdier, }
\textit{Th\'eorie des topos et cohomologie \'etale des sch\'emas, }
Lecture Notes in Math. \textbf{269}, \textbf{270}, \textbf{305}, 
Springer-Verlag (1972).

\bibitem[Si]{Si}{C. Simpson, }
\textit{Moduli of representations of the fundamental group of a smooth
projective variety I, }
Publ. Math. I.H.E.S. \textbf{79} (1994), 47--129.


\bibitem[Sp]{Sp}{E.H. Spanier, }
Algebraic topology,
Corrected reprint. 
Springer-Verlag, New York-Berlin, 1981. 


\end{thebibliography}
\end{document}